\documentclass[11pt,fullpage]{article}

\usepackage{amsmath}
\usepackage{amssymb}
\usepackage{amsfonts}
\usepackage{xr}

\numberwithin{equation}{section}
\newtheorem{thm}{Theorem}[section] 
\newtheorem{prp}[thm]{Proposition}
\newtheorem{lmm}[thm]{Lemma}   
\newtheorem{crl}[thm]{Corollary} 
\newtheorem{dfn}[thm]{Definition} 

\def\e_ref#1{(\ref{#1})}
\def\under#1{\underline{#1}}
\def\ov#1{\overline{#1}}
\def\lra{\longrightarrow}
\def\Lra{\Longrightarrow}
\def\Llra{\Longleftrightarrow}

\def\lan{\langle}
\def\ran{\rangle}

\def\al{\alpha}

\def\de{\delta}
\def\ep{\epsilon}
\def\ga{\gamma}
\def\io{\iota}
\def\ka{\kappa}
\def\la{\lambda}

\def\ve{\varepsilon}
\def\ups{\upsilon}
\def\vp{\varpi}
\def\vt{\vartheta}

\def\lap{\Delta}

\def\Ga{\Gamma}
\def\La{\Lambda}
\def\Om{\Omega}
\def\Si{\Sigma}

\def\rk{\text{rk}}
\def\ev{\text{ev}}

\def\RT{\text{RT}}

\def\P{\Bbb{P}^n}
\def\PP{\Bbb{P}^2}

\def\i{\infty}
\def\eset{\emptyset}

\begin{document}

\title{Enumeration of Genus-Three Plane Curves\\ 
with a Fixed Complex Structure}

\author{Aleksey Zinger
\thanks{Partially supported by NSF grant DMS-9803166}}

\maketitle

\begin{abstract}
\noindent
We give a practical formula for counting
irreducible nodal genus-three plane curves that 
a fixed generic complex structure on the normalization.
As an intermediate step, we enumerate rational
plane curves that have a $(3,4)$-cusp.
\end{abstract}

\tableofcontents

\section{\bf Introduction} 

\subsection{Background and Results}
\label{back}

\noindent
Let $(\Si,j_{\Si})$ be a nonsingular Riemann surface of genus-$3$
and let $d$ be a positive integer. 
Denote by ${\cal H}_{\Si,d}(\PP)$ the set of simple holomorphic maps 
from $\Si$ to $\PP$ of degree~$d$. 
Let $\mu\!=\!(p_1,\ldots,p_{3d-4})$
be a tuple of points in~$\PP$ in general position.
Then the~set
\begin{equation}
\label{intro_1}
{\cal H}_{\Si,d}(\mu)=
\big\{(y_1,\ldots,y_{3d-4};u)\!: u\!\in\!{\cal H}_{\Si,d}(\PP);~
y_l\!\in\!\Si,~u(y_l)\!=\!\mu_l~\forall l\!=\!1,\ldots,3d\!-\!4\big\}
\end{equation}
is finite. 
For a dense open subset of complex structures $j_{\Si}$ on~$\Si$,
the cardinality of this set is the same number~$n_{3,d}$.
More intrinsically, $n_{3,d}$ is
the number of genus-3 degree-$d$ plane curves
that pass through $3d\!-\!4$ points  in general position
and have a pre-specified generic complex structure on the normalization.\\

\noindent
Enumerative numbers, such as $n_{3,d}$, have been of interest
in algebraic geometry for a long time.
The low-genus numbers, $n_d\!\equiv\!n_{0,d}$, $n_{1,d}$,
and $n_{2,d}$ are computed in~\cite{KM}, \cite{RT},
\cite{I}, \cite{P1}, \cite{Z2}, and~\cite{KQR} 
with completion in~\cite{Z1}.
In this paper, we apply the machinery developed in~\cite{Z2}
and~\cite{Z3} to compute the numbers~$n_{3,d}$.\\

\noindent
It is shown in \cite{RT} that 
$$n_d=\RT_{0,d}(\mu_1,\mu_2,\mu_3;\mu_4,\ldots,\mu_{3d-1}),$$
where $\RT_{0,d}(\cdot;\cdot)$ denotes the symplectic invariant of $\PP$
defined as in~\cite{RT}.
In \cite{I}, the~difference 
$$\RT_{1,d}(\mu_1;\mu_2,\ldots,\mu_{3d-1})-2n_{1,d}$$
is shown to be a certain multiple of $n_d$.
Extending the general approach of~\cite{I}, in~\cite{Z2},
the difference
$$\RT_{2,d}(;\mu_1,\ldots,\mu_{3d-2})-2n_{2,d}$$
is expressed in terms of the numbers $n_{d'}$ with $d'\!\le\! d$.
Due to the two composition laws of~\cite{RT},
the symplectic invariants $\RT_{g,d}(\cdot,\cdot)$ of~$\PP$
are easily computable.
Thus, comparing enumerative invariants of~$\PP$ to
the symplectic ones as above is sufficient 
for computing the enumerative invariants.
In this paper, we prove

\begin{thm}
\label{g3n2_thm}
If $d$ is a positive integer and $\mu$ is a tuple of $3d\!-\!4$ points
in general position in $\PP$,
\begin{gather*}
n_{3,d}=RT_{3,d}(\cdot;\mu)-CR_3(\mu),
\qquad\hbox{where}\\
\begin{split}
&\frac{1}{12}CR_3(\mu)\!=\!\big\lan 413 a^2c_1^2({\cal L}^*)
+210ac_1^3({\cal L}^*)+44c_1^4({\cal L}^*),[\bar{\cal V}_1(\mu)]\big\ran\\
&\quad~-\!\big\lan 217a^2+
  84a\big(c_1({\cal L}_1^*)\!+\!c_1({\cal L}_2^*)\big)
+16\big(c_1^2({\cal L}_1^*)\!+\!c_1^2({\cal L}_2^*)\big)
+10c_1({\cal L}_1^*)c_1({\cal L}_2^*),[\bar{\cal V}_2(\mu)]\big\ran
+18|\bar{\cal V}_3(\mu)|.
\end{split}
\end{gather*}
\end{thm}

\begin{center}
\begin{tabular}{||c|c|c|c|c|c|c||}
\hline\hline
$d$&  2&3&4&5&6&7\\
\hline
$n_d$& 0&0&14,280& 9,469,152& 6,573,686,112&6,289,178,278,656
\\ 
\hline\hline
\end{tabular}
\end{center}

\noindent
The notation used in Theorem~\ref{g3n2_thm} is the same as in~\cite{Z2}; 
see Subsection~\ref{notation_sec} for more details.
For now, it is sufficient to say that the intersection numbers
in Theorem~\ref{g3n2_thm} can be expressed in terms of intersection
numbers of tautological classes in the space of stable rational
maps into~$\PP$.
The latter are shown to be computable in~\cite{P2};
see also Subsection~\ref{enum1-chern_class} in~\cite{Z2}.
Thus, Theorem~\ref{g3n2_thm} gives a practical formula for 
computing the numbers~$n_{3,d}$.
Our number~$n_{3,4}$ agrees with that of~\cite{AF}.\\

\noindent
Along the way, we also enumerate rational curves with 
certain singularities.
In particular, let ${\cal S}_{1;2}(\mu)$ denote
the set of plane rational degree-$d$ curves 
that have a $(3,4)$-cusp and pass through the $3d\!-\!4$ points
$\mu_1,\ldots,\mu_{3d-4}$.
Corollary~\ref{n2cusps2c} expresses the number $|{\cal S}_{1;2}(\mu)|$ 
in terms of intersections of tautological classes
in the space of rational maps.
See also Lemmas~\ref{n2_2cusps}, \ref{n2tacnodes}, and~\ref{n2cusps}.
The degree-four numbers of Lemma~\ref{n2_2cusps} and
Corollary~\ref{n2cusps2c} agree with the numbers computed
by P.~Aluffi from a formula in~\cite{AF2}.\\

\noindent
The argument of this paper can be modified to enumerate 
genus-3 plane curves with a fixed non-generic smooth complex
structure on the normalization.
In particular, suppose $(\Si,j)$
is not hyperelleptic and has $n$ ``hyperflexes'' in the sense of~\cite{AF},
i.e.~exactly $n$ Weierstrass points with gap values $(1,4)$
and $24\!-\!2n$ Weierstrass points with gap values $(1,3)$;
see~\cite[p273]{GH}.
Then the number of genus-3 degree-$d$ plane curves
passing through $3d\!-\!4$ points and with normalization~$(\Si,j)$
is $\big(n_{3,d}\!-\!2n|S_{1;2}(\mu)|\big)/\hbox{Aut}(\Si,j)$;
see the remarks following Corollary~\ref{CR3_str}.
In the $d\!=\!4$ case, \cite{AF} obtain the same correction.\\

\noindent
The author thanks P.~Aluffi, T.~Mrowka, and R.~Vakil
for helpful conversations.
In particular, the question of enumerating genus-$3$ curves with 
a fixed non-generic complex structure was posed to the author by P.~Aluffi,
who also pointed out that the conclusion of the original limiting argument
for obtaining counts of such curves directly from the numbers~$n_{3,d}$
could not have been be valid.

\subsection{Summary}
\label{summary}

\noindent
We now outline the proof of Theorem~\ref{g3n2_thm}.
If $\nu\!\in\!\Ga(\Si\!\times\!\PP;
            \La^{0,1}\pi_{\Si}^*T^*\Si\otimes\pi_{\PP}^*T\PP)$, 
let ${\cal M}_{\Si,\nu,d}$ denote the set of all smooth maps~$u$ 
from $\Si$ to $\PP$ of degree $d$ such that
$\bar{\partial}u|_z\!=\!\nu|_{(z,u(z))}$ for all $z\!\in\!\Si$.
If $\mu$ is as above and $N\!=\!3d\!-\!4$, put
$${\cal M}_{\Si,\nu,d}(\mu)=
\big\{(y_1,\ldots,y_N;u)\!: u\!\in\!{\cal M}_{\Si,\nu,d};~
y_l\!\in\!\Si,~ u(y_l)\!=\!\mu_l ~\forall l\!=\!1,\ldots,N\big\}.$$
For a generic $\nu$, ${\cal M}_{\Si,\nu,d}$ is a smooth
finite-dimensional oriented manifold, and 
${\cal M}_{\Si,\nu,d}(\mu)$
is a zero-dimensional finite submanifold of 
${\cal M}_{\Si,\nu,d}\!\times\!\Si^N$, 
whose cardinality (with sign) is the symplectic invariant~$\RT_{3,d}(;\mu)$.
This number depends only on the degree~$d$; see~\cite{RT}.\\

\noindent
If $\|\nu_i\|_{C^0}\!\lra\! 0$ and 
   $(\under{y}_i;u_i)\!\in\! {\cal M}_{\Si,\nu_i,d}(\mu)$, 
then a subsequence of $\{(\under{y}_i;u_i)\}_{i=1}^{\i}$ must converge 
in the Gromov topology to one of the following:\\
(1) an element of ${\cal H}_{\Si,d}(\mu)$;\\
(2) $(\Si_b,\under{y},u_b)$, where $\Si_b$ is a bubble tree of 
$S^2$'s attached to $\Si$ with marked points $y_1,\ldots,y_N$,
and \hbox{$u_b\!:\Si_b\!\lra\!\PP$} is a holomorphic map such that
$u_b(y_l)\!=\!\mu_l$ for $l\!=\!1,\ldots,N$, and\\
(2a) $u_b|\Si$ is simple and the tree contains at least one $S^2$;\\
(2b) $u_b|\Si$ is multiply-covered;\\
(2c) $u_b|\Si$ is constant and the tree contains at least one $S^2$.\\
By an argument similar to the proof of
Proposition~\ref{enum1-empty_spaces} in~\cite{Z2},
the cases (2a) and (2b) cannot occur.
As in~\cite{Z2}, our approach will be to take $t\!>\!0$ very small and 
to determine the number $CR_3(\mu)$ of elements of 
${\cal M}_{\Si,t\nu,d}(\mu)$ that lie near the maps of type~(2c).
The rest of the elements of ${\cal M}_{\Si,t\nu,d}(\mu)$
must lie near the space ${\cal H}_{\Si,d}(\mu)$.
By Proposition~\ref{gl-holom_prp} in~\cite{Z3} and 
Corollary~\ref{enum1-reg_crl} in~\cite{Z2},
there is a one-to-one correspondence between
the elements of ${\cal H}_{\Si,d}(\mu)$
and the nearby elements of ${\cal M}_{\Si,t\nu,d}(\mu)$,
at least if $d\!\ge\!5$.
In fact, a standard argument shows that
Corollary~\ref{enum1-reg_crl} remains valid
if $d\!=\!4$ and $(\Si,j)$ is generic.
If $d\!=\!1,2,3$, ${\cal H}_{\Si,d}(\mu)\!=\!\eset$
\hbox{by~\cite[p116]{ACGH}}.
Thus, we are able to compute the cardinality of 
${\cal H}_{\Si,d}(\mu)$ by computing the total number~$CR_3(\mu)$
of elements of ${\cal M}_{\Si,t\nu,d}(\mu)$ that lie
near the maps of type~(2c).\\

\noindent
The space of stable maps of type (2c) is stratified by smooth,
usually noncompact, manifolds~${\cal M}_{\cal T}(\mu)$.
The set of elements of ${\cal M}_{\Si,t\nu,d}(\mu)$
that lie near each space~${\cal M}_{\cal T}(\mu)$
corresponds to the zero set of a map between two bundles 
over~${\cal M}_{\cal T}(\mu)$.
By extracting dominant terms from each such map,
the signed cardinality $N({\cal T})$ of the zero set of 
the map can be identified with the signed cardinality of 
the zero set of an affine map
between vector bundles over a closure of~${\cal M}_{\cal T}(\mu)$
or of a certain submanifold of~${\cal M}_{\cal T}(\mu)$.
The argument is nearly the same as in 
Sections~\ref{enum1-analysis_sec} and~\ref{enum1-resolvent_sec} 
of~\cite{Z2}.
It is summarized briefly at the end of Subsection~\ref{str_thm_sec}.
The number $CR_3(\mu)$, the sum of the numbers~$N({\cal T})$, 
can then be expressed as the sum of the cardinalities 
of the zero sets of affine vector bundles over compact manifolds;  
see Corollary~\ref{CR3_str}.
Topological formulas for the six numbers~$n_m^{(k)}(\mu)$
of Corollary~\ref{CR3_str} are obtained in Section~\ref{comp_sect};
see Lemmas~\ref{m1k1}, \ref{m1k2}, \ref{m1k3},
\ref{m2k1}, \ref{m2k2}, and~\ref{m3k1}, respectively.
Using the results of Lemmas~\ref{n2_2cusps}, \ref{n2cusps},
and~\ref{n2cusps2}, 
we obtain the expression for $CR_3(\mu)$ given in Theorem~\ref{g3n2_thm}.\\

\noindent
In Subsection~\ref{top_sec}, we review the topological tools
to be used in Sections~\ref{rational_sect} and~\ref{comp_sect}.
We summarize our notation for spaces of bubble maps and vector bundle 
over them in Subsection~\ref{notation_sec}.
Subsection~\ref{str_global_sec} describes 
the structure of spaces of rational maps and
the behavior of certain bundle sections over them
near the boundary strata.
These descriptions are needed to implement the topological tools
of Subsection~\ref{top_sec} in Sections~\ref{rational_sect}
and~\ref{comp_sect}.
In Subsection~\ref{str_thm_sec}, we describe the number of elements
of ${\cal M}_{\Si,t\nu,d}(\mu)$ near each given strata 
${\cal M}_{\cal T}(\mu)$ of bubble maps of type~(2c)
in terms of the zero sets of affine maps 
between finite-rank vector bundles over relatively 
simple topological spaces.\\

\noindent
In Section~\ref{rational_sect}, we enumerate rational curves
with certain singularities and also deal 
the intersection numbers used in Section~\ref{comp_sect}.
Rational curves with singularities can be identified 
with the zeros of bundle sections over spaces of 
stable rational maps that lie in the main stratum.
We use the topological tools of Subsection~\ref{top_sec}
to determine the contribution from the boundary strata of such
spaces to the euler class of the bundle.
Finally, in Section~\ref{comp_sect}, we derive topological formulas
for the six numbers of Corollary~\ref{CR3_str}
using the same approach as in Section~\ref{rational_sect}.

\section{The Computational Setting}

\subsection{Topology}
\label{top_sec}

\noindent
We begin by describing the topological tools
used in the next two sections.
In particular, we review the notion of 
contribution to the euler class of a vector bundle from 
a (not necessarily closed) subset of the zero set of a section.
We also recall how one can enumerate the zeros of 
an affine map between vector bundles.
These concepts are closely intertwined.
Details can be found in Section~\ref{enum1-top_sect} of~\cite{Z2}.
\\

\noindent
Throughout this paper, all vector bundles 
are assumed to be complex and normed.
If $F\!\lra\!{\cal M}$ is a smooth vector bundle,
closed subset $Y$ of $F$ is {\it small}
if it contains no fiber of~$F$ and is preserved 
under scalar multiplication.
If ${\cal Z}$ is a compact oriented zero-dimensional manifold,
we denote the signed cardinality of~${\cal Z}$
by~$^{\pm}|{\cal Z}|$.

\begin{dfn}
\label{top_dfn1b}
Suppose $F,{\cal O}\!\lra\!{\cal M}$ are smooth vector bundles.\\
(1) If $F\!=\!\bigoplus\limits_{i=1}^{i=k}F_i$,
bundle map $\al\!:F\!\lra\!{\cal O}$ 
is a \under{polynomial of degree~$d_{[k]}$} if
for each $i\!\in\![k]$ there~exists
$$p_i\!\in\!\Ga({\cal M};F_i^{*\otimes d_i}\!\otimes\!{\cal O})
\hbox{~~for~}i\!\in\![k]
\qquad\hbox{s.t.}\qquad
\al(\ups)=\sum_{i=1}^{i=k}p_i\big(\ups_i^{d_i}\big)
\quad\forall \ups=(\ups_i)_{i\in[k]}\in\bigoplus_{i=1}^{i=k}F_i.$$
(2) If $\al\!:F\!\lra\!{\cal O}$ is a polynomial, 
the \under{rank of~$\al$} is the number
$$\rk~\al\equiv\max\{\rk_b\al\!: b\!\in\!{\cal M}\},
\quad\hbox{where}\quad
\rk_b\al=\dim_{\Bbb{C}}\big(\hbox{Im}~\al_b\big).$$
Polynomial $\al\!:F\!\lra\!{\cal O}$ is \under{of constant rank}
if $\rk_b\al=\rk~\al$ for all $b\!\in\!{\cal M}$;
$\al$ is \under{nondegenerate}
if \hbox{$\rk_b\al=rk~F$} for all $b\!\in\!{\cal M}$.\\
(3) If $\Om$ is an open subset of $F$ and 
$\phi\!: \Om\!\lra\!{\cal O}$ is a smooth bundle map, 
bundle map $\al\!:F\!\lra\!{\cal O}$ is
a \under{dominant term of~$\phi$} 
if there exists $\ve\!\in\! C^0(F; \Bbb{R})$ such that
$$\big|\phi(\ups)-\al(\ups)\big|\le
\ve(\ups)|\al(\ups)|
\quad\forall\ups\!\in\!\Om
\quad\hbox{and}\quad
\lim_{\ups\lra0}\ve(\ups)=0.$$
Dominant term $\al\!:F\!\lra\!{\cal O}$ of $\phi$ is 
the \under{resolvent} of $\phi$ if $\al$ is a polynomial
of constant rank.
(4) $\phi\!: \Om\!\lra\!{\cal O}$ is \under{hollow}
if there exist dominant term $\al$ of $\phi$ 
and splittings $F\!=\! F^-\!\oplus F^+$
and ${\cal O}\!=\!{\cal O}^-\!\oplus {\cal O}^+$ 
such that $\al(F^+)\!\subset\!{\cal O}^+$,
$\al^-\!\equiv\!\pi^-\!\circ(\al|F^-)$ is a constant-rank polynomial,
where $\pi^-\!:{\cal O}\!\lra\!{\cal O}^-$ is the projection map,
and \hbox{$(\rk~\al^-\!+\!\frac{1}{2}\dim{\cal M})\!<\!\rk~{\cal O}^-$}.
\end{dfn}

\noindent
The base spaces we work with in the next two sections
are closely related to 
spaces of rational maps into $\PP$ of total degree~$d$
that pass through the $N$ points $\mu_1,\ldots,\mu_N$,
where $N\!=\!3d\!-\!4$ as before.
From the algebraic geometry point of view,
spaces of rational maps are algebraic stacks, 
but with a fairly obscure local structure.
We view these spaces as {\it mostly smooth}, or {\it ms-}, manifolds:
compact oriented topological manifolds
stratified by smooth manifolds, such that the boundary strata
have (real) codimension at least two.
Subsection~\ref{str_global_sec} gives explicit descriptions
of neighborhoods of boundary strata and
of the behavior of certain bundle sections near such strata.
We call the main stratum ${\cal M}$ of ms-manifold $\bar{\cal M}$
the {\it smooth base} of~$\bar{\cal M}$.
Definition~\ref{enum1-euler_dfn1} in~\cite{Z2}
also introduces the natural notions of
{\it ms-maps} between ms-manifolds,
{\it ms-bundles} over ms-manifolds, and 
{\it ms-sections} of ms-bundles.

\begin{dfn}
\label{euler_dfn2}
Let  $\bar{\cal M}\!=\!{\cal M}_n\sqcup\bigsqcup_{i=0}^{n-2}{\cal M}_i
\!=\!{\cal M}\sqcup\bigsqcup_{i=0}^{n-2}{\cal M}_i$ 
be an ms-manifold of dimension~$n$.\\
(1) If ${\cal Z}\!\subset\!{\cal M}_i$ is a smooth oriented submanifold, 
a \under{normal-bundle model for ${\cal Z}$} is a tuple
$(F,Y,\vt)$, where\\
(1a) $F\!\lra\!{\cal Z}$ is a smooth vector bundle
and $Y$ is a small subset of~$F$;\\
(1b) for some $\de\!\in\!C^{\i}({\cal Z};\Bbb{R}^+)$,
$\vt\!:F_{\de}\!-\!(Y\!-\!{\cal Z})\lra\!\bar{\cal M}$ 
is a continuous map such that\\ 
(1b-i) $\vt\!:F_{\de}\!-\!(Y\!-\!{\cal Z})\!\lra\!\bar{\cal M}$ 
is a homeomorphism
onto an open neighborhood of ${\cal Z}$ in ${\cal M}\cup{\cal Z}$;\\
(1b-ii) $\vt|_{\cal Z}$ is the identity map, and
$\vt\!:F_{\de}\!-\!Y\!-\!{\cal Z}\!\lra\!{\cal M}$
is an orientation-preserving diffeomorphism on an open subset 
of~${\cal M}$.\\
(2) A \under{closure} of normal-bundle model $(F,Y,\vt)$ for ${\cal Z}$
is a tuple $(\bar{\cal Z},\tilde{F},\pi)$, where\\
(2a) $\bar{\cal Z}$ is an ms-manifold with 
smooth base ${\cal Z}$;\\
(2b) $\pi\!:\bar{\cal Z}\!\lra\!\bar{\cal M}$ is an ms-map
such that $\pi|_{\cal Z}$ is the identity;\\
(2c) $\tilde{F}\!\lra\!\bar{\cal Z}$ is an ms-bundle such that
$\tilde{F}|_{\cal Z}\!=\!F$.
\end{dfn}

\noindent
We use a normal-bundle model for ${\cal Z}$ to describe
the behavior of bundle sections over $\bar{\cal M}$ near~${\cal Z}$.
In particular, if $\al\!:E\!\lra\!{\cal O}$ is an ms-polynomial,
we call ${\cal Z}$ an {\it $\al$-regular} subset of $\bar{\cal M}$
if for some normal-bundle model $(F,Y,\vt)$ for ${\cal Z}$,
$\vt^*\al$ can be approximated, by a constant-rank polynomial
$p\!:F\oplus E\!\lra\!{\cal O}$;
see Definition~\ref{enum1-euler_dfn3} in~\cite{Z2}.
Polynomial $\al\!:E\!\lra\!{\cal O}$ is {\it regular}
if $\bar{\cal M}$ can be decomposed into finitely many
$\al$-regular subsets.
If $\rk~E\!+\!\frac{1}{2}\dim\bar{\cal M}\!=\!\rk~{\cal O}$,
for a generic $\nu\!\in\!\Ga(\bar{\cal M};{\cal O})$,
the zero set of the polynomial map
$$\psi_{\al,\nu}\!:E\lra{\cal O},\qquad 
\psi_{\al,\nu}(\ups)=\nu_{\ups}+\al(\ups),$$
is a zero-dimensional oriented submanifold of~$E|{\cal M}$.
By Lemma~\ref{enum1-euler_l1l} in~\cite{Z2},
if $\al$ is a regular polynomial,
$\psi_{\al,\nu}^{-1}(0)$ is a finite set for a generic choice of~$\nu$,
and $N(\al)\equiv^{\pm}\!\big|\psi_{\al,\nu}^{-1}(0)\big|$
is independent of such a choice of~$\nu$.\\

\noindent
In Section~\ref{rational_sect}, 
rational curves with pre-specified singularities are identified with 
the zero set of a section $s$ of a vector bundle~$V$
over a smooth manifold ${\cal M}$.
The section~$s$ extends over a compactification~$\bar{\cal M}$ of~${\cal M}$.
Thus, the subset of  ${\cal M}$ we are interested in
can be identified with the euler class $e(V)$ of $V$
minus the $s$-contribution ${\cal C}_{\bar{\cal M}\!-\!{\cal M}}(s)$
to $e(V)$ from $\bar{\cal M}\!-\!{\cal M}$.
In the cases we encounter in Sections~\ref{rational_sect}
and~\ref{comp_sect}, $s^{-1}(0)\cap(\bar{\cal M}\!-\!{\cal M})$
decomposes into disjoint, and usually non-compact, manifolds~${\cal Z}_i$ 
near which the behavior of $s$ can be understood.
Then 
${\cal C}_{\bar{\cal M}\!-\!{\cal M}}(s)\!=\!\sum{\cal C}_{{\cal Z}_i}(s)$,
where ${\cal C}_{{\cal Z}_i}(s)$ is
the {\it $s$-contribution of ${\cal Z}_i$ to $e(V)$}.
This is the signed number of elements of $\{s\!+\!\nu\}^{-1}(0)$
that lie very close to~${\cal Z}_i$,
where $\nu\!\in\!\Ga(\bar{\cal M};V)$ is a small generic perturbation of~$s$.
The manifolds ${\cal Z}_i$ we encounter fall in one of 
the two categories described below.

\begin{dfn}
\label{euler_dfn4}
Suppose $\bar{\cal M}$ is an ms-manifold of dimension~$2n$,
$V\!\lra\!\bar{\cal M}$ is an ms-bundle of rank~$n$,
$s\!\in\!\Ga(\bar{\cal M};V)$, and 
${\cal Z}\!\subset\!s^{-1}(0)$.\\
(1) ${\cal Z}$ is \under{$s$-hollow} 
if there exist a normal-bundle model $(F,Y,\vt)$ for ${\cal Z}$
and a bundle isomorphism 
\hbox{$\vt_V\!\!:\vt^*V\!\lra\!\pi_F^*V$},
covering the identity on $F_{\de}\!-\!(Y\!-\!{\cal Z})$, such~that\\
(1a) $\vt_V|_{F_{\de}-Y-{\cal Z}}$ is smooth and
$\vt_V|_{\cal Z}$ is the identity;\\
(1b) the map $\phi\equiv\vt_V\circ\vt^*s\!: 
          F_{\de}\!-\!(Y\!-\!{\cal Z})\!\lra\!V$ is hollow.\\
(2) ${\cal Z}$ is \under{$s$-regular} if there exist 
a normal-bundle model $(F,Y,\vt)$ for ${\cal Z}$
with closure~$(\bar{\cal Z},\tilde{F},\pi)$,
regular polynomial $\al\!:\tilde{F}\!\lra\!\pi^*V$, and 
a bundle isomorphism 
$\vt_V\!\!:\vt^*V\!\lra\!\pi_F^*V$
covering the identity on~$F_{\de}\!-\!(Y\!-\!{\cal Z})$, such~that\\
(2a) $\vt_V|_{F_{\de}-Y-{\cal Z}}$ is smooth and
$\vt_V|_{\cal Z}$ is the identity;\\
(2b)  $\al|_{\cal Z}$ is nondegenerate and is the resolvent for
$\phi\!\equiv\!\vt_V\!!\circ\vt^*s\!: 
          F_{\de}\!-\!(Y\!-\!{\cal Z})\!\lra\!V$, 
and $Y$ is preserved under  scalar multiplication
in each of the components of $F$ for the splitting 
corresponding to $\al$ as in (1) of Definition~\ref{top_dfn1b}.
\end{dfn}

\begin{prp}
\label{euler_prp}
Let $V\!\lra\!\bar{\cal M}$ be an ms-bundle of rank~$n$
over an ms-manifold of dimension~$2n$.
Suppose ${\cal U}$ is an open subset of ${\cal M}$
and $s\!\in\!\Ga(\bar{\cal M};V)$ is such that
$s|_{\cal U}$ is transversal to the zero~set.\\
(1) If $s^{-1}(0)\cap{\cal U}$ is a finite set,
$^{\pm}|s^{-1}(0)\cap{\cal U}|
=\lan e(V),[\bar{\cal M}]\ran-{\cal C}_{\bar{\cal M}-{\cal U}}(s)$.\\
(2) If $\bar{\cal M}-{\cal U}=
  \bigsqcup\limits_{i=1}\limits^{i=k}\!{\cal Z}_i$,
where each ${\cal Z}_i$ is $s$-regular or $s$-hollow, then 
$s^{-1}(0)\cap{\cal U}$ is finite, and 
$$^{\pm}\big|s^{-1}(0)\cap{\cal U}\big|
=\big\lan e(V),[\bar{\cal M}]\big\ran-{\cal C}_{\bar{\cal M}-{\cal U}}(s)
=\big\lan e(V),[\bar{\cal M}]\big\ran-
\sum_{i=1}^{i=k}{\cal C}_{{\cal Z}_i}(s).$$
If ${\cal Z}_i$ is $s$-hollow, ${\cal C}_{{\cal Z}_i}(s)\!=\!0$.
If ${\cal Z}_i$ is $s$-regular 
and \hbox{$\al_i\!:\tilde{F}_i\!\lra\!V$} is the corresponding polynomial,
$${\cal C}_{{\cal Z}_i}(s)=
^{\pm}\!\big|\{\ups\!\in\!\tilde{F}_i\!:
 \bar{\nu}_{\ups}\!+\!\al_i(\ups)=0\}\big|
\equiv N(\al_i),$$
where $\bar{\nu}\!\in\!\Ga(\bar{\cal Z}_i;V)$ is a generic section.
Finally, if $\al_i\!\in\!
\Ga(\bar{\cal Z}_i;\tilde{F}_i^{*\otimes k}\otimes\pi^* V)$  
has constant rank  over~$\bar{\cal Z}_i$
and factors through a $\tilde{k}$-to-$1$ cover
$\rho_i\!:\tilde{F}_i\!\lra\!\tilde{F}_i^{\otimes k}$,
$${\cal C}_{{\cal Z}_i}(s)=\tilde{k}
\big\lan e\big(\pi^*V/\al_i(\tilde{F}_i)\big),[\bar{\cal Z}_i]\big\ran.$$
\end{prp}

\noindent
This is Corollary~\ref{enum1-euler_crl} in~\cite{Z2}.
Proposition~\ref{euler_prp} reduces the problem of 
computing ${\cal C}_{{\cal Z}_i}(s)$ for
an $s$-regular manifold ${\cal Z}_i$
to counting the zeros of a polynomial map between two vector bundles. 
The general setting for the latter problem is the following.
Suppose  $E,{\cal O}\!\lra\!\bar{\cal M}$ are ms-bundles such that
$\rk~E\!+\!\frac{1}{2}\dim\bar{\cal M}\!=\!\rk~{\cal O}$,
and $\al\!:E\!\lra\!{\cal O}$ is a regular polynomial.
Let $\bar{\nu}\!\in\!\Ga(\bar{\cal M};{\cal O})$
be such that the~map
$$\psi_{\al,\bar{\nu}}\!\equiv\!\bar{\nu}\!+\!\al\!: E\!\lra\!{\cal O}$$
is transversal to the zero set in ${\cal O}$ 
on~$E|{\cal M}$, and all its zeros are contained in~$E|{\cal M}$.
Then \hbox{$N(\al)\!\equiv^{\pm}\!|\psi_{\al,\bar{\nu}}^{-1}(0)|$}
depends only on~$\al$.
If the rank of $E$ is zero, then clearly
$$N(\al)=^{\pm}\!\big|\psi_{\al,\bar{\nu}}^{-1}(0)\big|
=\big\lan e({\cal O}),[\bar{\cal M}]\big\ran.$$
If the rank of $E$ is positive and $\bar{\nu}$ is generic,
it does not vanish and thus determines 
a trivial line subbundle $\Bbb{C}\bar{\nu}$ of~${\cal O}$.
Let ${\cal O}^{\perp}\!=\!{\cal O}/\Bbb{C}\bar{\nu}$
and denote by $\al^{\perp}$ the composition of $\al$ with
the quotient projection map.
If $E$ is a line bundle and $\al$ is~linear,
$$N(\al)=^{\pm}\!\big|\psi_{\al,\bar{\nu}}^{-1}(0)\big|
=\big\lan e(E^*\otimes{\cal O}^{\perp}),[\bar{\cal M}]\big\ran
-{\cal C}_{\al^{-1}(0)}(\al^{\perp}).$$
By Proposition~\ref{euler_prp}, computation of 
${\cal C}_{\al^{-1}(0)}(\al^{\perp})$ again involves counting
the zeros of polynomial maps, but with the rank of the new target bundle, 
i.e.~$E^*\otimes{\cal O}^{\perp}$, one less than the rank
of the original one, i.e.~${\cal O}$.
Subsection~\ref{enum1-zeros_sec} in~\cite{Z2}
reduces the problem of determining $N(\al)$ in all other cases
to the case $E$ is a line bundle and $\al$ is linear.
Thus, at least in reasonably good cases, $N(\al)$ can be determined
after a finite number of steps.\\

\noindent
The next lemma summarizes the results of 
Subsection~\ref{enum1-zeros_sec} in~\cite{Z2}
in the case the original map $\al\!:E\!\lra\!{\cal O}$ is linear.
This case suffices for our purposes.
We denote by $\al^E\!\in\!\Ga(\Bbb{P}E;
\ga_E^*\otimes\pi_E^*{\cal O})$
the section induced by~$\al$.
\hbox{Let $\la_E\!=\!c_1(\ga_E^*)$}.

\begin{lmm}
\label{zeros_main}
Suppose  $\bar{\cal M}$ is  an ms-manifold and
$E,{\cal O}\lra\bar{\cal M}$ are ms-bundles such that
$$\rk~E+\frac{1}{2}\dim\bar{\cal M}= \rk~{\cal O}\equiv m.$$
If
$\al\!\in\!\Ga(\bar{\cal M};E^*\otimes{\cal O})$ and
$\bar{\nu}\!\in\!\Ga(\bar{\cal M};{\cal O})$
are such that $\al$ is regular,
$\bar{\nu}$ has no zeros, the map
$$\psi_{\al,\bar{\nu}}\!\equiv\!\bar{\nu}\!+\!\al\!: E\lra{\cal O}$$
is transversal to the zero set on~$E|{\cal M}$, and all its zeros
are contained in~$E|{\cal M}$, then
$\psi_{\al,\bar{\nu}}^{-1}(0)$ is a finite set,
$^{\pm}\!|\psi_{\al,\bar{\nu}}^{-1}(0)|$ depends only on $\al$, and
$$N(\al)\equiv ^{\pm}\!|\psi_{\al,\bar{\nu}}^{-1}(0)|=
\sum_{k=0}^{k=m-1}
\big\lan c_k({\cal O})\la_E^{m-1-k},[\Bbb{P}E]\big\ran
-{\cal C}_{\al^{E~-1}(0)}(\al^{E\perp}).$$
Furthermore, if the rank of $E$ is $n$,
\begin{equation}\label{zeros_main_e}
\la_E^n+\sum_{k=1}^{k=n}c_k(E)\la_E^{n-k}=0\in H^{2n}(\Bbb{P}E)
\quad\hbox{and}\quad
\big\lan \mu\la_E^{n-1},[\Bbb{P}E]\big\ran=
\big\lan\mu,[\bar{\cal M}]\big\ran~~
\forall\mu\!\in\!H^{2m-2n}(\bar{\cal M}).
\end{equation}
Finally, if $\al$ has constant rank, 
$$N(\al)=^{\pm}\!|\psi_{\al,\bar{\nu}}^{-1}(0)|=
\big\lan e\big({\cal O}/(Im~\al)\big),[\bar{\cal M}]\big\ran.$$
\end{lmm}

\subsection{Review of Notation}
\label{notation_sec}

\noindent
In this subsection, we give a brief description
of the most important notation used in this paper.
See Section~\ref{gl-bubble_sect} in~\cite{Z3} for more details.\\

\noindent
Let $q_N,q_S\!: \Bbb{C}\!\lra\! S^2\!\subset\!\Bbb{R}^3$ 
be the stereographic
projections mapping the origin in $\Bbb{C}$ to the north and 
south poles, respectively.
Denote the south pole of $S^2$, i.e.~the point 
$(0,0,-1)\!\in\!\Bbb{R}^3$, by~$\i$. 
We identify $\Bbb{C}$ with  $S^2-\{\i\}$ via the map~$q_N$.\\

\noindent
If $N$ is any nonnegative integer, let $[N]\!=\!\{1,\ldots,N\}$.
If $I_1$ and $I_2$ are two sets, we denote the disjoint
union of $I_1$ and $I_2$ by $I_1\!+\!I_2$.

\begin{dfn}
\label{rooted_tree}
A finite partially ordered set $I$ is a \under{linearly ordered set}
if for all \hbox{$i_1,i_2,h\!\in\! I$} such that $i_1,i_2\!<\!h$, 
either $i_1\!\le\! i_2$ \hbox{or $i_2\!\le\! i_1$.}\\
A linearly ordered set $I$ is a \under{rooted tree} if
$I$ has a unique minimal element, 
i.e.~there exists \hbox{$\hat{0}\!\in\! I$} such that $\hat{0}\!\le\! i$ 
for {all $i\!\in\! I$}.
\end{dfn}

\noindent
If $I$ is a linearly ordered set, let $\hat{I}$ be
the subset of the non-minimal elements of~$I$.
For every $h\!\in\!\hat{I}$,  denote by $\io_h\!\in\!I$
the largest element of $I$ which is smaller than~$h$.
We call $\io\!:\hat{I}\!\lra\!I$ the attaching map of~$I$.
Suppose $I\!=\!\bigsqcup\limits_{k\in K}\!I_k$
is the splitting of $I$ into rooted trees 
such that $k$ is the minimal element of~$I_k$.
If $\hat{1}\!\not\in\!I$, we define
the linearly ordered set $I\!+_k\!\hat{1}$ to
be the set $I\!+\!\hat{1}$ with all partial-order relations of~$I$
along with the relations
$$k<\hat{1},\qquad \hat{1}<h\hbox{~~if~~}h\!\in\!\hat{I}_k.$$
If $I$ is a rooted tree, we write $I\!+\!\hat{1}$ 
\hbox{for $I\!+_k\!\hat{1}$}.
\\

\noindent
If $S=\Si$ or $S=S^2$ and $M$ is a finite set, 
a {\it $\PP$-valued bubble map with $M$-marked 
points} is a tuple
$$b=\big(S,M,I;x,(j,y),u\big),$$
where $I$ is a linearly ordered set, and
$$x\!:\hat{I}\!\lra\!S\cup S^2,~~j\!:M\!\lra\! I,~~ 
y\!:M\!\lra\!S\cup S^2,\hbox{~~and~~} 
u\!:I\!\lra\!C^{\i}(S;\PP)\cup C^{\i}(S^2;\PP)$$
are maps such that
$$x_h\in
\begin{cases}
S^2-\{\i\},&\hbox{if~}\io_h\!\in\!\hat{I};\\
S,&\hbox{if~}\io_h\!\not\in\!\hat{I},
\end{cases}\qquad
y_l\in
\begin{cases}
S^2-\{\i\},&\hbox{if~}j_l\!\in\!\hat{I};\\
S,&\hbox{if~}j_l\!\not\in\!\hat{I},\qquad
\end{cases}
u_i\in
\begin{cases}
C^{\i}(S^2;\PP),&\hbox{if~}i\!\in\!\hat{I};\\
C^{\i}(S;\PP),&\hbox{if~}i\!\not\in\!\hat{I},
\end{cases}$$  
and  $u_h(\i)\!=\!u_{\io_h}(x_h)$ for all 
$h\!\in\!\hat{I}$.
We associate such a tuple with Riemann surface
$$\Si_b=
\Big(\bigsqcup_{i\in I}\Si_{b,i}\Big)\Big/\!\sim,
\hbox{~~where}\qquad 
\Si_{b,i}=
\begin{cases}
\{i\}\times S^2,&\hbox{if~}i\!\in\!\hat{I};\\
\{i\}\times S,&\hbox{if~}i\!\not\in\!\hat{I},
\end{cases}
\quad\hbox{and}\quad
(h,\i)\sim (\io_h,x_h)
~~\forall h\!\in\!\hat{I},$$
with marked points $(j_l,y_l)\!\in\!\Si_{b,j_l}$,
and continuous map $u_b\!:\Si_b\!\lra\!\PP$,
given by $u_b|\Si_{b,i}\!=\!u_i$ for \hbox{all $i\!\in\! I$}.
We require that all the singular points of $\Si_b$,
i.e.~$(\io_h,x_h)\!\in\!\Si_{b,\io_h}$ for $h\!\in\!\hat{I}$,
and all the marked points be distinct.
Furthermore, if $S\!=\!S^2$, all these points are to be different
from the special marked point $(\hat{0},\i)\!\in\!\Si_{b,\hat{0}}$.
In addition, if $\Si_{b,i}\!=\!S^2$ and 
$u_{i*}[S^2]\!=\!0\!\in\! H_2(\PP;\Bbb{Z})$,
then $\Si_{b,i}$ must contain at least two singular and/or
marked points of~$\Si_b$ other than~$(i,\i)$.
Two bubble maps $b$ and $b'$ are {\it equivalent} if there
exists a homeomorphism $\phi\!:\Si_b\!\lra\!\Si_{b'}$
such that $u_b\!=\!u_{b'}\!\circ\phi$,
$\phi(j_l,y_l)\!=\!(j_l',y_l')$ for \hbox{all $l\!\in\! M$},
$\phi|_{\Si_{b,i}}$ is holomorphic for \hbox{all $i\!\in\!I$},
and $\phi|_{\Si_{b,i}}\!=\!Id$ if $S\!=\!\Si$ 
\hbox{and $i\!\in\! I\!-\!\hat{I}$.}\\

\noindent
The general structure of bubble maps is described
by tuples ${\cal T}\!=\!(S,M,I;j,\under{d})$,
with \hbox{$d_i\!\in\!\Bbb{Z}$} 
describing the degree of the map $u_b$ on~$\Si_{b,i}$.
We call such tuples {\it bubble types}.
Bubble type ${\cal T}$ is {\it simple} if  $I$ is a rooted tree;
${\cal T}$ is {\it is basic} if $\hat{I}\!=\!\eset$;
${\cal T}$ is {\it semiprimitive} if 
$\io_h\!\not\in\!\hat{I}$ for all $h\!\in\!\hat{I}$.
We call semiprimitive bubble type ${\cal T}$
{\it primitive} if $j_l\!\in\!\hat{I}$ for \hbox{all $j_l\!\in\! M$}.
The above equivalence relation on the set of bubble maps
induces an equivalence relation on the set of bubble types.
For each $h,i\!\in\! I$, let
\begin{gather*}
D_i{\cal T}=\{h\!\in\!\hat{I}\!:i\!<\!h\},\quad
\bar{D}_i{\cal T}=D_i{\cal T}\cup\{i\},\quad
H_i{\cal T}=\{h\!\in\!\hat{I}\!:\io_h\!=\!i\},\quad
M_i{\cal T}=\{l\!\in\! M\!:j_l\!=\!i\},\\
\chi_{\cal T}h=\begin{cases}
0,&\hbox{if~} \forall i\!\in\! 
             I\hbox{~s.t.~}h\!\in\!\bar{D}_i{\cal T}, d_i\!=\!0;\\
1,&\hbox{if~} d_h\!\neq\! 0,\hbox{~but~}
\forall i\!\in\! I\hbox{~s.t.~}h\!\in\! D_i{\cal T}, d_i\!=\!0;\\
2,&\hbox{otherwise};
\end{cases}\qquad
\chi({\cal T})=\big\{h\!\in\!I\!:\chi_{\cal T}h\!=\!1\big\}.
\end{gather*}
Denote by ${\cal H}_{\cal T}$ the space of all holomorphic
bubble maps with structure~${\cal T}$.\\

\noindent
The automorphism group of every bubble type ${\cal T}$ we encounter
in the next two sections is trivial.
Thus, every bubble type discussed below is presumed to be 
automorphism-free.\\

\noindent
If $S\!=\!\Si$, we denote by ${\cal M}_{\cal T}$ the set of
equivalence classes of bubble maps in~${\cal H}_{\cal T}$.
Then there exists ${\cal M}_{\cal T}^{(0)}\!\subset\!{\cal H}_{\cal T}$
such that ${\cal M}_{\cal T}$ is the quotient of ${\cal M}_{\cal T}^{(0)}$ 
by an~$(S^1)^{\hat{I}}$-action.
Corresponding to this action, we obtain 
$|\hat{I}|$ line orbi-bundles 
$\{L_h{\cal T}\!\!\lra\!{\cal M}_{\cal T}\!:h\!\in\!\hat{I}\}$.
The bundle of gluing parameters in the case~$S\!=\!\Si$~is
$$F{\cal T}=\bigoplus_{h\in\hat{I}}F_h{\cal T},
\hbox{~~where}\qquad
F_{h,[b]}{\cal T}=\begin{cases}
L_{h,[b]}{\cal T}\otimes L_{\io_h,[b]}^*{\cal T},&
\hbox{if~}\io_h\!\in\!\hat{I};\\
L_{h,[b]}{\cal T}\otimes T_{x_h}\Si,&
\hbox{if~}\io_h\!\not\in\!\hat{I}.
\end{cases}$$
Let
$F^{\eset}{\cal T}\!=
\!\{\ups\!=\!(\ups_h)_{h\in{\hat{I}}}\!\in\! F{\cal T}\!: 
\ups_h\!\neq\!0~\forall h\!\in\!\hat{I}\}$.\\

\noindent
For each bubble type ${\cal T}=(S^2,M,I;j,d)$,
let
$${\cal U}_{\cal T}=\big\{[b]\!:
b\!=\!\big(S^2,M,I;x,(j,y),u\big)\!\in\!{\cal H}_{\cal T},~
u_{i_1}(\i)=u_{i_2}(\i)~\forall i_1,i_2\!\in\! I\!-\!\hat{I}
\big\}.$$
Similarly to the $S\!=\!\Si$ case above,
${\cal U}_{\cal T}$ is the quotient of a subset ${\cal B}_{\cal T}$
of ${\cal H}_{\cal T}$ by a
\hbox{$\tilde{G}_{\cal T}\!\equiv\!(S^1)^I$}-action.
Denote by ${\cal U}_{\cal T}^{(0)}$ the quotient of ${\cal B}_{\cal T}$
by 
\hbox{$G_{\cal T}\!\equiv\!(S^1)^{\hat{I}}\!\subset\!\tilde{G}_{\cal T}$}.
Then ${\cal U}_{\cal T}$ is the quotient of
${\cal U}_{\cal T}^{(0)}$ the residual 
$G_{\cal T}^*\!\equiv\!(S^1)^{I-\hat{I}}\!\subset\!\tilde{G}_{\cal T}$
action.
Corresponding to these quotients, we obtain line orbi-bundles
$\{L_h{\cal T}\!\!\lra\!{\cal U}_{\cal T}^{(0)}\!\!:h\!\in\!\hat{I}\}$
and $\{L_i{\cal T}\!\!\lra\!{\cal U}_{\cal T}\!\!:i\!\in\! I\}$.
Let 
\begin{gather*}
F{\cal T}=\bigoplus_{h\in\hat{I}}F_h{\cal T}\lra{\cal U}_{\cal T}^{(0)},
\hbox{~~where}\quad
F_{h,[b]}{\cal T}=\begin{cases}
L_{h,[b]}{\cal T}\otimes L_{\io_h,[b]}^*{\cal T},&
\hbox{if~}\io_h\!\in\!\hat{I};\\
L_{h,[b]}{\cal T},&
\hbox{if~}\io_h\!\not\in\!\hat{I};
\end{cases}\\
{\cal FT}=\bigoplus_{h\in\hat{I}}{\cal F}_h{\cal T}\lra{\cal U}_{\cal T},
\hbox{~~where}\quad
{\cal F}_{h,[b]}{\cal T}=L_{h,[b]}{\cal T}\otimes L_{\io_h,[b]}^*{\cal T}.
\end{gather*}
The bundle of gluing parameters in the case $S\!=\!S^2$
is~${\cal FT}$.\\

\noindent
Gromov topology on the space of equivalence classes of bubble maps
induces a partial ordering on the set of bubble types
and their equivalence classes such that the spaces
$$\bar{\cal M}_{\cal T}=
\bigcup_{{\cal T}'\le{\cal T}}{\cal M}_{{\cal T}'},\quad
\bar{\cal U}_{\cal T}^{(0)}=
\bigcup_{{\cal T}'\le{\cal T}}{\cal U}_{{\cal T}'}^{(0)}\quad
\bar{\cal U}_{\cal T}=
\bigcup_{{\cal T}'\le{\cal T}}{\cal U}_{{\cal T}'}$$
are compact and Hausdorff.
The $G_{\cal T}^*$-action on ${\cal U}_{\cal T}^{(0)}$ extends
to an action on $\bar{\cal U}_{\cal T}^{(0)}$,
and thus line orbi-bundles $L_{{\cal T},i}\lra{\cal U}_{\cal T}$
with $i\!\in\! I\!-\!\hat{I}$ extend over~$\bar{\cal U}_{\cal T}$.
The evaluation maps
$$\ev_l\!:{\cal H}_{\cal T}\lra\PP,\quad
\ev_l\big((S,M,I;x,(j,y),u)\big)=
u_{j_l}(y_l),$$
descend to all the quotients and induce continuous
maps on $\bar{\cal M}_{\cal T}$, $\bar{\cal U}_{\cal T}$,
and~$\bar{\cal U}_{\cal T}^{(0)}$.
If $\mu\!=\!\mu_M$ is an $M$-tuple of submanifolds of~$\PP$,
let
$${\cal M}_{\cal T}(\mu)=
\big\{b\!\in\!{\cal M}_{\cal T}\!:\ev_l(b)\!\in\!\mu_l~\forall 
l\!\in\! M\big\}$$
and define spaces ${\cal U}_{\cal T}(\mu)$,
$\bar{\cal U}_{\cal T}(\mu)$, etc.~in a similar way.
If $S=S^2$, we define another evaluation~map,
$$\ev\!: 
{\cal B}_{\cal T}\lra\PP\quad\hbox{by}\quad
\ev\big((S^2,M,I;x,(j,y),u)\big)=u_{\hat{0}}(\i),$$
where $\hat{0}$ is any minimal element of $I$.
This map descends to ${\cal U}_{\cal T}^{(0)}$ and~${\cal U}_{\cal T}$.
If $\mu=\mu_{\tilde{M}}$ is an $\tilde{M}$-tuple of constraints, let
$${\cal U}_{\cal T}(\mu)=
\big\{b\!\in\!{\cal U}_{\cal T}\!:
\ev_l(b)\!\in\!\mu_l~\forall l\!\in\!M\cap\tilde{M},~
\ev(b)\!\in\!\mu_l~\forall l\!\in\!M\!-\!\tilde{M}\big\}$$
and define ${\cal U}_{\cal T}^{(0)}(\mu)$, etc.~similarly.
If $S\!=\!\Si$, ${\cal T}$ is a simple bubble type,
and $d_{\hat{0}}\!=\!0$, define
$$\ev\!: 
{\cal H}_{\cal T}\lra\PP\quad\hbox{by}\quad
\ev\big((\Si,M,I;x,(j,y),u)\big)=u_{\hat{0}}(\Si).$$
This map is well-defined, since $u_{\hat{0}}$ is a degree-zero
holomorphic map and thus is constant.\\

\noindent
Suppose ${\cal T}\!=\!(S^2,M,I;j,\under{d})$ is a bubble type, 
$\big\{{\cal T}_k\!=\!(S^2,M_k,I_k;j_k,\under{d}_k)\big\}$
are the corresponding simple types (see~\cite{Z3}),
$k\!\in\! I\!-\!\hat{I}$, and $M_0$ is nonempty subset of~$M_k{\cal T}$.
Let
$${\cal T}/M_0=\big(S^2,\hat{I},M-M_0;j|(M-M_0),d|\hat{I}\big).$$
Define ${\cal T}(M_0)\equiv(S^2,M,\hat{I}+_k\hat{1};j',d')$  by
$$j'_l=\begin{cases}
k,&\hbox{if~}l\in M_0;\\
\hat{1},&\hbox{if~}l\in M_k{\cal T}-M_0;\\
j_l,&\hbox{otherwise};
\end{cases}\qquad
d'_i=\begin{cases}
0,&\hbox{if~}i=k;\\
d_k,&\hbox{if~}i=\hat{1};\\
d_i,&\hbox{otherwise}.
\end{cases}$$
The tuples ${\cal T}/M_0$ and ${\cal T}(M_0)$ are bubble types as long as 
$d_k\!\neq\!0$ or $M_0\!\neq\! M_{\hat{0}}{\cal T}$.
Then,
\begin{equation}
\label{cart_split2}
\bar{\cal U}_{{\cal T}(M_0)}(\mu)=
 \bar{\cal M}_{0,\{\hat{1}\}+M_0}\times \bar{\cal U}_{{\cal T}/M_0}(\mu),
\end{equation}
where $\bar{\cal M}_{0,\{\hat{1}\}+M_0}$ denotes the
Deligne-Mumford moduli space of rational curves with 
$(\{\hat{0},\hat{1}\}+M_0)$-marked points. 
If $l\!\in\! M_k{\cal T}$ for some $k\!\in\! I\!-\!\hat{I}$,
we denote ${\cal T}(\{l\})$ by~${\cal T}(l)$.
Let
\begin{equation}\label{normal_bundle1}
c_1({\cal L}_{{\cal T},k}^*)\equiv c_1(L_{{\cal T},k}^*)-
\sum_{M_0\subset M_k,M_0\neq\eset}
PD_{\bar{\cal U}_{\cal T}(\mu)}\big[\bar{\cal U}_{{\cal T}(M_0)}(\mu)\big]
\in H^2\big(\bar{\cal U}_{\cal T}(\mu)\big).
\end{equation}
If the constraints~$\mu$ are disjoint, 
\hbox{$\bar{\cal U}_{{\cal T}(M_0)}(\mu)=\eset$} if~$|M_0|\!\ge\!2$
and
\begin{equation}\label{normal_bundle2}
\big[\bar{\cal U}_{{\cal T}(l)}(\mu)\big] 
                        \cap c_1({\cal L}_{{\cal T},k}^*)=
\big[\bar{\cal U}_{{\cal T}(l)}(\mu)\big]
                        \cap c_1(L_{{\cal T}(l),\hat{1}}^*)=
\big[\bar{\cal U}_{{\cal T}(l)}(\mu)\big]
                    \cap c_1({\cal L}_{{\cal T}(l),\hat{1}}^*).
\end{equation}
See~Subsection~\ref{enum1-comp_sec2} in~\cite{Z2}.\\

\noindent
We are now ready to explain the statement of Theorem~\ref{g3n2_thm}.
Let $d$ and $\mu$ be as in Subsection~\ref{back}.
If $k$ is a positive integer, let $\bar{\cal V}_k(\mu)$ denote 
the disjoint union of the spaces $\bar{\cal U}_{\cal T}(\mu)$
taken over equivalence classes of basic bubble types
${\cal T}\!=\!(S^2,[N],I;j,\under{d})$ \hbox{with $|I|\!=\!k$}
and $\sum d_k\!=\!d$.
Similarly, we denote by ${\cal V}_k(\mu)$ 
the subspace of $\bar{\cal V}_k(\mu)$ consisting 
of the spaces ${\cal U}_{\cal T}(\mu)$ with ${\cal T}$ as above.
Note that the dimension of $\bar{\cal V}_k(\mu)$ over~$\Bbb{R}$
is~$12\!-\!4k$.
Let 
$$a=\ev^*c_1(\ga_{\PP}^*)\in 
H^2(\bar{\cal V}_k(\mu);\Bbb{Z}),\qquad
c_1({\cal L}^*)=c_1\big({\cal L}_{(S^2,[N],\{\hat{0}\};\hat{0},d),\hat{0}}^*
\big) \in H^2(\bar{\cal V}_1(\mu);\Bbb{Z}).$$
While the components of  $\bar{\cal V}_2(\mu)$ are unordered,
we can still define the chern classes 
$$c_1({\cal L}_1^*)+c_1({\cal L}_2^*),~
c_1^2({\cal L}_1^*)+c_1^2({\cal L}_2^*),~
c_1({\cal L}_1^*)c_1({\cal L}_2^*)
\in H^*\big(\bar{\cal V}_2(\mu)\big).$$
In the notation of the previous paragraph,
$c_1({\cal L}_i^*)$ denotes the cohomology 
class~$c_1({\cal L}_{{\cal T}_{k_i},{k_i}}^*)$,
where we write $I\!=\!\{k_1,k_2\}$.\\

\noindent
There are generalizations of the splitting~\e_ref{cart_split2}
that are useful in computations.
Let ${\cal T}$ and $\{{\cal T}_k\}$ be as above.
Suppose $k\!\in\!I\!-\!\hat{I}$ and $d_k\!=\!0$.
Denote by $\bar{\cal T}$ the bubble type obtained
from ${\cal T}$ by removing $k$ from $I$ and $M_k{\cal T}$ from~$M$.
Then,
\begin{equation}\label{cart_split1}
{\cal U}_{\cal T}(\mu)=
{\cal M}_{0,H_k{\cal T}+M_k{\cal T}}\times
{\cal U}_{\bar{\cal T}}(\mu),
\end{equation}
where ${\cal M}_{0,H_k{\cal T}+M_k{\cal T}}$ denotes the main
stratum of $\bar{\cal M}_{0,H_k{\cal T}+M_k{\cal T}}$.
It is shown in \cite{Z3}, that 
${\cal M}_{0,H_k{\cal T}+M_k{\cal T}}$ is a quotient of a subset 
${\cal M}_{0,H_k{\cal T}+M_k{\cal T}}^{(0)}$
of $\Bbb{C}^{H_k{\cal T}+M_k{\cal T}}$ by the diagonal $S^1$-action.
The closure $\tilde{\cal M}_{0,H_k{\cal T}+M_k{\cal T}}^{(0)}$
of ${\cal M}_{0,H_k{\cal T}+M_k{\cal T}}^{(0)}$
is $S^1$-equivarently diffeomorphic 
to~$S^{2(|H_k{\cal T}+M_k{\cal T}|)-3}$;
see Subsection~\ref{enum1-chern_class}
for the case $|H_k{\cal T}\!+\!M_k{\cal T}|\!=\!3$.
Thus, ${\cal M}_{0,H_k{\cal T}+M_k{\cal T}}$
admits a compactification
$$\tilde{\cal M}_{0,H_k{\cal T}+M_k{\cal T}}=
\tilde{\cal M}_{0,H_k{\cal T}+M_k{\cal T}}^{(0)}/S^1
\approx\Bbb{P}^{|H_k{\cal T}+M_k{\cal T}|-2}.$$
Via the splitting~\e_ref{cart_split1},
we obtain a compactification of ${\cal U}_{\cal T}(\mu)$:
\begin{equation}\label{cart_split}
\tilde{\cal U}_{\cal T}(\mu)\equiv
\tilde{\cal M}_{0,H_k{\cal T}+M_k{\cal T}}\times
\bar{\cal U}_{\bar{\cal T}}(\mu).
\end{equation}
Note that 
$\tilde{\cal U}_{\cal T}(\mu)\!=\!\bar{\cal U}_{\cal T}(\mu)$
if the cardinality of $H_k{\cal T}\!+\!M_k{\cal T}$ is two or three.
Furthermore, in all cases, the restrictions of
$L_{{\cal T},k}$ and the tautological line bundle of 
$\tilde{\cal M}_{0,H_k{\cal T}+M_k{\cal T}}$
to ${\cal U}_{\bar{\cal T}}(\mu)$ agree.\\

\noindent
Finally, if $X$ is any space, $F\!\lra\! X$ is a normed vector bundle,
and $\de\!: X\!\lra\!\Bbb{R}$ is any function,~let
$$F_{\de}=\big\{(b,v)\!\in\! F\!: |v|_b<\de(b)\big\}.$$
Similarly, if $\Om$ is a subset of $F$, let 
$\Om_{\de}=F_{\de}\cap\Om$.
If $\ups=(b,v)\!\in\! F$, denote by $b_{\ups}$ the image of $\ups$
under the bundle projection map, i.e.~$b$ in this case.

\subsection{Spaces of Rational Maps}
\label{str_global_sec}

\noindent
In this subsection, we describe the structure of various spaces
of bubble maps passing through the points $\mu_1,\ldots,\mu_N$.
The main goal is to describe the behavior of certain bundle sections 
over such spaces near the boundary strata.

\begin{lmm}
\label{flat_metrics}
There exist $r_{\PP}\!>\!0$ and a smooth family of Kahler metrics
$\{g_{\PP,q}\!: q\!\in\!\P\}$ on~$\PP$ with the following property.
If $B_q(q',r)\!\subset\!\PP$ denotes the $g_{\PP,q}$-geodesic ball about~$q'$
of radius~$r$, the triple $(B_q(q,r_{\PP}),J,g_{\PP,q})$ is isomorphic
to a ball in $\Bbb{C}^2$ for \hbox{all $q\!\in\!\PP$}.
\end{lmm}

\noindent
This is the $n\!=\!2$ case of Lemma~\ref{enum1-flat_metrics} in~\cite{Z2}.
If $b\!=\!\big(S^2,M,I;x,(j,y),u\big)\!\in\!
{\cal B}_{\cal T}$, $m\!\ge\!1$, and $k\!\in\! I$, let
$${\cal D}_{{\cal T},k}^{(m)}b=
\frac{2}{(m-1)!}
\frac{D^{m-1}}{ds^{m-1}}\frac{d}{ds}(u_k\circ q_S)
\Big|_{(s,t)=0},$$
where the covariant derivatives are taken with respect to the metric
$g_{\PP,b}\!\equiv\!g_{\PP,\ev(b)}$ and $s\!+\!it\!\in\!\Bbb{C}$.
If ${\cal T}^*$ is a basic bubble type, 
the maps ${\cal D}_{{\cal T},k}^{(m)}$ with 
${\cal T}\!<\!{\cal T}^*$ and $k\!\in\!I\!-\!\hat{I}$
induce a continuous section of $\ev^*T\PP$ over 
$\bar{\cal U}_{{\cal T}^*}^{(0)}$ and a continuous section of the bundle 
$L_{{\cal T}^*,k}^{*\otimes m}\otimes\ev^*T\PP$ over 
$\bar{\cal U}_{{\cal T}^*}$, described~by 
$${\cal D}_{{\cal T}^*,k}^{(m)}[b,c_k]
=c_k^m{\cal D}_{{\cal T},k}^{(m)}b,
\quad\hbox{if}~~b\in{\cal U}_{\cal T}^{(0)},~c_k\in\Bbb{C}.$$ 
We will often write 
${\cal D}_{{\cal T},k}$ instead of~${\cal D}_{{\cal T},k}^{(1)}$.
If~${\cal T}$ is simple, we will abbreviate 
${\cal D}_{{\cal T},k}^{(m)}$ as~${\cal D}^{(m)}$.
If ${\cal T}\!=\!(\Si,[N],I;j,\under{d})$ is a simple bubble type
and $k\!\in\!\hat{I}$, 
let ${\cal D}_{{\cal T},k}^{(m)}$ denote 
the section~${\cal D}_{\bar{\cal T},k}^{(m)}$.
Finally, fix a real number $p\!>\!2$.

\begin{thm}
\label{str_global}
Suppose $d$ is a positive integer, $N\!=\!3d\!-\!4$, 
$M_0$ is a subset of~$[N]$, and
\hbox{$\mu\!=\!(\mu_1,\ldots,\mu_N)$} is an $N$-tuple
of points in general position in~$\PP$.
If ${\cal T}^*\!\!=\!(S^2,[N]\!-\!M_0,I^*;j^*,\under{d}^*)$ is a 
basic bubble type such that~$\sum d_i^*\!=\!d$,
the space $\bar{\cal U}_{{\cal T}^*}(\mu)$ 
is an  ms-manifold of dimension $2(6\!-\!2|I^*|\!-\!|M_0|)$
and $L_{{\cal T}^*,k}$ for $k\!\in\!I^*$ and $\ev^*T\PP$
are ms-bundles over~$\bar{\cal U}_{{\cal T}^*}(\mu)$.
If ${\cal T}\!=\!(S^2,[N]\!-\!M_0,I;j,\under{d})\!<\!{\cal T}^*$, 
there exist $\de,C\in C^{\i}\big({\cal U}_{\cal T}(\mu);\Bbb{R}^+\big)$
and a homeomorphism
$$\ga_{\cal T}^{\mu}\!: 
     {\cal FT}_{\de}\lra \bar{\cal U}_{{\cal T}^*}(\mu),$$
onto an open neighborhood of ${\cal U}_{\cal T}(\mu)$
in $\bar{\cal U}_{{\cal T}^*}(\mu)$ such that
$\ga_{\cal T}^{\mu}|{\cal U}_{\cal T}(\mu)$ is the identity
and $\ga_{\cal T}^{\mu}|{\cal F}^{\eset}{\cal T}_{\de}$
is an orientation-preserving diffeomorphism onto an open subset 
of~${\cal U}_{\cal T}(\mu)$.
Furthermore, with appropriate identifications,
\begin{gather*}
\Big|{\cal D}_{{\cal T}^*,k}\ga_{\cal T}^{\mu}(\ups)
-\al_{{\cal T},k}\big(\rho_{\cal T}(\ups)\big)\Big|
\le C(b_{\ups})|\ups|^{\frac{1}{p}} \big|\rho_{\cal T}(\ups)\big|
\quad\forall\ups\!\in\!{\cal FT}_{\de},
\qquad\hbox{where}\\
\rho_{\cal T}(\ups)=
\big(b,(\tilde{v}_h)_{h\in\chi({\cal T})}\big)
\in\tilde{\cal F}{\cal T}\equiv
\bigoplus_{h\in\chi({\cal T})}
L_h{\cal T}\otimes L_{\tilde{\io}_h}^*{\cal T};~~
\tilde{v}_h=\prod_{i\in\hat{I},h\in\bar{D}_i{\cal T}}\!\!\!\!v_i;~~
\tilde{\io}_h\!\in\!I-\hat{I},~ 
h\!\in\!\bar{D}_{\tilde{\io}_h}{\cal T}\\
\al_{{\cal T},k}\big(b,(\tilde{v}_h)_{h\in\chi({\cal T})}\big)
=\sum_{h\in I_k\cap\chi({\cal T})}
{\cal D}_{{\cal T},h}\tilde{v}_h.
\end{gather*}
and $I_k\!\subset\!I$ is the rooted tree containing~$k$.
\end{thm}

\noindent
This is a special case of Theorem~\ref{enum1-str_global} in~\cite{Z2};
see also the remark following the theorem.
The analytic estimate on~${\cal D}_{{\cal T}^*,k}$
is used frequently in the next two sections.
If ${\cal T}$ is semiprimitive,
the bundle ${\cal FT}\!=\!\tilde{\cal F}{\cal T}$
and the section $\al_{\cal T}\!=\al_{\cal T}\!\circ\!\rho_{\cal T}$
extend over $\bar{\cal U}_{\cal T}(\mu)$
via the decomposition~\e_ref{cart_split}.
In terms of the notions of Subsection~\ref{top_sec},
$({\cal FT},{\cal FT}-F^{\eset}{\cal T},\ga_{\cal T}^{\mu})$
is a normal-bundle model for 
${\cal U}_{\cal T}(\mu)\subset\bar{\cal U}_{{\cal T}^*}(\mu)$.
This normal-bundle model admits a closure if ${\cal T}$ is semiprimitive.\\

\noindent
We will need a similar description for spaces of stable maps
corresponding to the rational degree-$d$ curves 
with certain singularities that pass through 
the $3d\!-\!4$ points $\mu_1,\ldots,\mu_N$.
If ${\cal T}\!=\!(S^2,M,I;j,\under{d})$ is a bubble type
and $\chi_{\cal T}h\!=\!1$, let
$$E_h{\cal T}=
\begin{cases}
L_h,&\hbox{if~}h\!\in\!I\!-\!\hat{I};\\
\tilde{\cal F}_h{\cal T},&\hbox{if~}h\!\in\!\hat{I};
\end{cases}
\qquad
E{\cal T}=\bigoplus_{h\in\chi({\cal T})}E_h{\cal T}.$$
If $M\!=\![N]\!-\!M_0$ and $\sum d_i\!=\!d$, put
$${\cal S}_{{\cal T};1}(\mu)=\big\{b\!\in\!{\cal U}_{\cal T}(\mu)\!: 
{\cal D}_{{\cal T},h}b\!=\!0\hbox{~for some~}h\!\in\!\chi({\cal T})\big\}.$$
If $|\chi({\cal T})|\!\ge\!2$, let
$${\cal S}_{{\cal T};2}(\mu)
=\big\{ 
\big[b,(\tilde{v})_{h\in\chi({\cal T})}\big]\!\in\!\Bbb{P}E{\cal T}\!:
\sum_{h\in\chi({\cal T})}{\cal D}_{{\cal T},h}\tilde{v}_h\!=\!0\big\}
-{\cal S}_{{\cal T};1}(\mu).$$
If ${\cal T}$ is basic and $|I|\!=\!1$, the set 
${\cal S}_1(\mu)\!\equiv\!{\cal S}_{{\cal T};1}(\mu)$
of maps can be identified with a dense open subset of
the set of irreducible rational degree-$d$ curves that pass through 
the $3d\!-\!4$ points and have a cusp.
We denote the closure of ${\cal S}_1(\mu)$ in 
$\bar{\cal V}_1(\mu)$ by~$\bar{\cal S}_1(\mu)$.
Let ${\cal S}_{2;2}(\mu)$ be the disjoint union of 
the spaces ${\cal S}_{{\cal T};2}(\mu)$ taken over
all equivalence classes of basic bubble types~${\cal T}$
\hbox{with $|I|\!=\!2$}.
The set ${\cal S}_{2;2}(\mu)$ can be identified with 
a dense open subset of the set of two-component rational degree-$d$ 
curves that pass through  the $3d\!-\!4$ points and 
have a tacnode as a node common to both components.
We denote by $\bar{\cal S}_{2;2}(\mu)$ the closure of 
${\cal S}_{2;2}(\mu)$ in~$\Bbb{P}E_2$,
where $E_2\!\lra\!\bar{\cal V}_2(\mu)$ is the bundle
such that $E_2|\bar{\cal U}_{\cal T}(\mu)\!=\!E{\cal T}$.
Similarly, we denote by ${\cal S}_{2;1}(\mu)$
the disjoint union of the spaces ${\cal S}_{{\cal T};1}(\mu)$ taken over
all equivalence classes of basic bubble types~${\cal T}$
\hbox{with $|I|\!=\!2$}.
This finite set can be identified with a subset of $\bar{\cal S}_{2;2}(\mu)$
as well as with the set of two-component rational degree-$d$ curves
passing through the $3d\!-\!4$ points such that the two components
meet at a node at which one of them has a cusp.\\

\noindent
Suppose ${\cal T}\!=\!(S^2,[N]\!-\!M_0,I;j,\under{d})$ 
is a simple bubble type \hbox{with $\sum d_i\!=\!d$}.
By Corollary~\ref{enum1-reg_crl2} in~\cite{Z2},
the space ${\cal S}_{{\cal T};1}(\mu)$ 
is a smooth complex submanifold of~${\cal U}_{\cal T}(\mu)$.
Let 
$${\cal NS}_{{\cal T};1}=
L_{{\cal T},h_1}^*\otimes\ev^*T\PP\!\lra\!{\cal S}_{{\cal T};1}(\mu)$$
denote its normal bundle, where $h_1\!\in\!\chi({\cal T})$
is such that ${\cal D}_{{\cal T},h_1}b\!=\!0$. 
Put
\begin{gather*}
{\cal FS}_{{\cal T};1}=
\bigoplus_{h\in\hat{I}-\chi({\cal T})}\!\!\!\!\!{\cal F}_h{\cal T}
~~\oplus
\begin{cases}
\{0\},&\hbox{if~}h_1\!\not\in\!\hat{I};\\
{\cal F}_{h_1}{\cal T},
&\hbox{if~}h_1\!\in\!\hat{I}~\&~\chi({\cal T})\!=\!\{h_1\};\\
{\cal F}_{h_1}{\cal T}\oplus L_{{\cal T},h_1}^*\!\!\otimes\! L_{{\cal T},h_2},
&\hbox{if~}\chi({\cal T})\!=\!\{h_1,h_2\};
\end{cases}\\
\tilde{\cal F}{\cal S}_{{\cal T};1}=
\begin{cases}
{\cal F}_{h_1}{\cal T}^{\otimes2}
&\hbox{if~}h_1\!\in\!\hat{I}~\&~\chi({\cal T})\!=\!\{h_1\};\\
{\cal F}_{h_1}{\cal T}^{\otimes2}\oplus
{\cal F}_{h_2}{\cal T},
&\hbox{if~}\chi({\cal T})\!=\!\{h_1,h_2\}.
\end{cases}
\end{gather*}
If $d_{\hat{0}}\!=\!0$,
we define $\rho_{{\cal T};1}\!:{\cal FS}_{{\cal T};1}\!\lra\!
\tilde{\cal F}{\cal S}_{{\cal T};1}$ and
$\al_{{\cal T};1}\!\in\!\Ga\big(\bar{\cal S}_{{\cal T};1}(\mu);
\hbox{Hom}(\tilde{\cal F}{\cal S}_{{\cal T};1};
L_{{\cal T},\hat{0}}^{*\otimes2}\otimes\ev^*T\PP)\big)$ by
\begin{gather*}
\rho_{{\cal T};1}(\ups)=
\begin{cases}
\ups_{h_1}\otimes\ups_{h_1},&
\hbox{if~}\chi({\cal T})\!=\!\{h_1\};\\
\big(\ups_{h_1}\otimes\ups_{h_1},\ups_{h_1}\otimes u\big),&
\hbox{if~}\chi({\cal T})\!=\!\{h_1,h_2\}~\&~
 u\!\in\!L_{{\cal T},h_1}^*\otimes L_{{\cal T},h_2};
\end{cases}\\
\al_{{\cal T};1}(\vp)=
\begin{cases}
{\cal D}_{{\cal T},h_1}^{(2)}\vp,&
\!\hbox{if~}\chi({\cal T})\!=\!\{h_1\};\\
{\cal D}_{{\cal T},h_1}^{(2)}\vp_1\!+\!
x_{h_2}{\cal D}_{{\cal T},h_1}^{(1)}\vp_2,&
\!\hbox{if~}\chi({\cal T})\!=\!\{h_1,h_2\},~
\vp\!=\!(\vp_1,\vp_2),~
b_{\vp}\!=\!\big(S^2,M,I;x,(j,y),u\big).
\end{cases}
\end{gather*}
If $\big|\chi({\cal T})\big|\!\ge\!2$, 
${\cal S}_{{\cal T};2}(\mu)$ is a smooth submanifold
of~$\Bbb{P}E{\cal T}\!\lra\!{\cal U}_{\cal T}(\mu)$.
We identify it with a subset of~${\cal U}_{\cal T}(\mu)$
via the bundle projection map
$\pi_{E{\cal T}}\!:\Bbb{P}E{\cal T}\!\lra\!{\cal U}_{\cal T}(\mu)$.
If $\big|\chi({\cal T})\big|\!=\!3$, 
${\cal S}_{{\cal T};2}(\mu)\!=\!{\cal U}_{\cal T}(\mu)$,
and we put ${\cal NS}_{{\cal T};2}\!=\!\{0\}$.
If $\chi({\cal T})\!=\!\{h_1,h_2\}$,
${\cal S}_{{\cal T};2}(\mu)$ is a smooth submanifold of
${\cal U}_{\cal T}(\mu)$ with normal bundle 
$${\cal NS}_{{\cal T};2}=L_{{\cal T};h_2}^*\otimes
\big(\hbox{Im}~{\cal D}_{{\cal T},h_1}\big)^{\perp}.$$
If $\io_{h_1}\!=\!\io_{h_2}$ for all $h_1,h_2\!\in\!\chi({\cal T})$, put 
$${\cal FS}_{{\cal T};2}=
\ga_{\!\!\bigoplus\limits_{h\in\chi({\cal T})}\!\!{\cal F}_h{\cal T}}
\oplus\bigoplus_{h\not\in\chi({\cal T})}F_h{\cal T},\qquad
\tilde{\cal F}{\cal S}_{{\cal T};2}=
\ga_{\!\!\bigoplus\limits_{h\in\chi({\cal T})}\!\!{\cal F}_h{\cal T}}\otimes
\begin{cases}
\Bbb{C},&\hbox{if~}\io_h\!=\!\hat{0}~\forall h\!\in\!\chi({\cal T});\\
{\cal F}_{\hat{1}}{\cal T}^{\otimes2}
,&\hbox{if~}\io_h\!=\!\hat{1}\!\neq\!\hat{0}~\forall h\!\in\!\chi({\cal T}).
\end{cases}$$
Let $\rho_{{\cal T};2}\!:{\cal FS}_{{\cal T};2}\!\lra\!
\tilde{\cal F}{\cal S}_{{\cal T};2}$ 
be the projection map followed by multiplication.
Define
\begin{gather*}
\al_{{\cal T};2}\!\in\!\Ga\big(\bar{\cal S}_{{\cal T};2}(\mu);
\hbox{Hom}(\tilde{\cal F}{\cal S}_{{\cal T};2};
L_{{\cal T},\hat{0}}^{*\otimes2}\otimes\ev^*T\PP)\big) 
\qquad\hbox{by}\\
\al_{{\cal T};2}\big((\tilde{\ups}_h)_{h\in\chi({\cal T})}\big)=
\sum_{h\in\chi({\cal T})}x_h{\cal D}_{{\cal T},h}^{(1)}\tilde{\ups}_h
\qquad\hbox{if~~}
b_{\ups}\!=\!\big(S^2,M,I;x,(j,y),u\big).
\end{gather*}
If $\chi_{\cal T}\hat{1}\!=\!0$ 
for some $\hat{1}\!\in\!\hat{I}\!-\!\chi({\cal T})$,
$H_{\hat{0}}{\cal T}\!=\!\{\hat{1},h_1\}$, and 
$h_2\!\in\!\chi({\cal T})\!-\!\{h_1\}$, let
$${\cal FS}_{{\cal T};2}={\cal F}_{\hat{1}}{\cal T}\!\oplus
{\cal F}_{h_2}{\cal T},\qquad
\tilde{\cal F}{\cal S}_{{\cal T};2}=
{\cal F}_{\hat{1}}{\cal T}\!\otimes\!{\cal F}_{h_2}{\cal T},
\quad\hbox{and}\quad
\rho_{{\cal T};2}\!:{\cal FS}_{{\cal T};2}\lra
\tilde{\cal F}{\cal S}_{{\cal T};2}$$
be the multiplication map.
Define  \hbox{$\al_{{\cal T};2}\!\in\!\Ga({\cal S}_{{\cal T};2}(\mu);
\hbox{Hom}(\tilde{\cal F}{\cal S}_{{\cal T};2};
L_{{\cal T},\hat{0}}^{*\otimes2}\otimes\ev^*T\PP)\big)$ by}
\begin{gather*}
\al_{{\cal T};2}\big(\ups_{\hat{1}}\otimes\ups_{h_2}\big)
=x_{\hat{1}}\ups_{\hat{1}}\!\!
\sum_{h\in\chi({\cal T})-\{h_1\}}\!\!{\cal D}_{{\cal T},h}\ups_h'
+x_{h_1}{\cal D}_{{\cal T},h_1}\ups_{h_1},\\
\hbox{if}\qquad
b_{\ups_{\hat{1}}\otimes\ups_{h_2}}\!=\!\big(S^2,M,I;x,(j,y),u\big),~~
\big[\ups_{\hat{1}}\otimes\ups',\ups_{h_1}\big]
\!\in\!{\cal S}_{{\cal T};2}(\mu),
\quad\hbox{where}\quad
\ups'=(\ups_h)_{h\in\chi({\cal T})-\{h_1\}}.
\end{gather*}
In all cases, denote by ${\cal F}^{\eset}{\cal S}_{{\cal T};k}$
the subset of ${\cal FS}_{{\cal T};k}$ consisting of vectors
with all components nonzero.

\begin{prp}
\label{n2cusps_str}
If ${\cal T}^*\!=\!(S^2,[N]\!-\!M_0,\{\hat{0}\};\hat{0},d)$
and ${\cal T}\!=\!(S^2,[N]\!-\!M_0,I;j,\under{d})\!<\!{\cal T}^*$
are simple bubble types,
$$\bar{\cal S}_{{\cal T}^*;1}(\mu)\cap{\cal U}_{\cal T}(\mu)=
\begin{cases}
{\cal S}_{{\cal T};1}(\mu),&\hbox{if~}|\chi({\cal T})|\!=\!1;\\
{\cal S}_{{\cal T};1}(\mu)\cup{\cal S}_{{\cal T};2}(\mu)
,&\hbox{if~}|H_{\hat{0}}{\cal T})|\!=\!|\hat{I}|\!=\!2~\&~
M_{\hat{0}}{\cal T}\!=\!\eset;\\
{\cal S}_{{\cal T};2}(\mu),&\hbox{otherwise}.
\end{cases}$$
In addition, there exist 
$\de,C\!\in\! C^{\i}({\cal S}_{{\cal T};k}(\mu);\Bbb{R}^+)$
and a continuous~map
$$\ga_{{\cal T};k}\!:
 {\cal FS}_{{\cal T};k;\de}\lra\bar{\cal S}_{{\cal T}^*;1}(\mu)$$
onto an open neighborhood of ${\cal S}_{{\cal T};k}(\mu)$
in $\bar{\cal S}_{{\cal T}^*;1}(\mu)$ such that 
$\ga_{{\cal T};k}|{\cal S}_{{\cal T};k}(\mu)$ is the identity
and $\ga_{{\cal T};k}|{\cal F}^{\eset}{\cal S}_{{\cal T};k;\de}$
is an orientation-preserving diffeomorphism onto
an open subset of~${\cal S}_{{\cal T}^*;1}(\mu)$.
Furthermore, if $d_{\hat{0}}\!\neq\!0$, 
${\cal D}_{{\cal T}^*,\hat{0}}^{(2)}$
does not vanish on ${\cal S}_{{\cal T};1}(\mu)$.
If $d_{\hat{0}}\!=\!0$, with appropriate identifications,
$$\Big|{\cal D}_{{\cal T}^*,\hat{0}}^{(2)}\big(\ga_{{\cal T};k}(\ups)\big)
-\al_{{\cal T};k}\big(\rho_{{\cal T};k}(\ups)\big)\Big|
\le C(b_{\ups})|\ups|^{\frac{1}{p}}\big|\rho_{{\cal T};k}(\ups)\big|
\qquad\forall\ups\!\in\!{\cal FS}_{{\cal T};k;\de}.$$
\end{prp}

\begin{lmm}
\label{n2cusps_lmm}
If ${\cal T}^*$ and ${\cal T}$ are as in Proposition~\ref{n2cusps_str},
there exist $\de\!\in\! C^{\i}({\cal S}_{{\cal T};k}(\mu);\Bbb{R}^+)$
and a continuous~map
$$\tilde{\ga}_{{\cal T};k}\!: 
\big({\cal NS}_{{\cal T};k}\oplus {\cal FT}\big)_{\de} \lra
\bar{\cal U}_{{\cal T}^*}(\mu)$$
onto an open neighborhood of ${\cal S}_{{\cal T};k}(\mu)$
in $\bar{\cal U}_{{\cal T}^*}(\mu)$ such that
$\ga_{{\cal T};k}|{\cal S}_{{\cal T};k}(\mu)$ is the identity
and $\ga_{{\cal T};k}$ is smooth and orientation-preserving
on  the preimage of~${\cal V}_{{\cal T}^*}(\mu)$.
Furthermore, with appropriate identifications,
$${\cal D}_{{\cal T}^*,\hat{0}}\tilde{\ga}_{{\cal T};k}(\ups)
=\al_{{\cal T};k}^{(1)}(\ups)
\qquad\forall\ups\!\in\!
\big({\cal NS}_{{\cal T};k}\oplus {\cal FT}\big)_{\de},$$
where $\al_{{\cal T};k}^{(1)}\!: 
{\cal NS}_{{\cal T};k}\!\oplus {\cal FT}\!\lra\!
L_{{\cal T},\hat{0}}^*\otimes\ev^*T\PP$
is a degenerate polynomial of constant rank.
\end{lmm}

\noindent
The proof of Proposition~\ref{n2cusps_str},
with the exception of the estimate on~${\cal D}_{{\cal T}^*,\hat{0}}^{(2)}$, 
is either the same or very similar to 
the proof of Lemma~\ref{enum1-cuspcurves_l1} in~\cite{Z2}, 
depending on the bubble type~${\cal T}$.
If $d_{\hat{0}}\!\neq\!0$, we apply 
the analytic estimate of Theorem~\ref{str_global}
and the Implicit Function Theorem to 
the section~$\ga_{\cal T}^{\mu*}{\cal D}_{{\cal T}^*,\hat{0}}$.
If $d_{\hat{0}}\!=\!0$ and $k\!=\!1$,
the two theorems are applied to a section of 
\hbox{$(L_{{\cal T},\hat{0}}\otimes{\cal F}_{h_1}{\cal T})^*
 \!\otimes\ev^*T\PP$}
induced by~$\ga_{\cal T}^{\mu*}{\cal D}_{{\cal T}^*,\hat{0}}$.
If $|\chi({\cal T})|\!\ge\!2$ and $H_{\hat{0}}{\cal T}\!=\!\chi({\cal T})$,
 we work with a section 
of \hbox{$(L_{{\cal T},\hat{0}}\otimes\ga_{E{\cal T}})^*\!\otimes\ev^*T\PP$}
over the blowup of~$E{\cal T}$ along~${\cal U}_{\cal T}(\mu)$.
The case  $\chi_{\cal T}h\!=\!\hat{1}$
for all $h\!\in\!\chi({\cal T})$ is similar.
If $\chi({\cal T})\!=\!\{h_1,h_2\}$ is a two-element set,
and $H_{\hat{0}}{\cal T}\!=\!\{\hat{1},h_1\}$,
we use the same section, but given a small element 
\hbox{$(\ups_{\hat{1}},\ups_{h_2})\!\in\!{\cal FS}_{{\cal T};2}$},
we start with the approximate solution 
$(\ups_{\hat{1}},\ka\ups_{\hat{1}}\ups_{h_2},\ups_{h_2})$, with
$$\ka\!\in\!(L_{{\cal T},\hat{1}}\otimes L_{{\cal T},h_2})^*
\!\otimes L_{{\cal T},h_1}
\quad\hbox{s.t.}\quad
\big[\ups_{\hat{1}}\ups_{h_2},\ka\ups_{\hat{1}}\ups_{h_2}\big]
\!\in\!{\cal S}_{{\cal T};2}(\mu).$$
The approach to the remaining case is analogous.
The estimate on ${\cal D}_{{\cal T}^*,\hat{0}}^{(2)}$
is obtained by the same argument as in the proof of 
Lemma~\ref{enum1-cuspcurves_l2} in~\cite{Z2}.
The proof makes use of the construction of~$\ga_{\cal T}^{\mu}$
in~\cite{Z3}, which involves a modification of
the pregluing step of standard gluing procedures
for pseudoholomorphic curves.
Finally, Lemma~\ref{n2cusps_lmm} is proved similarly
to Corollary~\ref{enum1-cuspcurves_l1c} in~\cite{Z2}.\\

\noindent
We next describe the behavior of the section
$${\cal D}_{2;2}\equiv c_1{\cal D}_1\!+c_2{\cal D}_2\!
\in\Ga\big(\bar{\cal S}_{2;2}(\mu);
\ga_{E_2}^*\otimes\ev^*T\PP\big),$$
for $c_1,c_2\!\in\!\Bbb{C}$ distinct, near 
${\partial}\bar{\cal S}_{2;2}(\mu)\!\equiv
\bar{\cal S}_{2;2}(\mu)\!-\!{\cal S}_{2;2}(\mu)$.
As before, we identify ${\cal S}_{2;1}(\mu)$ 
with a subset of~$\bar{\cal S}_{2;2}(\mu)$.
Similarly, if
${\cal T}\!=\!\big(S^2,[N],I;j,\under{d})$ is a bubble type
such that $I\!-\!\hat{I}\!=\!\{k_1,k_2\}$ is a two-element set
and $\sum d_i\!=\!d$, let
$${\cal S}_{{\cal T};2}(\mu)\!=\!
\big\{\big[b,L_{{\cal T},k_1}\big]\!: 
b\!\in\!{\cal U}_{\cal T}(\mu),~{\cal D}_{{\cal T},k_1}b\!=\!0\big\}
\cup
\big\{\big[b,L_{{\cal T},k_2}\big]\!: 
b\!\in\!{\cal U}_{\cal T}(\mu),~{\cal D}_{{\cal T},k_2}b\!=\!0\big\}
\subset\Bbb{P}E_2.$$

\begin{prp}
\label{n2tacnodes_str}
Suppose $d$ is a positive integer and
$\mu$ is a tuple of $3d\!-\!4$ points in general position in~$\PP$.
Then
$$\partial\bar{\cal S}_{2;2}(\mu)
={\cal S}_{2;1}(\mu)\cup\bigcup_{[{\cal T}]}{\cal S}_{{\cal T};2}(\mu),$$
where the union is taken over all equivalence classes of non-basic types
${\cal T}\!=\!\big(S^2,[N],I;j,\under{d})$ such that
$I\!-\!\hat{I}\!=\!\{k_1,k_2\}$ is a two-element set 
\hbox{and $\sum d_i\!=\!d$}.
Furthermore, there exist $\de,C\!>\!0$ and homeomorphism
$$\ga_{2;2}\!:
\big\{u\!\in\!\ga_{E_2}^*\otimes(E_2/\ga_{E_2})\big|
\partial\bar{\cal S}_{2;2}(\mu)\!: |u|<\de\big\}\lra\bar{\cal S}_{2;2}(\mu)$$
onto an open neighborhood of $\partial\bar{\cal S}_{2;2}(\mu)$
in $\bar{\cal S}_{2;2}(\mu)$ such that 
$\ga_{2;2}|\partial\bar{\cal S}_{2;2}(\mu)$ is the identity and
$\ga_{2;2}$ restricts to an orientation-preserving diffeomorphism
from the complement of $\partial\bar{\cal S}_{2;2}(\mu)$
onto an open subset of~${\cal S}_{2;2}(\mu)$.
Finally, with appropriate identifications,
$$\big|{\cal D}_{2;2}\ga_{2;2}(u)-\al_{2;2}(u)\big|
\le C|u|^{1+\frac{1}{p}}  \qquad\forall 
u\!\in\! \big\{u\!\in\!\ga_{E_2}^*\otimes(E_2/\ga_{E_2})\big|
\partial\bar{\cal S}_{2;2}(\mu)\!: |u|<\de\big\}\lra\bar{\cal S}_{2;2}(\mu),$$
where $\al_{2;2}\!\in\!\Ga\big(\partial\bar{\cal S}_{2;2}(\mu);
\hbox{Hom}(\ga_{E_2}^*\otimes(E_2/\ga_{E_2}),
\ga_{E_2}^*\otimes\ev^*T\PP)\big)$ is an injection on every fiber.
\end{prp}

\noindent
This proposition follows from Theorem~\ref{str_global}
and the Implicit Function Theorem by an argument similar
to the proof of Proposition~\ref{n2cusps_str}. 
Note that with our choice of constraints, 
$\partial\bar{\cal S}_{2;2}(\mu)$
is a finite set. Thus, we are able to take $\de$ and $C$ to be positive
real numbers rather than continuous functions 
$\partial\bar{\cal S}_{2;2}(\mu)\!\lra\!\Bbb{R}^+$.

\subsection{Description of $CR_3(\mu)$}
\label{str_thm_sec}

\noindent
In this subsection, we describe the number of elements
of ${\cal M}_{\Si,t\nu,d}(\mu)$ that lie near each 
strata ${\cal M}_{\cal T}(\mu)$ of bubble maps of type (2c)
in terms of the zeros of affine maps between vector bundles
over closures of certain subspaces of~${\cal M}_{\cal T}(\mu)$.
These results are proved by an argument similar to
Sections~\ref{enum1-analysis_sec} and~\ref{enum1-resolvent_sec}
in~\cite{Z2}, which is outlined briefly at the end of 
this subsection.\\

\noindent
We start by recalling more notation used in~\cite{Z2}.
If ${\cal T}\!=\!(\Si,[N],I;j,\under{d})$ is a simple bubble type,
let $\big\{g_{b,\hat{0}}\!:b\!\in\!{\cal M}_{\cal T}(\mu)\big\}$ be
a smooth family of Kahler metrics on $(\Si,j)$ such that\\
(1) for all $b\!=\!(\Si,[N],I;x,(j,y),u)\!\in\!{\cal M}_{\cal T}(\mu)$
and $h\!\in\!H_{\hat{0}}{\cal T}$,
$g_{b,\hat{0}}$ is flat on a neighborhood of $x_h$~in~$\Si$;\\
(2) the metric $g_{b,\hat{0}}$ is determined by the tuple
$(x_h)_{h\in H_{\hat{0}}{\cal T}}$.\\
We denote by ${\cal H}_{\Si}^{0,1}$ the 3-dimensional space
of harmonic $(0,1)$-forms on~$\Si$.\\

\noindent
If $\psi\!\in\!{\cal H}_{\Si}^{0,1}$,
$b\!\in\!{\cal M}_{\cal T}$, $m\!\ge\! 1$, and
the metric $g_{b,\hat{0}}$ is flat near~$x$, 
we define  $D_{b,x}^{(m)}\psi\!\in\! T_x^{0,1}\Si^{\otimes m}$
as follows. 
If $(s,t)$ are conformal coordinates centered at $x$
such that $s^2\!+\!t^2$  is the square of the $g_{b,\hat{0}}$-distance 
to $x$, let  
$$\{D_{b,x}^{(m)}\psi\}\Big(\frac{\partial}{\partial s}\Big)\equiv
\{D_{b,x}^{(m)}\psi\}\bigg(
\underset{m}{\underbrace{\frac{\partial}{\partial s},\ldots,
\frac{\partial}{\partial s}}}\bigg)=\frac{\pi}{m!}
\Big\{\frac{D^{m-1}}{ds^{m-1}}\psi_j\Big|_{(s,t)=0}\Big\}
\Big(\frac{\partial}{\partial s}\Big),$$
where the covariant derivatives are taken with respect to the 
metric~$g_{b,\hat{0}}$. 
If $\{\psi_j\}$ is an orthonormal basis for ${\cal H}_{\Si}^{0,1}$,
let $s_{b,x}^{(m)}\!\in\! T_x^*\Si^{\otimes m}\otimes{\cal H}_{\Si}^{0,1}$
be given~by
$$s_{b,x}^{(m)}(v)\equiv
s_{b,x}^{(m)}(\underset{m}{\underbrace{v,\ldots,v}})=
\sum\ov{\Big\{D_{b,x}^{(m)}\psi_j\Big\}(v)}\psi_j.$$
The section $s_{b,x}^{(m)}$ is always independent of the choice
of  a basis for ${\cal H}_{\Si}^{0,1}$, but is dependent
on the choice of the metric $g_{b,\hat{0}}$ \hbox{if $m\!>\!1$.}
By \cite[p246]{GH}, $s_x\!\equiv\!s_{b,x}$  does not vanish and thus 
determines a line subbundle ${\cal H}_{\Si}^+$ of
$\Si\!\times\!{\cal H}_{\Si}^{0,1}\!\lra\!\Si$.
We denote its orthogonal complement by ${\cal H}_{\Si}^-$.
Let 
$$\pi^-\in\Ga\big(\Si;(\Si\times{\cal H}_{\Si}^{0,1})^*\otimes
{\cal H}_{\Si}^-\big)$$
be the orthogonal projection map onto~${\cal H}_{\Si}^-$.
While the section~$s_{b,x}$ depends on the choice of
the metric~$g_{b,\hat{0}}$, 
$s_x^{(2)}\!\equiv\!\pi_x^-\circ s_{b,x}^{(2)}$ does not
and thus is globally defined on~$\Si$.
If $\Si$ is not hyperelleptic, as we assume to be the case,
$s_x^{(2)}$ does not vanish and thus determines 
a line subbundle ${\cal H}_{\Si}^{-+}$ of~${\cal H}_{\Si}^-$.
We denote its orthogonal complement by ${\cal H}_{\Si}^{- -}$
and the corresponding orthogonal projection map by~$\pi^{- -}$.
The composition $s_x^{(3)}\!\equiv\!\pi_x^{- -}\circ s_{b,x}^{(3)}$
is again independent of the choice of the metric~$g_{b,\hat{0}}$.
If $\Si$ is generic, the section 
$$s^{(3)}\!\in\!\Ga
\big(\Si;T^*\Si^{\otimes3}\otimes{\cal H}_{\Si}^{- -}\big)$$
vanishes transversally at $24$ distinct points $z_1,\ldots,z_{24}$ of~$\Si$.
These points correspond to the flexes of~$\Si$ under
the canonical embedding into~$\PP$.

\begin{thm}
\label{si_str}
Suppose $d$ is a positive integer, $N\!=\!3d\!-\!4$, 
$\mu$ is an $N$-tuple of points in general position in~$\PP$,
${\cal T}\!=\!(\Si,[N],I;j,\under{d})$ is a simple bubble type
such $d_{\hat{0}}\!=\!0$ and $\sum d_i\!=\!d$.
If 
$$\nu\in\Ga(\Si\!\times\!\PP;
            \La^{0,1}\pi_{\Si}^*T^*\Si\otimes\pi_{\PP}^*T\PP)$$
is a generic section,
there exist a compact subset $K_{{\cal T},\nu}$ of 
${\cal M}_{\cal T}(\mu)$ and integer~$N({\cal T})$
with the following property.
If $K$ is a compact subset of ${\cal M}_{\cal T}(\mu)$
containing~$K_{{\cal T},\nu}$,
there exist a neighborhood $U_K$ of $K$ in $\bar{C}^{\i}_{(d;[N])}(\Si;\mu)$
and $\ep_K\!>\!0$ such that for all $t\!\in\!(0,\ep_K)$,
$$^{\pm}\big|U_K\cap{\cal M}_{\Si,d,t\nu}(\mu)\big|
=N({\cal T}).$$
If ${\cal T}$ is not primitive, 
$U_K\cap{\cal M}_{\Si,d,t\nu}(\mu)\!=\!\eset$.
If ${\cal T}$ is primitive, $N({\cal T})$ is the numbers
of zeros of the affine maps between vector bundles described below.
\end{thm}

\noindent
Above $\bar{C}^{\i}_{(d;[N])}(\Si;\mu)$ denotes the space
of all stable maps from $\Si$ to~$\PP$ that map the marked points
to $\mu_1,\ldots,\mu_N$.
For each primitive bubble type~${\cal T}$,
we now describe the number $N({\cal T})$ as 
the sum of numbers $N(\al)$, where each $\al$ is a regular ms-polynomial
between two ms-bundles over an ms-manifold; see Subsection~\ref{top_sec}.\\

\noindent
If $|\hat{I}|\!>\!3$, ${\cal M}_{\cal T}(\mu)\!=\!\eset$.
If $|\hat{I}|\!=\!1,2,3$, define
\begin{gather}
\al_{|\hat{I}|}\!\in\!\Ga\big(\Si^{\hat{I}}\!\times
                           \!\bar{\cal U}_{\bar{\cal T}}(\mu);
\hbox{Hom}(\bigoplus_{h\in\hat{I}}T\Si_h\otimes L_{{\cal T},h},
{\cal H}_{\Si}^{0,1}\otimes\ev^*T\PP)\big),\notag\\
\label{k1contr}
\hbox{by}\qquad
\al_{|\hat{I}|}\big(b,(x_h)_{h\in\hat{I}};
(v_h\otimes\ups_h)_{h\in\hat{I}}\big)=\sum_{h\in\hat{I}}
\big({\cal D}_{{\cal T},h}\ups_h\big)\big(s_{x_h}v_h\big).
\end{gather}
If $|\hat{I}|\!=\!3$, $N({\cal T})\!=\!n^{(1)}({\cal T})\!=\!N(\al_3)$.
If $|\hat{I}|\!=\!2$, 
$N({\cal T})\!=\!n^{(1)}({\cal T})\!+\!2n^{(2)}({\cal T})
\!=\!N(\al_2)\!+\!2N(\al_{2;1})$,
with $\al_{2;1}$ defined as follows.
If $\hat{I}\!=\!\{\hat{1},\hat{2}\}$, 
$b\!\in\!{\cal S}_{\bar{\cal T};1}(\mu)$, and 
${\cal D}_{\bar{\cal T};\hat{1}}(b)\!=\!0$,
\begin{gather}
\al_{2;1}\big(b,(x_h)_{h\in\hat{I}};
(v_h\otimes\ups_h)_{h\in\hat{I}}\big)=
\big({\cal D}_{{\cal T},\hat{1}}^{(2)}\ups_{\hat{1}}\big)
\big(s_{x_{\hat{1}}}^{(2)}v_{\hat{1}})+
\big({\cal D}_{{\cal T},\hat{2}}\ups_{\hat{2}}\big)
\big(\pi_{x_{\hat{1}}}^-s_{x_{\hat{2}}}v_{\hat{2}}\big)
\in {\cal H}_{\Si}^-(x_{\hat{1}})\otimes T_{\ev(b)}\PP\notag,\\
\label{m2k2_e}
\hbox{if}\qquad
v_{\hat{1}}\otimes \ups_{\hat{1}}\!\in\!
T\Si_{\hat{1}}^{\otimes2}\otimes L_{{\cal T},\hat{1}}^{\otimes2}
\big|_{(b,x_{\hat{1}})},\quad
v_{\hat{2}}\otimes \ups_{\hat{2}}\!\in\!
T\Si_{\hat{2}}\otimes L_{{\cal T},\hat{2}}\big|_{(b,x_{\hat{2}})}.
\end{gather}
If $|\hat{I}|\!=\!1$,
\begin{equation}\label{m1_sum}
N({\cal T})=
n^{(1)}({\cal T})\!+\!2n^{(2)}({\cal T})\!+\!3n_3({\cal T})\!+\!4n_4({\cal T})
=N(\al_1)\!+\!2N(\al_{1;1})\!+\!3N(\al_{1;2})\!+\!96|{\cal S}_{1;2}(\mu)|,
\end{equation}
where 
\begin{gather*}
\al_{1;1}\!\in\!\Ga\big(\Si\!\times\!\bar{\cal S}_1(\mu);
\hbox{Hom}(T\Si^{\otimes2}\otimes L_{{\cal T},\hat{1}}^{\otimes2},
{\cal H}_{\Si}^-\otimes\ev^*T\PP)\big)\quad\hbox{and}\\
\al_{1;2}\!\in\!\Ga\big(\Si\!\times\!{\cal S}_{1;2}(\mu);
\hbox{Hom}(T\Si^{\otimes3}\otimes L_{{\cal T},\hat{1}}^{\otimes3},
{\cal H}_{\Si}^{- -}\otimes\ev^*T\PP)\big)
\end{gather*}
are defined by
\begin{equation}\label{m1k2and3}
\al_{1;k}\big(x,b,v\otimes\ups)=
\big({\cal D}_{{\cal T},\hat{1}}^{(k+1)}\ups\big)\big(s_x^{(k+1)}v\big).
\end{equation}
Finally, for each $m\!=\!1,2,3$, we denote by $n^{(k)}_m(\mu)$
the sum of the numbers $n^{(k)}({\cal T})$ taken over all
equivalence classes of primitive bubble types ${\cal T}$
\hbox{with $|\hat{I}|\!=\!m$}.

\begin{crl}
\label{CR3_str}
With notation as above,
$$CR_3(\mu)\!=\!
\big(n_1^{(1)}(\mu)\!+\!2n_1^{(2)}(\mu)\!+\!3n_1^{(3)}(\mu)\!
+\!96|{\cal S}_{1;2}(\mu)|\big)
+\big(n_2^{(1)}(\mu)\!+\!2n_2^{(2)}(\mu)\big)+n_3^{(1)}(\mu).$$
\end{crl}

\noindent
Corollary~\ref{CR3_str} follows immediately from 
the preceding paragraph, Theorem~\ref{si_str}, and 
the definition of $CR_3(\mu)$ in Subsection~\ref{summary}.\\

\noindent
{\it Remark~1:} The multiplicity $k$ for $n_m^{(k)}(\mu)$
is the degree of a polynomial map between two vector spaces
of small dimensions.
In the cases under consideration, these degrees are clear.
In the genus-two case, they are described in 
Section~\ref{enum1-resolvent_sec} of~\cite{Z2}.\\

\noindent
{\it Remark~2:}
The number $96$ in~\e_ref{m1_sum} arises
because each element of ${\cal M}_{\cal T}(\mu)$
corresponding to a simple flex of~$\Si$
and a rational map in ${\cal S}_{2;2}(\mu)$ enters with 
a multiplicity of~$4$.
A hyperflex of $\Si$ would result in a multiplicity of~$10$,
at least \hbox{if $d\!\ge\!4$}.
Thus, if $\Si$ has $n$~hyperflexes and $24\!-\!2n$ simple flexes,
the number~$96$ in equation~\e_ref{m1_sum} and in Corollary~\ref{CR3_str}
should be replaced \hbox{by $96\!+\!2n$}.
No other changes are needed.
The analogue of the term $n_1^{(4)}(\mu)\!=\!24|{\cal S}_{1;2}(\mu)$
in the genus-two case is $n_1^{(3)}(\mu)\!=\!6|{\cal S}_1(\mu)|$,
as each of the six hyperelleptic points of~$\Si$ and a cuspidal map
through the fixed $3d\!-\!2$ points enters with a multiplicity of~$3$;
see Subsection~\ref{enum1-order3_case1} in~\cite{Z2}.\\

\noindent
{\it Remark~3:} It should be possible to adapt this approach
to the case $\Si$ is a hyperelleptic genus-three surface,
but significant changes would be required.
In particular, there will likely be a contribution to $CR_3(\mu)$
from an affine map over the space $\Si\!\times\!\bar{\cal S}_{2;2}(\mu)$,
as is the case in the genus-two case in~$\Bbb{P}^3$;
see Subsection~\ref{enum1-order2_case2a} in~\cite{Z2}.
Furthermore, higher-order contributions $n_m^{(k)}(\mu)$, $k\!\ge\!3$,
will have a very different description, which will involve
the hyperelleptic and Weierstrass points of~$\Si$.\\

\noindent
We now outline the proof of Theorem~\ref{si_str} following~\cite{Z2}.
For each $b\!\in\!{\cal M}_{\cal T}$ and 
$\ups\!\in\!F^{\eset}{\cal T}$ sufficiently small, 
we first construct a  nearly holomorphic map 
$u_{\ups}\!:\Si\!\lra\!\P$.
We then attempt to solve the equation
\begin{equation}\label{pert_e}
\bar{\partial}\exp_{u_{\ups}}\xi=t\nu\Llra
\bar{\partial}u_{\ups}+D_{\ups}\xi+N_{\ups}\xi=t\nu
\in\Ga^{0,1}(u_{\ups})
\end{equation}
for a small vector $\xi\!\in\!\Ga(u_{\ups})$ along $u_{\ups}$.
Since we need to count the number of elements of
${\cal M}_{\Si,t\nu,d}(\mu)$ that lie nearly~${\cal M}_{\cal T}(\mu)$,
we require that $\xi$ lie in a subspace $\tilde{\Ga}_+(\ups)$
complementary to the ``tangent bundle'' $\Ga_-(\ups)$
of the space~$\big\{u_{\ups}\!:\ups\!\in\!F^{\eset}{\cal T}_{\de}\big\}$.
The cokernel of the operator $D_b$ is 
${\cal H}_{\Si}^{0,1}\!\otimes T_{\ev(b)}\PP$.
It induces an orthogonal splitting 
$$\Ga^{0,1}(\ups)=\Ga_-^{0,1}(\ups)\oplus \Ga_+^{0,1}(\ups)$$
such that $\pi^{0,1}_+\!\circ D_{\ups}\!:\tilde{\Ga}_+(\ups)
\!\lra\!\Ga_+^{0,1}(\ups)$ is an isomorphism and
$\Ga_-^{0,1}(\ups)$ is isomorphic to 
${\cal H}_{\Si}^{0,1}\!\otimes T_{\ev(b)}\PP$.
If $\ups$ and $t$ are sufficiently small, depending on~$b$,
the $\pi^{0,1}_+$-part of equation~\e_ref{pert_e} has 
a unique solution~$\xi_{\ups}$, which also solves 
the entire equation~\e_ref{pert_e} if and only~if
\begin{equation}\label{pert_e2}
\pi^{0,1}_-\big(t\nu-
\bar{\partial}u_{\ups}-D_{\ups}\xi_{\ups}-N_{\ups}\xi_{\ups}\big)
=0\in\Ga^{0,1}_-(u_{\ups})\approx 
{\cal H}_{\Si}^{0,1}\!\otimes T_{\ev(b)}\PP.
\end{equation}
It then remains to adjust for the constraints and 
extract the leading-order term from~\e_ref{pert_e2}.
The latter part depends on the choices of the above splittings
of $\Ga(\ups)$ and~$\Ga^{0,1}(\ups)$.
The spaces $\Ga_-(\ups)$ and $\Ga_-^{0,1}(\ups)$ are constructed 
from the kernel and cokernel of $D_b$ fairly explicitly
in Subsections~\ref{enum1-obs_sec} and~\ref{enum1-tangent_sec} of~\cite{Z2}.
In order to extract the leading-order term from~\e_ref{pert_e2},
we need the composite $\pi^{0,1}_-\!\circ\!D_{\ups}$ to 
be sufficiently small on~$\tilde{\Ga}_+(\ups)$.
In~\cite{Z2}, this is insured by choosing $\tilde{\Ga}_+(\ups)$
so that its image under~$D_{\ups}$ is orthogonal to
${\cal H}_{\Si}^+(x_{\hat{1}})\otimes\ev^*T\PP$ if
$b\!=\!(S^2,[N],I;x,(j,y),u)$ and $\hat{1}$ is an element 
of~$H_{\hat{0}}{\cal T}$.
By an argument similar to Subsection~\ref{enum1-tangent_sec} in~\cite{Z2},
in the given case we can  choose $\tilde{\Ga}_+(\ups)$
so that its image under~$D_{\ups}$ is orthogonal to
$\big({\cal H}_{\Si}^+(x_{\hat{1}})\oplus
{\cal H}_{\Si}^{-+}(x_{\hat{1}})\big)\otimes\ev^*T\PP$,
provided $d_{\hat{1}}\!\ge\!2$.
Then in all cases, $\pi^{0,1}_-\!\circ\!D_{\ups}|\tilde{\Ga}_+(\ups)$ 
will be sufficiently small for the purposes of extracting
dominant terms from~\e_ref{pert_e2} as in 
Section~\ref{enum1-resolvent_sec} of~\cite{Z2}.
Polynomial maps between vector bundles arise from 
the power series expansion for~$\bar{\partial}u_{\ups}$
given in Proposition~\ref{enum1-dbar_gen2prp} of~\cite{Z2}.

\section{Rational Curves with Singularities}
\label{rational_sect}

\subsection{Intersections in ${\cal S}_{2;1}(\mu)$,
$\bar{\cal S}_{2;2}(\mu)$, and $\bar{\cal S}_1(\mu)$}

\noindent
This section is dedicated to computing the intersection numbers
of spaces of stable rational maps that are needed in~Section~\ref{comp_sect}.
We start with ``codimension-one'' and ``-two'' cases.

\begin{lmm}
\label{n2_2cusps}
If $d$ is a positive integer,
the number of two-component rational degree-$d$ curves
passing through a tuple~$\mu$ of $3d\!-\!4$ points 
in general position in~$\PP$
such that the two components meet at a node at which one of them has a cusp
is given by
$$\big|{\cal S}_{2;1}(\mu)\big|=
\big\lan 6a^2+3a\big(c_1({\cal L}_1^*)\!+\!c_1({\cal L}_2^*)\big)+
\big(c_1^2({\cal L}_1^*)\!+\!c_1^2({\cal L}_2^*)\big),
\big[\bar{\cal V}_2(\mu)\big]\big\ran-3\tau_3(\mu),$$
where $\tau_3(\mu)=\big|{\cal V}_3(\mu)\big|$.
\end{lmm}

\begin{center}
\begin{tabular}{||c|c|c|c|c|c|c||}
\hline\hline
$d$&  2&3&4&5&6&7\\
\hline
$|{\cal S}_{2;1}(\mu)|$& 
0&0&528& 91,872& 26,055,360& 12,596,219,904\\ 
\hline\hline
\end{tabular}
\end{center}

\noindent
{\it Proof:} The proof is essentially the same as that of 
Lemma~\ref{enum1-n2cusps} in~\cite{Z2}, which enumerates 
irreducible cuspidal curves through $3d\!-\!2$ points.
The argument uses the analytic estimate of Theorem~\ref{str_global}
and the topological tools of Subsection~\ref{top_sec}.
In fact, the above formula can be deduced from 
the formula of Lemma~\ref{enum1-n2cusps} in~\cite{Z2},
since its proof applies with no change to enumerate
irreducible curves through $3d\!-\!3$ points with a cusp on
a fixed line in~$\PP$.

\begin{lmm}
\label{n2tacnodes}
If $d$ is a positive integer and $\mu$ is
a tuple of $3d\!-\!4$ points in general position in~$\PP$,
\begin{gather*}
\big\lan a,\big[\bar{\cal S}_{2;2}(\mu)\big]\big\ran
=\big\lan 3a^2+a\big(c_1({\cal L}_1^*)\!+\!c_1({\cal L}_2^*)\big),
\big[\bar{\cal V}_2(\mu)\big]\big\ran,\\
\big\lan \la_{E_2},\big[\bar{\cal S}_{2;2}(\mu)\big]\big\ran=
\big\lan 3a^2+3a\big(c_1({\cal L}_1^*)\!+\!c_1({\cal L}_2^*)\big)
+\big(c_1^2(L_1^*)\!+\!c_1^2(L_2^*)\big)+c_1(L_1^*)c_1(L_2^*),
\big[\bar{\cal V}_2(\mu)\big]\big\ran,
\end{gather*}
with chern classes $c_1^2(L_1^*)\!+\!c_1^2(L_2^*)$ and
$c_1(L_1^*)c_1(L_2^*)$ defined similarly to
$c_1^2({\cal L}_1^*)\!+\!c_1^2({\cal L}_2^*)$ and 
$c_1({\cal L}_1^*)c_1({\cal L}_2^*)$.
\end{lmm}

\noindent
{\it Proof:} We only sketch the argument, since
the proof is analogous to that of Lemma~\ref{enum1-n3cusps} in~\cite{Z2}.
Let ${\cal D}\!\in\!\Ga\big(\Bbb{P}E_2;\ga_{E_2}^*\!\otimes\ev^*T\P)$
be the section induced by the section
$${\cal D}_1+{\cal D}_2\in\Ga\big(
\bar{\cal V}_2(\mu);E_2^*\otimes\ev^*T\PP\big),$$
defined similarly to $c_1({\cal L}_1^*)\!+\!c_1({\cal L}_2^*)$;
see Subsection~\ref{notation_sec}. Then
$${\cal S}_{2;2}(\mu)={\cal D}^{-1}(0)\cap
\big(\Bbb{P}E_2|\big({\cal V}_2(\mu)-{\cal S}_{2;1}(\mu)\big)\big).$$
Let $s$ be a section of $\ev^*{\cal O}(1_{\PP})\!\lra\!\Bbb{P}E_2$ such that
$s$ is smooth and transversal to the zero set of the bundle
on all smooth strata of $\Bbb{P}E_2$ and of~$\bar{\cal S}_{2;2}(\mu)$.
In particular, $s^{-1}(0)\cap\partial\bar{\cal S}_{2;2}(\mu)\!=\!\eset$.
Then, by Proposition~\ref{euler_prp},
\begin{equation*}\begin{split}
\big\lan a,\big[\bar{\cal S}_{2;2}(\mu)\big]\big\ran
&=^{\pm}\!\big|s^{-1}(0)\cap{\cal S}_{2;2}(\mu)\big|
= ^{\pm}\!\big|{\cal D}^{-1}(0)\cap
\big(\Bbb{P}\big(E_2|{\cal V}_2(\mu)-{\cal S}_{2;1}(\mu)\big)\big)
\cap s^{-1}(0)\big|\\
&=\big\lan e\big(\ga_{E_2}^*\!\otimes\ev^*T\PP\big),
\big[s^{-1}(0)\big]\big\ran -{\cal C}_{\big(\Bbb{P}E_2|
\big(\bar{\partial}{\cal V}_2(\mu)-{\cal S}_{2;1}(\mu)\big)
\cap s^{-1}(0)}({\cal D}).
\end{split}\end{equation*}
By our assumptions on $s$, the last term above is zero.
Thus,
$$\big\lan a,\big[\bar{\cal S}_{2;2}(\mu)\big]\big\ran
=\big\lan ac_2\big(\ga_{E_2}^*\!\otimes\ev^*T\PP\big),
\big[\Bbb{P}E_2\big]\big\ran.$$
The first claim then follows from~\e_ref{zeros_main_e}
and from definitions of $E_2$ and of the relevant chern classes.
Note that
$$a\big(c_1(L_1^*)+c_1(L_2^*)\big)=
a\big(c_1({\cal L}_1^*)+c_1({\cal L}_2^*)\big).$$
The second identity is proved similarly.

\begin{lmm}
\label{n2cusps}
If $d$ is a positive integer and $\mu$ is
a tuple of $3d\!-\!4$ points in general position in~$\PP$,
\begin{gather*}
\big\lan a^2,\big[\bar{\cal S}_1(\mu)\big]\big\ran
=\big\lan a^2c_1^2({\cal L}^*),\big[\bar{\cal V}_1(\mu)\big]\big\ran
-\big\lan a^2,\big[\bar{\cal V}_2(\mu)\big]\big\ran,\\
\big\lan ac_1({\cal L}^*),\big[\bar{\cal S}_1(\mu)\big]\big\ran
=\big\lan 3a^2c_1^2({\cal L}^*)+ac_1^3({\cal L}^*),
\big[\bar{\cal V}_1(\mu)\big]\big\ran,\\
\big\lan c_1^2({\cal L}^*),\big[\bar{\cal S}_1(\mu)\big]\big\ran
=\big\lan 3a^2c_1^2({\cal L}^*)+3ac_1^3({\cal L}^*)+c_1^4({\cal L}^*),
\big[\bar{\cal V}_1(\mu)\big]\big\ran.
\end{gather*}
\end{lmm}

\noindent
{\it Proof:}
(1) This lemma is proved similarly to Lemma~\ref{n2tacnodes}.
Let 
$$E\!=\!\ev^*{\cal O}(1_{\PP})\oplus\ev^*{\cal O}(1_{\PP}),\qquad
E\!=\!\ev^*{\cal O}(1_{\PP})\oplus{\cal L}^*, 
\quad\hbox{or}\quad E\!=\!{\cal L}^*\oplus{\cal L}^*,$$ 
depending on which of the three identities we are proving.
Let $s$ be a section of $E\!\lra\!\bar{\cal V}_1(\mu)$
such that $s$ is smooth and transversal to the zero set in~$E$
on all smooth strata of $\bar{\cal V}_1(\mu)$ and 
of~$\bar{\cal S}_1(\mu)$. 
Then
\begin{equation}\label{n2cusps_e1}
\big\lan c_2(E),\big[\bar{\cal S}_1(\mu)\big]\big\ran
=\big\lan c_2(E)c_2\big(L^*\otimes\ev^*T\PP),
\big[\bar{\cal V}_1(\mu)\big]\big\ran-
{\cal C}_{\partial\bar{\cal V}_1(\mu)\cap s^{-1}(0)}({\cal D}).
\end{equation}
The bundle $L\!\lra\!\bar{\cal V}_1(\mu)$ and section 
${\cal D}\!\in\!\Ga\big(\bar{\cal V}_1(\mu);L^*\otimes\ev^*T\PP\big)$
in~\e_ref{n2cusps_e1} are defined as follows.
Let $N\!=\!3d\!-\!4$ and 
${\cal T}^*\!=\!(S^2,[N],\{\hat{0}\};\hat{0},d)$.
Then $L\!=\!L_{{\cal T}^*,\hat{0}}$ and 
${\cal D}\!=\!{\cal D}_{{\cal T}^*,\hat{0}}$.
Suppose 
$${\cal T}=(S^2,[N],I;j,\under{d})<{\cal T}^*.$$
If $d_{\hat{0}}\!\neq\!0$, ${\cal D}$ is transversal 
to the zero set on ${\cal U}_{\cal T}(\mu)$, and 
$s^{-1}(0)\cap{\cal D}^{-1}(0)\cap{\cal U}_{\cal T}(\mu)\!=\!\eset$
by our assumptions on~$s$.
Thus, from now on, we assume \hbox{that $d_{\hat{0}}\!=\!0$}.
If ${\cal T}$ is a semiprimitive bubble type,
either $\ev^*{\cal O}(1_{\PP})\big|\bar{\cal U}_{\cal T}(\mu)$
or ${\cal L}\big|\bar{\cal U}_{\cal T}(\mu)$ is trivial.
It follows that if $E\!=\!\ev^*{\cal O}(1_{\PP})\oplus{\cal L}^*$,
we can choose $s$ so that it does not vanish 
on~$\partial\bar{\cal V}_1(\mu)$.
The second equality is then immediate from~\e_ref{n2cusps_e1} and
\begin{equation}\label{n2cusps_e3}
ac_1(L^*)=ac_1({\cal L}^*)\in H^4\big(\bar{\cal V}_1(\mu)\big).
\end{equation}
(2) Suppose $E\!=\!\ev^*{\cal O}(1_{\PP})\oplus\ev^*{\cal O}(1_{\PP})$.
If the complex dimension of ${\cal U}_{\cal T}(\mu)$
is at least two and 
$\ev^*{\cal O}(1_{\PP})\big|\bar{\cal U}_{\cal T}(\mu)$ is not trivial,
$|\hat{I}|\!=\!\big|H_{\hat{0}}{\cal T}\big|\!=\!2$.
For a good choice of~$s$, the map $\ga_{\cal T}^{\mu}$ of 
Theorem~\ref{str_global} identifies neighborhoods 
of $s^{-1}(0)\cap{\cal U}_{\cal T}(\mu)$ in ${\cal FT}$ and in~$s^{-1}(0)$.
Since \hbox{$s^{-1}(0)\cap{\cal S}_{{\cal T};k}(\mu)\!=\!\eset$},
the section $\al_{\cal T}$ of Theorem~\ref{str_global}
defines an isomorphism between ${\cal FT}$ and $L^*\!\otimes\ev^*T\PP$
over every point of $s^{-1}(0)\cap{\cal U}_{\cal T}(\mu)$.
Thus, by Proposition~\ref{euler_prp} and the analytic estimate 
of Theorem~\ref{str_global},
$${\cal C}_{{\cal U}_{\cal T}(\mu)\cap s^{-1}(0)}({\cal D})
=^{\pm}\!\big|{\cal U}_{\cal T}(\mu)\cap s^{-1}(0)\big|
=\big\lan a^2,\big[\bar{\cal U}_{\cal T}(\mu)\big]\big\ran.$$
Summing these contributions over all equivalence classes
of such bubble types ${\cal T}$, we obtain
\begin{equation}\label{n2cusps_e5}
{\cal C}_{\partial\bar{\cal V}_1(\mu)\cap s^{-1}(0)}({\cal D})
=\big\lan a^2,\big[\bar{\cal V}_2(\mu)\big]\big\ran.
\end{equation}
The first identity follows from equations~\e_ref{n2cusps_e1},
\e_ref{n2cusps_e3}, and~\e_ref{n2cusps_e5}.\\
(3)  Suppose $E\!=\!{\cal L}^*\oplus{\cal L}^*$.
If the complex dimension of ${\cal U}_{\cal T}(\mu)$
is at least two and
${\cal L}\big|\bar{\cal U}_{\cal T}(\mu)$ is not trivial,
$H_{\hat{0}}{\cal T}\!=\!\{\hat{1}\}$ is a one-element set 
and $|\hat{I}|\!\in\!\{1,2\}$.
If $|\hat{I}|\!=\!2$, by an argument similar to (2) above,
${\cal U}_{\cal T}(\mu)\cap s^{-1}(0)$ is 
${\cal D}$-hollow in the sense of Definition~\ref{euler_dfn4},
where ${\cal D}$ is viewed as a section over 
$s^{-1}(0)\!\subset\!\bar{\cal V}_1(\mu)$.
Thus, by Proposition~\ref{euler_prp}, 
${\cal C}_{{\cal U}_{\cal T}(\mu)\cap s^{-1}(0)}({\cal D})\!=\!0$.
If $|\hat{I}|\!=\!1$, ${\cal T}\!=\!{\cal T}^*(\{l\})$
for some $l\!\in\![N]$.
For the purposes of computing 
${\cal C}_{{\cal U}_{\cal T}(\mu)\cap s^{-1}(0)}({\cal D})\!=\!0$,
it can be assumed that the map of Theorem~\ref{str_global} still identifies 
neighborhoods of $s^{-1}(0)\cap{\cal U}_{\cal T}(\mu)$ 
in ${\cal FT}$ and in~$s^{-1}(0)$; 
see the proof of Lemma~\ref{enum1-n3cusps} in~\cite{Z2}.
Since $\al_{\cal T}$ does not vanish on 
${\cal U}_{\cal T}(\mu)\cap s^{-1}(0)$, 
by Proposition~\ref{euler_prp} and the analytic estimate 
of Theorem~\ref{str_global},
$${\cal C}_{{\cal U}_{\cal T}(\mu)\cap s^{-1}(0)}({\cal D})
=\big\lan c_1\big(L^*\!\otimes\ev^*T\PP\big)-c_1({\cal FT}),
\big[\bar{\cal U}_{\cal T}(\mu)\cap s^{-1}(0)\big]\big\ran
=\big\lan c_1^2({\cal L}^*)c_1(L_{{\cal T},\hat{1}}^*),
\big[\bar{\cal U}_{\cal T}(\mu)\big]\big\ran,$$
since the restrictions of the bundles $L^*$ and $\ev^*T\PP$
to $\bar{\cal U}_{\cal T}(\mu)$ are trivial.
Summing up these contributions and using equation~\e_ref{normal_bundle2}, 
we obtain
\begin{equation}\label{n2cusps_e7}
{\cal C}_{\partial\bar{\cal V}_1(\mu)\cap s^{-1}(0)}({\cal D})
=\sum_{l\in[N]}
\big\lan c_1^3({\cal L}^*),
\big[\bar{\cal U}_{{\cal T}^*(l)}(\mu)\big]\big\ran.
\end{equation}
The final identity of the lemma follows from equations~\e_ref{n2cusps_e1},
\e_ref{n2cusps_e3}, and~\e_ref{n2cusps_e7};
see also equation~\e_ref{normal_bundle1}.

\subsection{Computation of $\big|{\cal S}_{1;2}(\mu)\big|$}

\noindent
In this subsection, we enumerate rational plane curves 
that have a $(3,4)$-cusp. 
We call point~$p$ a $(3,4)$-{\it cusp} of plane curve~${\cal C}$ 
if for a choice of local coordinates near~$p$, 
${\cal C}$ is parameterized by a map of the form
$$t\lra\big(t^3,t^4+o(t^4)\big),\qquad 0\lra p.$$

\begin{lmm}
\label{n2cusps2}
If $d$ is a positive integer,
the number of rational degree-$d$ curves that pass 
through a tuple~$\mu$ of $3d\!-\!4$ points in general position in~$\PP$
and have a $(3,4)$-cusp is given~by
\begin{equation*}\begin{split}
\big|{\cal S}_{1;2}(\mu)\big|=
\big\lan 3a^2+6ac_1({\cal L}^*)+4c_1^2({\cal L}^*),
\big[\bar{\cal S}_1(\mu)\big]\big\ran-2\big|{\cal S}_{2;1}(\mu)\big|
-3\tau_3(\mu)\\
-\big\lan 6a^2+3a\big(c_1({\cal L}_1^*)\!+\!c_1({\cal L}_2^*)\big)
+c_1({\cal L}_1^*)c_1({\cal L}_2^*),
\big[\bar{\cal V}_2(\mu)\big]\big\ran&.
\end{split}\end{equation*}
\end{lmm}

\noindent
{\it Proof:}
(1) We continue with the notation used in the proof of Lemma~\ref{n2cusps}.
By definition, Proposition~\ref{euler_prp}, and equation~\e_ref{n2cusps_e3},
\begin{equation}\label{n2cusps2_e1}\begin{split}
\big|{\cal S}_{1;2}(\mu)\big|&=
^{\pm}\!\big|{\cal D}^{(2)~\!-1}(0)\cap{\cal S}_1(\mu)\big|
=\big\lan e\big(L^{*\otimes2}\!\otimes\ev^*T\PP\big),
\big[\bar{\cal S}_1(\mu)\big]\big\ran-
{\cal C}_{\partial\bar{\cal S}_1(\mu)}\big({\cal D}^{(2)}\big)\\
&=\big\lan 3a^2 + 6ac_1({\cal L}^*) + 4c_1^2(L^*),
\big[\bar{\cal S}_1(\mu)\big]\big\ran-
{\cal C}_{\partial\bar{\cal S}_1(\mu)}\big({\cal D}^{(2)}\big).
\end{split}\end{equation}
Suppose  ${\cal T}\!=\!(S^2,[N],I;j,\under{d})\!<\!{\cal T}^*$.
If $d_{\hat{0}}\!\neq\!0$, ${\cal D}^{(2)}$ does not vanish 
on~${\cal U}_{\cal T}(\mu)$ by Proposition~\ref{n2cusps_str}.
Thus, from now on we consider only bubble types ${\cal T}$
such \hbox{that $d_{\hat{0}}\!=\!0$}.\\
(2) Suppose $\chi({\cal T})\!=\!\{\hat{1}\}$ is a one-element set.
Then $\bar{\cal S}_1(\mu)\cap{\cal U}_{\cal T}(\mu)\!=\!
{\cal S}_{{\cal T};1}(\mu)$ and with appropriate identifications
$$\Big|{\cal D}_{{\cal T}^*,\hat{0}}^{(2)}\big(\ga_{{\cal T};1}(\ups)\big)
-{\cal D}_{{\cal T},\hat{1}}^{(2)}
\big(\tilde{\ups}_{\hat{1}}\otimes\tilde{\ups}_{\hat{1}}\big)\Big|
\le C(b_{\ups})|\ups|^{\frac{1}{p}}\big|\tilde{\ups}_{\hat{1}}\big|^2
\qquad\forall\ups\!\in\!{\cal FS}_{{\cal T};1;\de}\!={\cal FT}_{\de},$$
where $\ga_{{\cal T};1}\!:{\cal FT}_{\de}\!\lra\!\bar{\cal S}_1(\mu)$
is the map of Proposition~\ref{n2cusps_str}.
Since $d_{\hat{1}}\!\neq\!0$, ${\cal D}_{{\cal T},\hat{1}}^{(2)}$
does not vanish on~${\cal S}_{{\cal T};1}(\mu)$.
Thus, if $\hat{I}\!\neq\!H_{\hat{0}}{\cal T}$,
${\cal T}$ is ${\cal D}^{(2)}$-hollow and
${\cal C}_{\bar{\cal S}_1(\mu)\cap{\cal U}_{\cal T}(\mu)}
\big({\cal D}^{(2)}\big)\!=\!0$.
If $\hat{I}\!=\!H_{\hat{0}}{\cal T}$, 
i.e.~${\cal T}\!=\!{\cal T}^*(l)$ for some $l\!\in\![N]$,
by Proposition~\ref{euler_prp} and the splitting~\e_ref{cart_split},
\begin{equation}\label{n2cusps2_e3a}
{\cal C}_{\bar{\cal S}_1(\mu)\cap{\cal U}_{\cal T}(\mu)}
\big({\cal D}^{(2)}\big)
=2N\big({\cal D}_{\bar{\cal T},\hat{1}}^{(2)}\big),
\end{equation}
where  ${\cal D}_{\bar{\cal T},\hat{1}}^{(2)}\!\in\!
\Ga(\bar{\cal S}_{\bar{\cal T};1}(\mu);
\hbox{Hom}(L_{\bar{\cal T},\hat{1}}^{\otimes2};\ev^*T\PP)\big)$.
By Lemma~\ref{zeros_main},
\begin{equation}\label{n2cusps2_e3b}\begin{split}
N\big({\cal D}_{\bar{\cal T},\hat{1}}^{(2)}\big)
&=\big\lan c_1\big(\ev^*T\PP\big)-
c_1\big(L_{\bar{\cal T},\hat{1}}^{\otimes2}\big),
\big[\bar{\cal S}_{\bar{\cal T};1}(\mu)\big]\big\ran
-{\cal C}_{{\cal D}_{\bar{\cal T},\hat{1}}^{(2)~\!-1}(0)}
\big({\cal D}_{\bar{\cal T},\hat{1}}^{(2)\perp}\big)\\
&=2\big\lan c_1\big(L_{{\cal T},\hat{1}}^*\big),
\big[\bar{\cal S}_{{\cal T};1}(\mu)\big]\big\ran
-{\cal C}_{{\cal D}_{\bar{\cal T},\hat{1}}^{(2)~\!-1}(0)}
\big({\cal D}_{\bar{\cal T},\hat{1}}^{(2)\perp}\big),
\end{split}\end{equation}
since $\ev\big|{\cal U}_{{\cal T}^*(l)}(\mu)\!=\!\mu_l$.
If ${\cal T}'\!=\!(S^2,[N]\!-\!\{l\},I';j',\under{d}')\!\le\!\bar{\cal T}$
and $d_{\hat{1}}'\!\neq\!0$, by Corollary~\ref{enum1-reg_crl2} in~\cite{Z2},
${\cal D}_{\bar{\cal T},\hat{1}}^{(2)}$ does not vanish on
${\cal U}_{{\cal T}'}(\mu)$ if the constraints~$\mu$ are in general position.
If $d_{\hat{1}}'\!=\!0$ and 
\hbox{$\bar{\cal S}_{\bar{\cal T};1}(\mu)
\cap{\cal U}_{{\cal T}'}(\mu)\!\neq\!\eset$},
Proposition~\ref{n2cusps_str} implies that
$H_{\hat{1}}{\cal T}'\!=\!\hat{I}'$ is a two-element set. 
Furthermore, in such a~case,
$$\bar{\cal S}_{\bar{\cal T};1}(\mu)\cap{\cal U}_{{\cal T}'}(\mu)
={\cal S}_{{\cal T}';2}(\mu),$$
is a finite set and there exists an identification~$\ga_{{\cal T}';2}$
of neighborhoods of ${\cal S}_{{\cal T}';2}(\mu)$ in
$\ga_{E{\cal T}'}$ and in~$\bar{\cal S}_{\bar{\cal T};1}(\mu)$ such that
\begin{equation}\label{n2cusps2_e3c}
\Big|{\cal D}_{\bar{\cal T},\hat{1}}^{(2)}\ga_{{\cal T}';2}(\ups)
-\al_{{\cal T}';k}(\ups)\Big|\le C|\ups|^{1+\frac{1}{p}}
\qquad\forall\ups\!\in\!\ga_{E{\cal T}',\de},
\end{equation}
where $\al_{{\cal T}';k}\!\in\!\Ga({\cal S}_{{\cal T}';2}(\mu);
\hbox{Hom}(\ga_{E{\cal T}'};
L_{{\cal T}',\hat{1}}^{*\otimes2}\!\otimes\ev^*T\PP)\big)$
is an injection on every fiber.
On the other hand, 
${\cal D}_{\bar{\cal T},\hat{1}}^{(2)\perp}\!=\!
\pi_{\bar{\nu}}^{\perp}\!\circ\!{\cal D}_{\bar{\cal T},\hat{1}}^{(2)}$,
where 
$$\pi_{\bar{\nu}}^{\perp}\!:L_{{\cal T}',\hat{1}}^{*\otimes2}
\!\otimes\ev^*T\PP\!
\lra L_{{\cal T}',\hat{1}}^{*\otimes2}\!\otimes\ev^*T\PP/\Bbb{C}\bar{\nu}$$
is the projection onto the quotient by a trivial 
subbundle~$\Bbb{C}\bar{\nu}$; see Subsection~\ref{top_sec}.
Since $\al_{{\cal T}';2}$ is an injection,
$\pi_{\bar{\nu}}^{\perp}\!\circ\!\al_{{\cal T}';2}$ is an isomorphism between
$\ga_{E{\cal T}'}$ and
$L_{{\cal T}',\hat{1}}^{*\otimes2}\!\otimes\ev^*T\PP/\Bbb{C}\bar{\nu}$
over every point of~${\cal S}_{{\cal T}';2}(\mu)$
if $\bar{\nu}$ is generic.
Thus, by Proposition~\ref{euler_prp} and the estimate~\e_ref{n2cusps2_e3c},
\begin{equation}\label{n2cusps2_e3d}
{\cal C}_{{\cal S}_{{\cal T}';2}(\mu)}
\big({\cal D}_{\bar{\cal T},\hat{1}}^{(2)\perp}\big)
=\big|{\cal S}_{{\cal T}';2}(\mu)\big|.
\end{equation}
Combining equations~\e_ref{n2cusps2_e3a}-\e_ref{n2cusps2_e3d},
using~\e_ref{normal_bundle2}, 
and summing over the bubble types ${\cal T}$ with 
\hbox{$\chi({\cal T})\!=\!\{\hat{1}\}$}, we obtain
\begin{equation}\label{n2cusps2_e3}
\sum_{|\chi({\cal T})|=1}
{\cal C}_{\bar{\cal S}_1(\mu)\cap{\cal U}_{\cal T}(\mu)}
\big({\cal D}^{(2)}\big)=
\sum_{l\in[N]}\big\lan 4c_1({\cal L}^*),
\big[\bar{\cal S}_{{\cal T}^*(l);1}(\mu)\big]\big\ran
-2\big|{\cal S}_{2,1;2}(\mu)\big|,
\end{equation}
where ${\cal S}_{2,1;2}(\mu)$ denotes the set of two-component
rational degree-$d$ curves that pass through the $3d\!-\!4$ points
and have a tacnode at one of the points, which is a node common
to both irreducible components.\\
(3) Suppose $\chi({\cal T})\!=\!\{\hat{1},\hat{2}\}$ is a two-element set.  
If $\chi({\cal T})\!=\!\hat{I}$ and $M_{\hat{0}}{\cal T}\!=\!\eset$,
$$\bar{\cal S}_1(\mu)\cap{\cal U}_{\cal T}(\mu)
={\cal S}_{{\cal T};1}(\mu)\cup{\cal S}_{{\cal T};2}(\mu);$$
see Proposition~\ref{n2cusps_str}.
By Corollary~\ref{enum1-reg_crl2} in~\cite{Z2},
the images of ${\cal D}_{{\cal T},\hat{1}}^{(2)}$
and ${\cal D}_{{\cal T},\hat{2}}^{(1)}$ are transversal
in $\ev^*T\PP\!\lra\!{\cal U}_{\cal T}(\mu)$.
Since ${\cal S}_{{\cal T};1}(\mu)$ is a finite set,
it follows that the section $\al_{{\cal T};1}$ of 
Proposition~\ref{n2cusps_str} defines an isomorphism
between ${\cal FS}_{{\cal T};1}$ and $L^*\!\otimes\ev^*T\PP$
over every point of~${\cal S}_{{\cal T};1}(\mu)$.
Thus, by Proposition~\ref{euler_prp} and
the analytic estimate of Proposition~\ref{n2cusps_str},
$${\cal C}_{{\cal S}_{{\cal T};1}(\mu)}\big({\cal D}^{(2)}\big)
=2\big|{\cal S}_{{\cal T};1}(\mu)\big|\Lra
{\cal C}_{{\cal S}_{2;1}(\mu)}\big({\cal D}^{(2)}\big)
=2\big|{\cal S}_{2;1}(\mu)\big|.$$
Similarly, since $\al_{{\cal T};2}$ does not vanish on 
${\cal S}_{{\cal T};2}(\mu)$ and extends naturally 
over~$\bar{\cal S}_{{\cal T};2}(\mu)$,
$${\cal C}_{{\cal S}_{{\cal T};2}(\mu)}\big({\cal D}^{(2)}\big)
=N\big(\al_{{\cal T};2}\big)
=N\big({\cal D}_{2;2}).$$
If $\chi({\cal T})\!\neq\!\hat{I}$ or $M_{\hat{0}}{\cal T}\!\neq\!\eset$,
by Proposition~\ref{n2cusps_str}, 
$\bar{\cal S}_1(\mu)\cap{\cal U}_{\cal T}(\mu)
\!=\!{\cal S}_{{\cal T};2}(\mu)$.
The section $\al_{{\cal T};2}$ has full rank on every fiber
in these cases.
Thus, by Proposition~\ref{euler_prp} and 
the analytic estimate of  Proposition~\ref{n2cusps_str},
if $\hat{I}\!\neq\!H_{\hat{0}}{\cal T}$,
${\cal T}$ is ${\cal D}^{(2)}$-hollow and
${\cal C}_{\bar{\cal S}_1(\mu)\cap{\cal U}_{\cal T}(\mu)}
\big({\cal D}^{(2)}\big)\!=\!0$.
If $\chi({\cal T})\!=\!\hat{I}$ and 
$M_{\hat{0}}{\cal T}\!\neq\!\eset$, $\al_{{\cal T};2}$ extends over 
$\bar{\cal S}_{{\cal T};2}(\mu)$ via the splitting~\e_ref{cart_split} by
$$\al_{{\cal T};2}
\big[x_{\hat{1}},x_{\hat{2}},y_l,b;\ups_{\hat{1}},\ups_{\hat{2}}\big]
=x_{\hat{1}}{\cal D}_{{\cal T},\hat{1}}\ups_{\hat{1}}
+x_{\hat{2}}{\cal D}_{{\cal T},\hat{2}}\ups_{\hat{2}}.$$
This extension vanishes only on the set $x_{\hat{1}}\!=\!x_{\hat{2}}$.
Thus, by Proposition~\ref{euler_prp} and Lemma~\ref{zeros_main},
$${\cal C}_{\bar{\cal S}_1(\mu)\cap{\cal U}_{\cal T}(\mu)}
\big({\cal D}^{(2)}\big)=
\big\lan c_1\big(L^{*\otimes2}\!\otimes\ev^*T\P\big)
-c_1\big(\ga_{L^*\otimes L_{\hat{1}}\oplus L^*\otimes L_{\hat{2}}}\big),
\big[\bar{\cal S}_{{\cal T};2}(\mu)]\big\ran
-{\cal C}_{\al_{{\cal T};2}^{-1}(0)}\big(\al_{{\cal T};2}^{\perp}\big)
=2\big|{\cal S}_{\bar{\cal T};2}(\mu)\big|.$$
Here we used $\lan c_1(L^*),[\bar{\cal M}_{0,4}]\ran\!=\!1$;
see Corollary~\ref{enum1-chern_mbar4} in~\cite{Z2}.
Summing over all bubble types ${\cal T}$ as above and
using Lemma~\ref{n2cusps2l}, we obtain
\begin{equation}\label{n2cusps2_e5}\begin{split}
\sum_{|\chi({\cal T})|=2}
{\cal C}_{\bar{\cal S}_1(\mu)\cap{\cal U}_{\cal T}(\mu)}
\big({\cal D}^{(2)}\big)=
\big\lan 6a^2\!+\!3a\big(c_1({\cal L}_1^*)\!+\!c_1({\cal L}_2^*)\big)
\!+\!c_1({\cal L}_1^*)c_1({\cal L}_2^*),
\big[\bar{\cal V}_2(\mu)\big]\big\ran&\\
+2\big|S_{2;1}(\mu)\big|+2\big|{\cal S}_{2,1;2}(\mu)\big|&.
\end{split}\end{equation}
(4) Finally, suppose $\chi({\cal T})\!=\!\{\hat{1},\hat{2},\hat{3}\}$
is a three-element set.
By Proposition~\ref{n2cusps_str},
$$\bar{\cal S}_1(\mu)\cap{\cal U}_{\cal T}(\mu)=
{\cal S}_{{\cal T};2}(\mu)={\cal U}_{\cal T}(\mu).$$
The section $\al_{{\cal T};2}$ again has full rank.
Thus, ${\cal S}_{{\cal T};2}(\mu)$ is ${\cal D}^{(2)}$-hollow
if $\chi({\cal T})\!\neq\!\hat{I}$.
If $\chi({\cal T})\!=\!\hat{I}$,  $\al_{{\cal T};2}$ extends over 
$\bar{\cal S}_{{\cal T};2}(\mu)$ via the splitting~\e_ref{cart_split} by
$$\al_{{\cal T};2}
\big[x_{\hat{1}},x_{\hat{2}},x_{\hat{3}},b;
\ups_{\hat{1}},\ups_{\hat{2}},\ups_{\hat{3}}\big]
=x_{\hat{1}}{\cal D}_{{\cal T},\hat{1}}\ups_{\hat{1}}
+x_{\hat{2}}{\cal D}_{{\cal T},\hat{2}}\ups_{\hat{2}}
+x_{\hat{3}}{\cal D}_{{\cal T},\hat{3}}\ups_{\hat{3}}.$$
This extension does not vanish on ${\cal FS}_{{\cal T};2}(\mu)$, since 
$x_{\hat{1}},x_{\hat{2}},x_{\hat{3}}$ are never all the same.
Thus, by Proposition~\ref{euler_prp},
$${\cal C}_{\bar{\cal S}_1(\mu)\cap{\cal U}_{\cal T}(\mu)}
\big({\cal D}^{(2)}\big)=
\big\lan c_1\big(L^{*\otimes2}\!\otimes\ev^*T\P\big)
-c_1\big(
\ga_{L^*\otimes L_{\hat{1}}\oplus L^*\otimes L_{\hat{2}}
\oplus L^*\otimes L_{\hat{3}}}\big),
\big[\bar{\cal S}_{{\cal T};2}(\mu)]\big\ran
=3\big|{\cal U}_{\bar{\cal T}}(\mu)\big|.$$
Summing over all such bubble types ${\cal T}$, we obtain
\begin{equation}\label{n2cusps2_e7}
\sum_{|\chi({\cal T})|=3}
{\cal C}_{\bar{\cal S}_1(\mu)\cap{\cal U}_{\cal T}(\mu)}
\big({\cal D}^{(2)}\big)=3\tau_3(\mu).
\end{equation}
The claim follows by
plugging the sum of equations~\e_ref{n2cusps2_e3},
\e_ref{n2cusps2_e5}, and~\e_ref{n2cusps2_e7}
into~\e_ref{n2cusps2_e1} and using~\e_ref{normal_bundle1}.

\begin{crl}
\label{n2cusps2c}
If $d$ is a positive integer,
the number of rational degree-$d$ curves that pass 
through a tuple~$\mu$ of $3d\!-\!4$ points in general position in~$\PP$
and have a $(3,4)$-cusp is given~by
\begin{equation*}\begin{split}
\big|{\cal S}_{1;2}(\mu)\big|=&
\big\lan 33a^2c_1^2({\cal L}^*)+18ac_1^3({\cal L}^*)+
4c_1^4({\cal L}^*),\big[\bar{\cal V}_1(\mu)\big]\big\ran+3\tau_3(\mu)\\
&-\big\lan 21a^2+9a\big(c_1({\cal L}_1^*)\!+\!c_1({\cal L}_2^*)\big)
+2\big(c_1^2({\cal L}_1^*)\!+\!c_1^2({\cal L}_2^*)\big)
+c_1({\cal L}_1^*)c_1({\cal L}_2^*),
\big[\bar{\cal V}_2(\mu)\big]\big\ran
\end{split}\end{equation*}
\end{crl}

\begin{center}
\begin{tabular}{||c|c|c|c|c|c|c||}
\hline\hline
$d$&  2&3&4&5&6&7\\
\hline
$|{\cal S}_{1;2}(\mu)|$& 
0&0&147& 54,612& 23,177,124& 14,617,373,280\\ 
\hline\hline
\end{tabular}
\end{center}

\noindent
This corollary is immediate from Lemmas~\ref{n2cusps2},
\ref{n2_2cusps},  and~\ref{n2cusps}.

\begin{lmm}
\label{n2cusps2l}
If $d$ is a positive integer and $\mu$ is a 
a tuple of $3d\!-\!4$ points in general position in~$\PP$,
$$N\big({\cal D}_{2;2}\big)=
\big\lan 6a^2+3a\big(c_1({\cal L}_1^*)\!+\!c_1({\cal L}_2^*)\big)
+c_1({\cal L}_1^*)c_1({\cal L}_2^*),
\big[\bar{\cal V}_2(\mu)\big]\big\ran.$$
\end{lmm}

\noindent
{\it Proof:} 
(1) Since ${\cal D}_{2;2}$ does not vanish on ${\cal S}_{2;2}(\mu)$,
by Lemmas~\ref{zeros_main} and~\ref{n2tacnodes},
\begin{gather}
N\big({\cal D}_{2;2}\big)
\!=\!\big\lan c_1(\ev^*T\PP)\!-\!c_1\big(\ga_{E_2}\big),
\big[\bar{\cal S}_{2;2}(\mu)\big]\big\ran-
{\cal C}_{\partial\bar{\cal S}_{2;2}(\mu)}\big({\cal D}_{2;2}^{\perp}\big)
\!=\!\big\lan 3a\!+\!\la_{E_2},
\big[\bar{\cal S}_{2;2}(\mu)\big]\big\ran\!-\!
{\cal C}_{\partial\bar{\cal S}_{2;2}(\mu)}\big({\cal D}_{2;2}^{\perp}\big)
\notag\\
\label{n2cusps2l_e1}
=\big\lan 12a^2\!+\!6a\big(c_1({\cal L}_1^*)\!+\!c_1({\cal L}_2^*)\big)
\!+\!\big(c_1^2(L_1^*)\!+\!c_1^2(L_2^*)\big)\!+\!c_1(L_1^*)c_1(L_2^*),
\big[\bar{\cal V}_2(\mu)\big]\big\ran-
{\cal C}_{\partial\bar{\cal S}_{2;2}(\mu)}\big({\cal D}_{2;2}^{\perp}\big)
\end{gather}
In \e_ref{n2cusps2l_e1}, 
${\cal D}_{2;2}^{\perp}\!=\!\pi_{\bar{\nu}}^{\perp}\!\circ\!{\cal D}_{2;2}$,
where $\pi_{\bar{\nu}}^{\perp}
        \!\!:\ev^*T\PP\!\!\lra\!\ev^*T\PP/\Bbb{C}\bar{\nu}$
is the projection onto the quotient by a trivial 
subbundle~$\Bbb{C}\bar{\nu}$; see Subsection~\ref{top_sec}.
If $\bar{\nu}$ is generic, by Proposition~\ref{n2tacnodes_str},
$\pi_{\bar{\nu}}^{\perp}\!\circ\!\al_{2;2}$ is an isomorphism 
between $\ga_{E_2}^*\!\otimes\big(E_2/\ga_{E_2}\big)$  and
$\ev^*T\PP/\Bbb{C}\bar{\nu}$ 
over every point of~$\partial\bar{\cal S}_{2;2}(\mu)$.
Thus, by Proposition~\ref{euler_prp} and the analytic estimate of
Proposition~\ref{n2tacnodes_str},
\begin{equation}\label{n2cusps2l_e3}
{\cal C}_{\partial\bar{\cal S}_{2;2}(\mu)}\big({\cal D}_{2;2}^{\perp}\big)
=\big|\partial\bar{\cal S}_{2;2}(\mu)\big|
=\big|{\cal S}_{2;1}(\mu)\big|+
\sum_{[{\cal T}]}\big|{\cal S}_{{\cal T};2}(\mu)\big|,
\end{equation}
where the sum is taken over all equivalence classes of non-basic types
${\cal T}\!=\!\big(S^2,[N],I;j,\under{d})$ such that
$I\!-\!\hat{I}\!=\!\{k_1,k_2\}$ is a two-element set 
\hbox{and $\sum d_i\!=\!d$}.\\
(2) Let ${\cal T}_i\!=\!\big(S^2,M_{k_i},I_{k_i};j,\under{d})$,
where $i=1,2$, be the simple bubble types corresponding 
to a bubble type~${\cal T}$ as above.
If ${\cal S}_{{\cal T};2}(\mu)\!\neq\!\eset$,
up to a re-ordering of indices, ${\cal T}$ must have one of two forms.
The first possibility is that ${\cal T}_2$ is basic, while 
$d_{k_1}\!=\!0$ and  $H_{k_1}{\cal T}\!=\!I_{k_1}$ is a two-element set.
Then ${\cal S}_{{\cal T};2}(\mu)\!=\!{\cal U}_{\cal T}(\mu)$.
The sum of the cardinalities of the sets ${\cal S}_{{\cal T};2}(\mu)$
taken over all equivalence classes of such bubble types
is then~$3\tau_3(\mu)$, since one of the three irreducible components
of the image of each map is distinguished.
The other possibility is that ${\cal T}_2$ is basic, while 
$d_{k_1}\!=\!0$,
$H_{k_1}{\cal T}\!=\!I_{k_1}$ is a one-element set,
and $j_l\!=\!k_1$ for some $l\!\in\![N]$.
Since $\ev^*{\cal O}(1_{\PP})$ is trivial on~$\bar{\cal U}_{\cal T}(\mu)$,
\begin{equation*}\begin{split}
\big|{\cal S}_{{\cal T};2}(\mu)\big|
=\big\lan c_2\big(\ga_{E\bar{\cal T}}^*\!\otimes\ev^*T\PP\big),
\big[E\bar{\cal T}\big]\big\ran 
&=\big\lan c_1(L_{\bar{\cal T},\hat{1}}^*)+c_1(L_{\bar{\cal T},k_2}^*),
\big[\bar{\cal U}_{\bar{\cal T}}(\mu)\big]\big\ran \\
&=\big\lan  c_1(L_{{\cal T},\hat{1}}^*)+c_1(L_{{\cal T},k_2}^*),
\big[\bar{\cal U}_{{\cal T}}(\mu)\big]\big\ran,
\end{split}\end{equation*}
if $I\!=\!\{k_1,k_2,\hat{1}\big\}$.
Summing over all equivalence classes of such bubble types 
and using equations~\e_ref{normal_bundle1} and~\e_ref{normal_bundle2},
we obtain
$$\sum_{[{\cal T}]}\big|{\cal S}_{{\cal T};2}(\mu)\big|
=\sum_{[{\cal T}^*]}\sum_{l\in[N]}
\big\lan c_1({\cal L}_1^*)+c_1({\cal L}_2^*),
\big[\bar{\cal U}_{{\cal T}^*(l)}(\mu)\big]\big\ran+3\tau_3(\mu),$$
where the second sum is taken over equivalence classes
of basic bubble types ${\cal T}^*\!=\!(S^2,[N],I^*;j^*,\under{d}^*)$
such that $|I^*|\!=\!2$ \hbox{and $\sum d_i^*\!=\!d$}.
Combing equations~\e_ref{n2cusps2l_e1} and~\e_ref{n2cusps2l_e3} thus gives
$$N\big({\cal D}_{2;2}\big)\!=\!
\big\lan 12a^2\!+\!6a\big(c_1({\cal L}_1^*)\!+\!c_1({\cal L}_2^*)\big)
\!+\!\big(c_1^2({\cal L}_1^*)+c_1^2({\cal L}_2^*)\big)
\!+\!c_1({\cal L}_1^*)c_1({\cal L}_2^*),
\big[\bar{\cal V}_2(\mu)\big]\big\ran
-\big|S_{2;1}(\mu)\big|-3\tau_3(\mu).$$
The claim now follows by using Lemma~\ref{n2_2cusps}.

\section{Computation of $CR_3(\mu)$}
\label{comp_sect}

\subsection{The Numbers $n_1^{(3)}(\mu)$, $n_2^{(2)}(\mu)$, and 
$n_3^{(1)}(\mu)$}

\noindent
The goal of this section is to give topological formulas
for the six numbers~$n_m^{(k)}(\mu)$ of Corollary~\ref{CR3_str}.
We start with the three numbers that involve zero-dimensional
spaces of rational maps, 
i.e.~${\cal S}_{1;2}(\mu)$,  ${\cal S}_{2;1}(\mu)$,
and~${\cal V}_3(\mu)$.

\begin{lmm}
\label{m1k3} $n_1^{(3)}(\mu)=12\big|{\cal S}_{1;2}(\mu)\big|$
\end{lmm}

\noindent
{\it Proof:} By Subsection~\ref{str_thm_sec}, 
$n_1^{(3)}(\mu)\!=\!N(\al_{1;2})$, where 
\begin{gather}
\al_{1;2}\!\in\!\Ga\big(\Si\!\times\!{\cal S}_{1;2}(\mu);
\hbox{Hom}(T\Si^{\otimes3}\!\otimes\! L^{\otimes3}\!,
{\cal H}_{\Si}^{- -}\!\otimes\!\ev^*T\PP)\big),\notag\\
\label{m1k3_e1}
\al_{1;2}\big(x,b,v\otimes\ups)=
\big({\cal D}^{(3)}\ups\big)\big(s_x^{(3)}v\big).
\end{gather}
The section 
$s^{(3)}\!\!\in\!\Ga\big(\Si;
 \hbox{Hom}(T\Si^{\otimes3},{\cal H}_{\Si}^{- -})\big)$ 
has simple zeros at $z_1,\ldots,z_{24}$.
Thus, $s^{(3)}$ induces a non-vanishing section
$$\tilde{s}^{(3)}\!\in\!\Ga\big(\Si;
\hbox{Hom}(\tilde{T}\Si,{\cal H}_{\Si}^{- -})\big),
\quad\hbox{where}\quad
\tilde{T}\Si=T\Si^{\otimes3}\otimes{\cal O}(z_1)\otimes\ldots\otimes
{\cal O}(z_{24}).$$
Furthermore, $N(\tilde{\al}_{1;2})\!=\!N(\al_{1;2})$, where
$$\tilde{\al}_{1;2}\!\in\!\Ga\big(\Si\!\times\!{\cal S}_{1;2}(\mu);
\hbox{Hom}(\tilde{T}\Si\!\otimes L^{\otimes3},
{\cal H}_{\Si}^{- -}\!\otimes\ev^*T\PP)\big)$$
is the section obtained by replacing $s^{(3)}$ 
by $\tilde{s}^{(3)}$ in~\e_ref{m1k3_e1}.
See Subsection~\ref{enum1-comp_sec2} in~\cite{Z2} for a similar argument
in the genus-two case.
Since ${\cal D}^{(3)}$ does not vanish on~${\cal S}_{1;2}(\mu)$
by Corollary~\ref{enum1-reg_crl2} in~\cite{Z2},
$\tilde{\al}_{1;2}$ does not vanish on~$\Si\!\times\!{\cal S}_{1;2}(\mu)$.
Thus, by Lemma~\ref{zeros_main},
$$n_1^{(3)}(\mu)=N\big(\tilde{\al}_{1;2}\big)=
\big\lan c_1\big({\cal H}_{\Si}^{- -}\!\otimes\ev^*T\PP\big)
-c_1(\tilde{T}\Si),\big[\Si\!\times\!{\cal S}_{1;2}(\mu)\big]\big\ran
=12\big|{\cal S}_{1;2}(\mu)\big|.$$
since the euler characteristic of $\Si$ is~$-4$.\\

\noindent
{\it Remark:} Note that this argument remains valid even
if $\Si$ has hyperflexes, 
i.e.~the points  $z_1,\ldots,z_{24}$ are not all distinct.

\begin{lmm}
\label{m2k2} $n_2^{(2)}(\mu)=36\big|{\cal S}_{2;1}(\mu)\big|$
\end{lmm}

\noindent
{\it Proof:} (1) By Subsection~\ref{str_thm_sec}, 
$n_2^{(2)}(\mu)\!=\!N(\al_{2;1})$, where 
\begin{gather*}
\al_{2;1}\!\in\!\Ga\big(
\Si_1\!\times\!\Si_2\!\times\!{\cal S}_{2;1}(\mu);
\hbox{Hom}(E,{\cal O})\big),
\quad
E=T\Si_1^{\otimes2}\!\otimes L_1^{\otimes2}\!
\oplus T\Si_2\!\otimes L_2,~~
{\cal O}={\cal H}_{\Si_1}^-\!\otimes\ev^*T\PP,\\
\al_{2;1}\big(x_1,x_2,b; v_1\otimes\ups_1,v_2\otimes\ups_2\big)
= \big({\cal D}_1^{(2)}\ups_1\big)\big(s_{x_1}^{(2)}v_1)
+ \big({\cal D}_2\ups_2\big)\big(\pi_{x_1}^-s_{x_2}v_2\big)
\in {\cal H}_{\Si}^-(x_1)\otimes T_{\ev(b)}\PP.
\end{gather*}
Here we define the line bundles $L_1,L_2\!\lra\!{\cal S}_{2;1}(\mu)$
and the sections ${\cal D}_1^{(2)}$ and  ${\cal D}_2^{(1)}$
as follows.
If $b\!\in\!{\cal U}_{{\cal T}^*}(\mu)\cap{\cal S}_{2;1}(\mu)$,
${\cal T}^*\!=\!\big(S^2,[N],I^*;j^*,\under{d}^*)$,
$I^*\!=\!\{k_1,k_2\}$, and ${\cal D}_{{\cal T}^*,k_1}b\!=\!0$,
we take
$$L_1\big|_b=L_{k_1}\big|_b,\quad
L_2\big|_b=L_{k_2}\big|_b,\quad
{\cal D}_1^{(2)}\big|_b={\cal D}_{{\cal T}^*,k_1}^{(2)}\big|_b,\quad
{\cal D}_2^{(1)}\big|_b={\cal D}_{{\cal T}^*,k_2}^{(1)}\big|_b.$$
(2) By Lemma~\ref{zeros_main},
\begin{equation}\label{m2k2_e3}
N(\al_{2;1})=\sum_{k=0}^{k=2}
\big\lan \la_E^{3-k}c_k\big({\cal H}_{\Si_1}^-\!\otimes\!\Bbb{C}^2\big),
\big[\Bbb{P}E\big]\big\ran
-{\cal C}_{\al^{\!E~\!\! -\!1}(0)}\big(\al^{E\perp}\big)
=64\big|{\cal S}_{2;1}(\mu)\big|-
{\cal C}_{\al^{\!E~\!\! -\!1}(0)}\big(\al^{E\perp}\big).
\end{equation}
Since $\Si$ is not hyperelleptic by assumption,
$s_{x_1}\!=\!\la s_{x_2}$ for some $\la\!\in\!C^*$ 
if and only if $x_1\!=\!x_2$.
Thus, $\pi_{x_1}^-s_{x_2}\!=\!0$  if and only if $x_1\!=\!x_2$.
Since the images of ${\cal D}_1^{(2)}\big|_b$ and ${\cal D}_2^{(1)}\big|_b$
in $T_{\ev(b)}\PP$ are linearly independent for all 
$b\!\in\!{\cal S}_{2;1}(\mu)$ by Corollary~\ref{enum1-reg_crl2} in~\cite{Z2},
it follows that
$$\al^{E~\!-1}(0)={\cal Z}\equiv
\big\{ \big(x,x,b;T\Si_2\otimes L_2\big)\!:
x\!\in\!\Si,~b\!\in\!{\cal S}_{2;1}(\mu)\big\}.$$
The normal bundle of ${\cal Z}$ in $\Bbb{P}E_2$ is
$${\cal NZ}=T\Si_2\oplus \big(T\Si_2\otimes L_2)^*\otimes 
 T\Si_1^{\otimes2}\otimes L_1^{\otimes2}
\approx T\Si\oplus T\Si \lra {\cal Z}.$$
With appropriate identifications,
\begin{gather}
\label{m2k2_e5}
\big|\al^E(x,x,b;w,u)-\al_{\cal Z}(x,b;w,u)\big|
\le C|w|^2\qquad\forall(x,x,b;w,u)\!\in\!{\cal NZ},\\
\hbox{where}\qquad
\al_{\cal Z}\in\Ga\big({\cal Z};
\hbox{Hom}({\cal NZ};\ga_E^*\!\otimes{\cal O})\big),\quad
\al_{\cal Z}(x,b;w,u)=
\big\{{\cal D}_1^{(2)}\otimes s_x^{(2)}\big\}\circ u
+{\cal D}_2^{(1)}\otimes s_x^{(2)}(w,\cdot)\notag.
\end{gather}
Since the images of ${\cal D}_1^{(2)}\big|_b$ and ${\cal D}_2^{(1)}\big|_b$
in $T_{\ev(b)}\PP$ are linearly independent for all 
$b\!\in\!{\cal S}_{2;1}(\mu)$,
$\al_{\cal Z}$~has full rank over all of~${\cal Z}$.
If $\bar{\nu}$ is generic, $\pi_{\bar{\nu}}^{\perp}\!\circ\!\al_{\cal Z}$
also has full rank on every fiber, where
$\pi_{\bar{\nu}}^{\perp}\!:{\cal O}\!\lra\!{\cal O}/\Bbb{C}\bar{\nu}$
is the quotient projection as before.
Then by the analytic estimate~\e_ref{m2k2_e5} and 
Proposition~\ref{euler_prp},
\begin{equation}\label{m2k2_e7}\begin{split}
{\cal C}_{\al^{\!E~\!\!-\!1}(0)}\big(\al^{E\perp}\big)&=
\big\lan c_1\big(T^*\Si\otimes{\cal O}^{\perp}\big)-
c_1\big({\cal NZ}\big),\big[\Si\!\times\!{\cal S}_{2;1}(\mu)\big]\big\ran\\
&=\big\lan \big( 3c_1(T^*\Si)\!+\!2c_1(T^*\Si)\big)+2c_1(T^*\Si),
[\Si]\big\ran \big|{\cal S}_{2;1}(\mu)\big|
=28\big|{\cal S}_{2;1}(\mu)\big|.
\end{split}\end{equation}
The claim is obtained by plugging \e_ref{m2k2_e7} into~\e_ref{m2k2_e3}.

\begin{lmm}
\label{m3k1} $n_3^{(1)}(\mu)=36\tau_3(\mu)$
\end{lmm}

\noindent
{\it Proof:} (1) By Subsection~\ref{str_thm_sec}, 
$n_3^{(1)}(\mu)\!=\!N(\al_3)$, where 
\begin{gather*}
\al_3\!\in\!\Ga\big(
\Si_1\!\times\!\Si_2\!\times\!\Si_3\!\times\!{\cal V}_3(\mu);
\hbox{Hom}(E,{\cal O})\big), \quad
E=\bigoplus_i T\Si_i\!\otimes L_i,~~
{\cal O}={\cal H}_{\Si}^{0,1}\!\otimes\ev^*T\PP,\\
\al_3\big(x_1,x_2,x_3,b; 
v_1\otimes\ups_1,v_2\otimes\ups_2,v_3\otimes\ups_3\big)
= \sum_i  \big({\cal D}_i\ups_i\big)\big(s_{x_i}v_i)
\in {\cal H}_{\Si}^{0,1}\otimes T_{\ev(b)}\PP.
\end{gather*}
Here the bundles $L_i\!\lra\!{\cal V}_3(\mu)$ and the sections 
${\cal D}_i\!\in\!\Ga\big({\cal V}_3(\mu);L_i^*\!\otimes\!\ev^*T\PP\big)$
are defined as follows.
If $b\!\in\!{\cal U}_{{\cal T}^*}(\mu)\!\subset\!{\cal V}_3(\mu)$,
${\cal T}^*\!=\!(S^2,[N],I^*;j^*,\under{d}^*)$,
and $I^*\!=\!\{k_1,k_2,k_3\}$, we let 
$L_i\big|_b\!=\!L_{{\cal T},k_i}$ and 
${\cal D}_i\!=\!{\cal D}_{{\cal T},k_i}$.
These bundles and sections are well-defined once we fix
a representative for each equivalence class of such bubble types~${\cal T}^*$
and order the elements of the corresponding set~$I^*$.\\
(2) By Lemma~\ref{zeros_main},
\begin{equation}\label{m3k1_e3}
N(\al_3)=\sum_{k=0}^{k=3}
\big\lan \la_E^{5-k}c_k\big(\Bbb{C}^5\big),
\big[\Bbb{P}E\big]\big\ran
-{\cal C}_{\al^{\!E~\!\! -\!1}(0)}\big(\al^{E\perp}\big)
=64\big|{\cal V}_3(\mu)\big|-
{\cal C}_{\al^{\!E~\!\! -\!1}(0)}\big(\al^{E\perp}\big).
\end{equation}
Since $\Si$ is not hyperelleptic,
$s_{x_1}\!=\!\la s_{x_2}$ for some $\la\!\in\!C^*$ 
if and only if $x_1\!=\!x_2$.
Since the images of ${\cal D}_i\big|_b$ and ${\cal D}_j\big|_b$
in $T_{\ev(b)}\PP$ are linearly independent for all 
$b\!\in\!{\cal V}_3(\mu)$ and $i\!\neq\!j$
by Corollary~\ref{enum1-reg_crl2} in~\cite{Z2},
it follows that
$$\al^{E~\!-1}(0)={\cal Z}\equiv
\big\{ \big(x,x,x,b;[v\otimes\ups_1,v\otimes\ups_2,v\otimes\ups_3]\big)
\!\in\!\Bbb{P}E\!:
x\!\in\!\Si,~b\!\in\!{\cal V}_3(\mu),~
\sum_i{\cal D}_i\ups_i\!=\!0\big\}.$$
The normal bundle of ${\cal Z}$ in $\Bbb{P}E_2$ is
$${\cal NZ}=T\Si_2\oplus T\Si_3\oplus 
\big(T\Si_1\otimes L_1\big)^*\!\otimes
\big(T\Si_2\otimes L_2\oplus T\Si_3\otimes L_3\big)
\approx T\Si\oplus T\Si\oplus\Bbb{C}^2\lra {\cal Z}.$$
With appropriate identifications,
\begin{gather}
\big|\al^E(x,x,x,b;w_2,w_3,u_2,u_3)\!-\!
\al_{\cal Z}(x,b;w_2,w_3,u_2,u_3)\big|\!\le\! C\big|(w_2,w_3)\big|^2
~~\forall\big(w_2,w_3,u_2,u_3\big)\!\in\!{\cal NZ},\notag\\
\label{m3k1_e5}
\hbox{where}\qquad
\al_{\cal Z}\in\Ga\big({\cal Z};
\hbox{Hom}({\cal NZ};(T\Si_1\otimes L_1)^*\!\otimes{\cal O})\big),\\
\al_{\cal Z}\big(x,b;w_2,w_3,u_2,u_3\big)=
\big\{{\cal D}_2\otimes s_x\big\}\!\circ\! u_2\!+\!
\big\{{\cal D}_3\otimes s_x\big\}\!\circ\! u_3\!+\!
{\cal D}_2\otimes s_{g_x,x}^{(2)}(w_2,\cdot)\!+\!
{\cal D}_3\otimes s_{g_x,x}^{(2)}(w_3,\cdot)\notag,
\end{gather}
and $\big\{g_x\!:x\!\in\!\Si\big\}$ is a smooth family of metrics
on $\Si$ such $g_x$ is flat on a neighborhood of~$x$.
Since the images of ${\cal D}_2\big|_b$ and ${\cal D}_3\big|_b$
are linearly independent in~$T_{\ev(b)}\PP$
for all $b\!\in\!{\cal V}_3(\mu)$ and 
the section $s_x^{(2)}\!=\!\pi_x^-\!\circ\!s_{g_x,x}$
does not vanish on~$\Si$,
the linear map $\al_{\cal Z}$ is injective over~${\cal Z}$.
Thus, by the analytic estimate~\e_ref{m3k1_e5} and 
Proposition~\ref{euler_prp},
\begin{equation}\label{m3k1_e7}\begin{split}
{\cal C}_{\al^{\!E~\!\!-\!1}(0)}\big(\al^{E\perp}\big)&=
\big\lan c_1\big(T^*\Si\otimes{\cal O}^{\perp}\big)-
c_1\big({\cal NZ}\big),\big[\Si\!\times\!{\cal V}_3(\mu)\big]\big\ran\\
&=\big\lan 5c_1(T^*\Si)\!+\!2c_1(T^*\Si),
[\Si]\big\ran \big|{\cal V}_3(\mu)\big|
=28\big|{\cal V}_3(\mu)\big|.
\end{split}\end{equation}
The claim is obtained by plugging \e_ref{m3k1_e7} into~\e_ref{m3k1_e3}.

\subsection{The Number $n_2^{(1)}(\mu)$}

\noindent
We now use the topological tools of Subsection~\ref{top_sec}
along with the analytic estimates of Subsection~\ref{str_global_sec}
to give a topological formula for the number
$n_2^{(1)}(\mu)$ of Corollary~\ref{CR3_str}.
The computation involved is long, but fairly straightforward.

\begin{lmm}
\label{m2k1}
If $d$ is a positive integer and $\mu$ is
a tuple of $3d\!-\!4$ points in general position in~$\PP$,
$$n_2^{(1)}(\mu)=12\big\lan 10a^2
+3a\big(c_1({\cal L}_1^*)\!+\!c_1({\cal L}_1^*)\big)
+c_1({\cal L}_1^*)c_1({\cal L}_1^*)
,\big[\bar{\cal V}_2(\mu)\big]\big\ran.$$
\end{lmm}

\noindent
{\it Proof:} (1) By Subsection~\ref{str_thm_sec}, 
$n_2^{(1)}(\mu)\!=\!N(\al_2)$, where 
\begin{gather*}
\al_2\!\in\!\Ga\big(
\Si_1\!\times\!\Si_2\!\times\!\bar{\cal V}_2(\mu);
\hbox{Hom}(\tilde{E},{\cal O})\big), \quad
\tilde{E}=T\Si_1\!\otimes L_1\oplus T\Si_2\oplus L_2,~~
{\cal O}={\cal H}_{\Si}^{0,1}\!\otimes\ev^*T\PP,\\
\al_2\big(x_1,x_2,b; v_1\otimes\ups_1,v_2\otimes\ups_2\big)
= \big({\cal D}_1\ups_1\big)\big(s_{x_1}v_1)+
\big({\cal D}_2\ups_2\big)\big(s_{x_2}v_2)
\in {\cal H}_{\Si}^{0,1}\otimes T_{\ev(b)}\PP,
\end{gather*}
with the bundles $L_i\!\lra\!\bar{\cal V}_2(\mu)$
and the sections 
${\cal D}_i\!\in\!\Ga\big(\bar{\cal V}_2(\mu);L_i^*\!\otimes\!\ev^*T\PP\big)$
defined as in the proof of Lemma~\ref{m3k1}.
By Lemma~\ref{zeros_main},
\begin{alignat}{1}\label{m2k1_e3}
&N(\al_2)=\sum_{k=0}^{k=3}
\big\lan \la_{\tilde{E}}^{5-k}c_k\big({\cal O}\big),
                           \big[\Bbb{P}\tilde{E}\big]\big\ran
-{\cal C}_{\al^{\!\tilde{E}~\!\! -\!1}(0)}\big(\al^{\tilde{E}\perp}\big)\\
&~~=16\big\lan 36a^2\!+\!18a\big(c_1({\cal L}_1^*)\!+\!c_1({\cal L}_2^*)\big)
\!+\!3\big(c_1^2(L_1^*)\!+\!c_1^2(L_2^*)\big)
\!+\!4c_1(L_1^*)c_1(L_2^*),
\big[\bar{\cal V}_2(\mu)\big]\big\ran
\!-\!{\cal C}_{\al^{\!\tilde{E}~\!\! -\!1}(0)}\big(\al^{\tilde{E}\perp}\big)
\notag.
\end{alignat}
Since $s_{x_1}\!=\!\la s_{x_2}$ for some $\la\!\in\!\Bbb{C}$
if and only if $x_1\!=\!x_2$,
$$\al^{\!\tilde{E}~\!\! -\!1}(0)=
\big\{\big(x,x,b;[\ups_1,\ups_2]\big)\!:
{\cal D}_1\ups_1\!+\!{\cal D}_2\ups_2\!=\!0\big\}\cup
\bigcup_{i=1,2}
\big\{\big(x_1,x_2,b;T_{x_i}\Si_i\otimes L_i|_b\big)\!:
 {\cal D}_ib\!=\!0\big\}.$$
We now partition these sets further and apply 
the topological tools of Subsection~\ref{top_sec}.
As usually, we denote by 
$\bar{\nu}\!\in\!\Ga\big(\Bbb{P}\tilde{E};{\cal O}\big)$
a generic non-vanishing section.\\
(2) We start with the spaces 
$(\Si^2\!-\!\lap)\!\times\!{\cal S}_{2;1}(\mu)$ 
and $\lap\!\times\!{\cal S}_{2;1}(\mu)$.
For notational simplicity, we assume that ${\cal D}_1b\!=\!0$
for $b\!\in\!{\cal S}_{2;1}(\mu)$.
The normal bundle of the subspace
$${\cal Z}_{2;1}(\mu)=
\big\{\big(x_1,x_2,b;T_{x_1}\Si_1\otimes L_1|_b\big)\!:
x_1\!\neq\! x_2,~b\!\in\!{\cal S}_{2;1}(\mu)\big\}$$
in $\Bbb{P}\tilde{E}$ is given by
$${\cal NZ}_{2;1}=
\pi_{\tilde{E}}^*\big(L_1^*\otimes\ev^*T\PP\oplus
T\Si_1^*\otimes L_1^*\otimes T\Si_2\otimes L_2\big)
\approx \Bbb{C}^2\otimes T\Si_1^*\otimes T\Si_2.$$
With appropriate identifications, 
${\cal D}_1(b,X)\!=\!X$ for all 
$X\!\in\! (L_1^*\otimes\ev^*T\PP)_b$ sufficiently small. 
Thus,
$$\al^{\tilde{E}}
\big(x_1,x_2,b;X,u\big)= \al_{2;1}\big(x_1,x_2,b;X,u\big)\equiv
X\otimes s_{\Si,x_1}+ \big\{{\cal D}_2\otimes s_{\Si,x_2}\big\}\circ u.$$
Since $\al_{2;1}$ has full rank on 
${\cal Z}_{2;1}(\mu)\!\approx\!(\Si^2\!-\!\lap)\!\times\!{\cal S}_{2;1}(\mu)$ 
and extends over $\Si^2\!\times\!{\cal S}_{2;1}(\mu)$,
$${\cal C}_{(\Si^2-\lap)\times{\cal S}_{2;1}(\mu)}
\big(\al^{\tilde{E}\perp}\big)=
N\big(\pi^{\perp}_{\bar{\nu}}\circ\al_{2;1}\big),$$
as long as $\bar{\nu}$ is generic.
In fact, it can be assumed the image of~$\bar{\nu}$ is disjoint
from~${\cal H}_{\Si}^+(x_1)\otimes\ev^*T\PP$.
Then $\pi^{\perp}_{\bar{\nu}}\big({\cal H}_{\Si}^+(x_1)\otimes\ev^*T\PP\big)$
is a rank-two subbundle of ${\cal O}^{\perp}$, 
and
$$\pi^{\perp}_{\bar{\nu}}\circ\al_{2;1}\!: \Bbb{C}^2\lra 
\ga_{\tilde{E}}^*\otimes
     \pi^{\perp}_{\bar{\nu}}\big({\cal H}_{\Si}^+(x_1)\otimes\ev^*T\PP\big)$$
is an isomorphism.
It follows that 
$N\big(\pi^{\perp}_{\bar{\nu}}\circ\al_{2;1}\big)\!=\!
N\big(\tilde{\al}_{2;1}\big)$, where
\begin{gather*}
\tilde{\al}_{2;1}\in\Ga\big(\Si^2\!\times\!{\cal S}_{2;1}(\mu);
\hbox{Hom}(F_2;{\cal O}_2)\big),\quad
\tilde{\al}_{2;1}(x_1,x_2,b;u)=
\pi^{\perp}_{\pi_{x_1}^-\bar{\nu}}\circ
\big\{{\cal D}_2\otimes \pi_{x_1}^-s_{x_2}\big\}\circ u,\\
F_2=T\Si_1^*\!\otimes\! L_1^*\!\otimes\! T\Si_2\!\otimes\! L_2\!\approx\!
T\Si_1^*\!\otimes\! T\Si_2,\quad
{\cal O}_2=
T\Si_1^*\!\otimes\! L_1^*\!\otimes\!
\big({\cal H}_{\Si_1}^-\!\otimes\!\ev^*T\PP\big)^{\perp}
\!\approx\! T\Si_1^*\!\otimes\!\big({\cal H}_{\Si_1}^-\!
  \otimes\Bbb{C}^2\big)^{\perp}\!,
\end{gather*}
By Lemma~\ref{zeros_main},
$$N\big(\tilde{\al}_{2;1}\big)=
\big\lan c_1^2(F_2^*)\!+\!c_1(F_2^*)c_1({\cal O}_2),
\big[\Si^2\!\times\!{\cal S}_{2;1}(\mu)\big]\big\ran
-{\cal C}_{\tilde{\al}_{2;1}^{-1}(0)}\big(\tilde{\al}_{2;1}^{\perp}\big)
=48\big|{\cal S}_{2;1}(\mu)\big|
-{\cal C}_{\tilde{\al}_{2;1}^{-1}(0)}\big(\tilde{\al}_{2;1}^{\perp}\big).$$
The zero set of $\tilde{\al}_{2;1}$ is $\lap\!\times\!{\cal S}_{2;1}(\mu)$; 
see the proof of Lemma~\ref{m2k2}.
Its normal bundle is $T\Si_2\!\approx\! T\Si$.
If $\bar{\nu}$ and $\bar{\nu}_2$ are generic, as 
in the proof of Lemma~\ref{m2k2}, we obtain
$$\big|\tilde{\al}_{2;1}^{\perp}(x,x,b;w)-
\tilde{\al}_{2;1;\lap}(x,b;w)\big|\le C|w|^2\quad\forall w\!\in\!T\Si_{\de},$$
where $\tilde{\al}_{2;1;\lap}\!:T\Si\!\lra\!F_2^*\!\otimes{\cal O}_2^{\perp}$
is an injection on every fiber.
Thus, by Proposition~\ref{euler_prp},
$${\cal C}_{\tilde{\al}_{2;1}^{-1}(0)}\big(\tilde{\al}_{2;1}^{\perp}\big)
=\big\lan c_1\big(F_2^*\!\otimes{\cal O}_2^{\perp}\big)-c_1\big(T\Si),
\big[\Si\!\times\!{\cal S}_{2;1}(\mu)\big]\big\ran
=24\big|{\cal S}_{2;1}(\mu)\big|.$$
We conclude that
\begin{equation}\label{m2k1_e5}
{\cal C}_{(\Si^2-\lap)\times{\cal S}_{2;1}(\mu)}\big(\al^{\tilde{E}\perp}\big)
=24\big|{\cal S}_{2;1}(\mu)\big|.
\end{equation}
On the other hand, the space
$\bar{\cal Z}_{2;1}\!-\!{\cal Z}_{2;1}\!\approx
\lap\!\times\!{\cal S}_{2;1}(\mu)$ is 
$\al^{\tilde{E}\perp}$-hollow, and thus
$${\cal C}_{\lap\times{\cal S}_{2;1}(\mu)}\big(\al^{\tilde{E}\perp}\big)=0.
$$
Indeed, its normal bundle in $\Bbb{P}\tilde{E}_2$ is given~by
$${\cal NZ}=\pi_{\tilde{E}}^*\big(
T\Si_2\oplus L_1^*\otimes\ev^*T\PP\oplus 
T\Si_1^*\otimes L_1^*\otimes T\Si_2\otimes L_2\big).$$
With appropriate identifications,
\begin{gather*}
\big|\al^{\tilde{E}}(x,x,b;w,X,u)-
\tilde{\al}(x,b;w,X,u)\big|\le C|w|^2|u|
\quad\forall (w,X,u)\!\in\!{\cal NZ}_{\de},
\quad\hbox{where}\\
\tilde{\al}(x,b;w,X,u)=
X\otimes s_x+\big\{{\cal D}_2\otimes s_x\big\}\circ u+
\big\{ {\cal D}_2\otimes s^{(2)}_{g_x,x}(w,\cdot)\big\}\circ u.
\end{gather*}
Since $\pi_x^-s^{(2)}_{g_x,x}$ does not vanish,
$\tilde{\al}$ is a dominant term for~$\al^{\tilde{E}}$;
the same holds for composites with projection maps.
Since 
$$\hbox{rk}\big({\cal H}_{\Si}^-\otimes\ev^*T\PP\big)^\perp
> \hbox{rk}~T\Si_2\otimes
\big(T\Si_1^*\otimes L_1^*\otimes T\Si_2\otimes L_2\big)
+\frac{1}{2}\dim\big(\lap\!\times\!{\cal S}_{2;1}(\mu)\big),$$
$\lap\!\times\!{\cal S}_{2;1}(\mu)$ is $\al^{\tilde{E}\perp}$-hollow.\\
(3) Suppose ${\cal T}\!=\!(S^2,[N],I;j,\under{d})$ is a non-basic
bubble type and ${\cal D}_1b\!=\!0$ for some 
\hbox{$b\!\in\!{\cal U}_{\cal T}(\mu)\!\subset\!\bar{\cal V}_2(\mu)$}.
Let $I_1$ and $I_2$ be the corresponding rooted trees and 
$k_1\!\in\!I_1$ and $k_2\!\in\!I_2$ the minimal elements.
Then $d_{k_1}\!=\!0$, $d_{k_2}\!\neq\!0$, and
$\big|H_{k_1}{\cal T}\big|\!\in\!\{1,2\}$.
Let 
$${\cal Z}_{\cal T}=\big\{
\big(x_1,x_2,b;T_{x_1}\Si_1\!\otimes L_1|_b\big)\!\in\!\Bbb{P}\tilde{E}\!:
b\!\in\!{\cal U}_{\cal T}(\mu)\big\}.$$
By Theorem~\ref{str_global}, the normal bundle of 
${\cal Z}_{\cal T}$ in $\Bbb{P}\tilde{E}$~is 
$${\cal NZ}_{\cal T}=\pi_{\tilde{E}}^*
\big({\cal FT} \oplus
T\Si_1^*\!\otimes L_1^*\!\otimes T\Si_2\!\otimes L_2\big)
\approx  {\cal FT}\oplus T\Si_1^*\!\otimes T\Si_2\!\otimes L_2.$$
First suppose  $H_{k_1}{\cal T}\!=\!\{\hat{1}\}$ is a one-element set.
Then, with appropriate identifications,
\begin{gather}\label{m2k1_e7b}
\big|\al^{\tilde{E}}\big(x_1,x_2,b;\ups,u\big)
-\al_{{\cal Z}_{\cal T}}\big(x_1,x_2,b;\ups_{\hat{1}},u\big)\big|
\le C(b)|\ups|^{\frac{1}{p}}\big(|\ups_{\hat{1}}|+|u|\big)
\quad\forall (\ups,u)\in{\cal N}_b{\cal Z}_{{\cal T},\de(b)},\\
\hbox{where}\qquad
\al_{{\cal Z}_{\cal T}}\big(x_1,x_2,b;\ups_{\hat{1}},u\big)=
\big({\cal D}_{{\cal T},\hat{1}}\ups_{\hat{1}}\big) \otimes s_{x_1}
+\big\{{\cal D}_2\otimes s_{x_2}\big\}\circ u.\notag
\end{gather}
If $H_{k_1}{\cal T}\!\neq\!\hat{I}$,
the images of ${\cal D}_{{\cal T},\hat{1}}\big|_b$ and 
of ${\cal D}_2\big|_b$ in $T_{\ev(b)}\PP$ are linearly independent
for all~$b\!\in\!{\cal U}_{\cal T}(\mu)$.
Thus, $\al_{{\cal Z}_{\cal T}}$ is injective on every fiber
and ${\cal Z}_{\cal T}$ is $\al^{\tilde{E}\perp}$-hollow
by~\e_ref{m2k1_e7b}, provided $\bar{\nu}$ is generic.
Then, by Proposition~\ref{euler_prp},
$${\cal C}_{{\cal Z}_{\cal T}}\big(\al^{\tilde{E}\perp}\big)=0
\qquad\hbox{if}\quad
\big|{\cal H}_{k_1}{\cal T}\big|=1<|\hat{I}|.$$
If ${\cal H}_{k_1}{\cal T}\!=\!\hat{I}$,
$\al_{{\cal Z}_{\cal T}}$ has full rank outside of the set
$$\tilde{\cal S}_{{\cal T};2}(\mu)
=\big\{\big(x,x,b;T_x\Si_1\!\otimes L_1|_b\big)\!\in\!{\cal Z}_{\cal T}\!:
{\cal D}_{{\cal T},\hat{1}}\ups_{\hat{1}}\!+\!{\cal D}_2\ups_2\!=\!0
\hbox{~for some~}(\ups_{\hat{1}},\ups_2)\!\neq\!0\big\}
\approx \Si\times{\cal S}_{{\cal T};2}(\mu).$$
Since $\al_{{\cal Z}_{\cal T}}$ extends naturally over 
$\bar{\cal Z}_{\cal T}\!\subset\!\Bbb{P}\tilde{E}$,
by  Proposition~\ref{euler_prp},
\begin{gather*}
{\cal C}_{{\cal Z}_{\cal T}-\tilde{\cal S}_{{\cal T};2}(\mu)}
\big(\al^{\tilde{E}\perp}\big)=
N\big(\tilde{\al}_{\cal T}\big),
\qquad\hbox{where}\quad
\tilde{\al}_{\cal T}=
\pi^{\perp}_{\bar{\nu}}\circ\al_{{\cal Z}_{\cal T}}\!\in\!
\Ga\big(\bar{\cal Z}_{\cal T};\hbox{Hom}(F_2;{\cal O}_2)\big),\\
\tilde{\al}_{\cal T}=\pi^{\perp}_{\bar{\nu}}\circ
\big( \big({\cal D}_{{\cal T},\hat{1}}\ups_{\hat{1}}\big) \otimes s_{x_1}
+\big\{{\cal D}_2\otimes s_{x_2}\big\}\circ u\big);\\
F_2=L_1^*\!\otimes\! L_{{\cal T},\hat{1}}\oplus
T\Si_1^*\!\otimes\! L_1^*\!\otimes\! T\Si_2\!\otimes\! L_2\approx
 L_{{\cal T},\hat{1}}\oplus T\Si_1^*\!\otimes T\Si_2\otimes\! L_2,\quad
{\cal O}_2=\ga_{\tilde{E}}^*\otimes{\cal O}^{\perp}
\approx T\Si_1^*\!\otimes\Bbb{C}^5\!.
\end{gather*}
By Lemma~\ref{zeros_main},
\begin{equation*}\begin{split}
N\big(\tilde{\al}_{\cal T}\big)=&
\sum_{k=0}^{k=4}\big\lan \la_{F_2}^{4-k}c_k({\cal O}_2),
     \big[\Bbb{P}F_2\big]\big\ran
-{\cal C}_{\tilde{\al}_{\cal T}^{F_2~\!\!-\!1}(0)}
 \big(\tilde{\al}_{\cal T}^{F_2\perp}\big)\\
=&16\big\lan 3c_1({\cal L}_1^*)\!+\!4c_1({\cal L}_2^*),
\big[\bar{\cal U}_{\cal T}(\mu)\big]\big\ran
-{\cal C}_{\tilde{\al}_{\cal T}^{F_2~\!\!-\!1}(0)}
 \big(\tilde{\al}_{\cal T}^{F_2\perp}\big),
\end{split}\end{equation*}
since $a\!=\!0$, $c_1(L_{{\cal T},\hat{1}}^*)\!=\!c_1({\cal L}_1^*)$,
and $c_1(L_2^*)\!=\!c_1({\cal L}_2^*)$
in $H^*\big(\bar{\cal U}_{\cal T}(\mu)\big)$.
Furthermore,
\begin{gather*}
\tilde{\al}_{\cal T}^{F_2~\!\!-\!1}(0)= \big\{
\big(x,x,b;[v\otimes\ups_{\hat{1}},v\otimes\ups_2]\big)\!\in\!\Bbb{P}F_2\!:
{\cal D}_{{\cal T},\hat{1}}\ups_{\hat{1}}\!+\!{\cal D}_2\ups_2\!=\!0\big\}
\approx \Si\times{\cal S}_{{\cal T};2}(\mu),\\
F_3\equiv{\cal N}\tilde{\al}_{\cal T}^{F_2~\!\!-\!1}(0)=
T\Si_2\oplus\ga_{E{\cal T}}\approx T\Si\oplus\Bbb{C}^2,\\
\big| \tilde{\al}_{\cal T}^{F_2\perp}(x,x,b;w,X)-
\tilde{\al}_{{\cal T},\lap}(x,b;w,X)\big|\le C|w|^2
~~~\forall (w,X)\!\in\!F_{3,\de},
\qquad\hbox{where}\\
\tilde{\al}_{{\cal T},\lap}\!\in\!
\Ga\big(\tilde{\al}_{\cal T}^{F_2~\!\!-\!1}(0);
\hbox{Hom}(F_3,{\cal O}_3)\big),
\quad
{\cal O}_3=\ga_{F_2}^*\otimes{\cal O}_2^{\perp}
\approx \big(T\Si^*\otimes\Bbb{C}^5\big)^{\perp},\\
\tilde{\al}_{{\cal T},\lap}(w,X)=
\pi^{\perp}_{\bar{\nu}_2}\circ \big(\pi^{\perp}_{\bar{\nu}}\circ
\big( X\otimes s_x + {\cal D}_2\otimes s^{(2)}_{g_x,x}(w,\cdot)\big)\big).
\end{gather*}
Since $\tilde{\al}_{{\cal T},\lap}$ has full rank on every fiber,
by Proposition~\ref{euler_prp},
$${\cal C}_{\tilde{\al}_{\cal T}^{F_2~\!\!-\!1}(0)}
\big(\tilde{\al}_{\cal T}^{F_2\perp}\big)=
\big\lan c_1({\cal O}_3)\!-\!c_1(F_3),
\big[\tilde{\al}_{\cal T}^{F_2~\!\!-\!1}(0)\big]\big\ran
=24\big|{\cal S}_{{\cal T};2}(\mu)\big|
=24\big\lan c_1({\cal L}_{{\cal T},\hat{1}}^*)\!+\!c_1({\cal L}_2^*),
\big[\bar{\cal U}_{\cal T}(\mu)\big]\big\ran.$$
On the other hand, $\tilde{\cal S}_{{\cal T};2}(\mu)$ is 
$\al^{\tilde{E}\perp}$-hollow, and thus
$${\cal C}_{\tilde{\cal S}_{{\cal T};2}(\mu)}
\big(\al^{\tilde{E}\perp}\big)=0.$$
Indeed, by Theorem~\ref{str_global},
\begin{gather*}
{\cal N}\tilde{\cal S}_{{\cal T};2}
=T\Si_2\oplus L_2^*\otimes
\big(\hbox{Im}~\!{\cal D}_{{\cal T},\hat{1}}\big)^{\perp}
\oplus L_1^*\otimes L_{{\cal T},\hat{1}}
\oplus T\Si_1^*\otimes L_1^*\otimes T\Si_2\otimes L_2,\\
\big|\al^{\tilde{E}}\big(x,x,b;w,X,\ups_{\hat{1}},u\big)- 
\tilde{\al}\big(x,b;w,X,\ups_{\hat{1}},u)\big|\le 
C\big(|w|+|X|)|w||u|,
\quad\forall \big(w,X,\ups_{\hat{1}},u\big)\!\in
                       \!{\cal N}\tilde{\cal S}_{{\cal T};2;\de},\\
\hbox{where}\qquad
\pi_x^{-1}\tilde{\al}\big(x,b;w,X,\ups_{\hat{1}},u\big)= 
\big\{ {\cal D}_2\otimes s^{(2)}_x(w,\cdot)\big\}\circ u.
\end{gather*}
The claim that $\tilde{\cal S}_{{\cal T};2}(\mu)$ is 
$\al^{\tilde{E}\perp}$-hollow then follows as in~(2).
Summing over bubble types as above, we conclude that
\begin{equation}\label{m2k1_e7}\begin{split}
\sum_{|\chi({\cal T})|=1}
{\cal C}_{\Si^2\times{\cal U}_{\cal T}(\mu)}
\big(\al^{\tilde{E}\perp}\big)=&
16\sum_{{\cal T}^*}\sum_{\{i,j\}=\{1,2\}}\sum_{l\in M_i^*}
\big\lan 3c_1({\cal L}_i^*)\!+\!4c_1({\cal L}_j^*),
\big[\bar{\cal U}_{{\cal T}^*(l)}\big]\big\ran\\
&\quad -24\sum_{{\cal T}^*}\sum_{l\in[N]}
\big\lan c_1({\cal L}_1^*)\!+\!c_1({\cal L}_2^*),
\big[\bar{\cal U}_{{\cal T}^*(l)}\big]\big\ran,
\end{split}\end{equation}
where the outer sums are taken over all equivalence classes
of basic bubble types ${\cal T}^*$ such that 
${\cal U}_{{\cal T}^*}(\mu)\!\subset\!{\cal V}_2(\mu)$.\\
(4) We next consider the case $H_{k_1}{\cal T}\!=\!\{\hat{1},\hat{2}\}$
is a two-element set. Then $\hat{I}\!=\!H_{k_1}{\cal T}$,
\begin{gather*}
\big| \al^{\tilde{E}}\big(x_1,x_2,b;\ups,u\big)
-\al_{{\cal Z}_{\cal T}}\big(x_1,x_2,b;\ups,u\big)\big|
\le C(b)|\ups|^{\frac{1}{p}}\big(|\ups|+|u|\big)
\quad\forall (\ups,u)\in{\cal N}_b{\cal Z}_{{\cal T},\de(b)},\\
\hbox{where}\qquad
\al_{{\cal Z}_{\cal T}}\big(x_1,x_2,b;\ups,u\big)=
\big({\cal D}_{{\cal T},\hat{1}}\ups_{\hat{1}}\big) \otimes s_{x_1}
+\big({\cal D}_{{\cal T},\hat{2}}\ups_{\hat{2}}\big) \otimes s_{x_1}
+\big\{{\cal D}_2\otimes s_{x_2}\big\}\circ u.\notag
\end{gather*}
The map $\al_{{\cal Z}_{\cal T}}$ has full rank outside of the set
$$\tilde{\cal S}_{{\cal T};2}(\mu)
=\big\{\big(x,x,b;T_x\Si_1\!\otimes L_1|_b\big)\!\in\!{\cal Z}_{\cal T}\!
\big\}
\approx \Si\times{\cal U}_{\cal T}(\mu).$$
Thus, by Proposition~\ref{euler_prp},
${\cal C}_{{\cal Z}_{\cal T}-\tilde{\cal S}_{{\cal T};2}(\mu)}
\big(\al^{\tilde{E}\perp}\big)=
N\big(\pi^{\perp}_{\bar{\nu}}\circ\al_{{\cal Z}_{\cal T}}\big)$.
By the same argument as in (2) above,
$N\big(\pi^{\perp}_{\bar{\nu}}\circ\al_{{\cal Z}_{\cal T}}\big)
\!=\!N\big(\tilde{\al}_{\cal T}\big)$, where
\begin{gather*}
\tilde{\al}_{\cal T}\in\Ga\big({\cal Z}_{\cal T};
\hbox{Hom}(F_2;{\cal O}_2)\big),\quad
\tilde{\al}_{\cal T}(u)= \pi^{\perp}_{\bar{\nu}}\circ
\big(\big\{{\cal D}_2\otimes\pi_{x_1}^-\circ s_{x_2}\big\}\circ u\big);\\
F_2=T\Si_1^*\!\otimes\! L_1^*\!\otimes\! T\Si_2\!\otimes\! L_2
\approx T\Si_1^*\!\otimes T\Si_2,\quad
{\cal O}_2=\ga_{\tilde{E}}^*\!\otimes\!
\big({\cal H}_{\Si_1}^-\!\otimes\!\ev^*T\PP\big)^{\perp}
\approx T\Si_1^*\!\otimes\!\big({\cal H}_{\Si_1}^-\!\otimes\!\Bbb{C}^2
\big)^{\perp}.
\end{gather*}
Thus, applying Lemma~\ref{zeros_main} and again Proposition~\ref{euler_prp},
similarly to (2) we obtain
\begin{gather*}
N\big(\tilde{\al}_{\cal T}\big)
=\big\lan c_1(F_2^*)c_1({\cal O}_2)\!+\!c_2({\cal O}_2),
\big[{\cal Z}_{\cal T}\big]\big\ran-
{\cal C}_{\tilde{\al}_{\cal T}^{-1}(0)}\big(\tilde{\al}_{\cal T}^{\perp}\big)
=48\big|{\cal U}_{\cal T}(\mu)\big|-
{\cal C}_{\tilde{\al}_{\cal T}^{-1}(0)}\big(\tilde{\al}_{\cal T}^{\perp}\big);
\\
{\cal C}_{\tilde{\al}_{\cal T}^{-1}(0)}\big(\tilde{\al}_{\cal T}^{\perp}\big)
=\big\lan c_1(T\Si^*)+c_1\big(F_2^*\otimes{\cal O}_2^{\perp}\big),
\big[\tilde{\cal S}_{{\cal T};2}(\mu)\big]\big\ran
=24\big|{\cal U}_{\cal T}(\mu)\big|.
\end{gather*}
On the other hand, by an argument similar to (3) above,
$\tilde{\cal S}_{{\cal T};2}(\mu)$ is $\al^{\tilde{E}\perp}$-hollow.
We conclude~that
\begin{equation}\label{m2k1_e9}
\sum_{|\chi({\cal T})|=2}{\cal C}_{\Si^2\times{\cal U}_{\cal T}(\mu)}
\big(\al^{\tilde{E}\perp}\big)=
24\cdot 3\tau_3(\mu)=72\tau_3(\mu).
\end{equation}
(5) We finally compute the $\al^{\tilde{E}\perp}$-contribution 
to $e\big(\ga_{\tilde{E}}^*\!\otimes{\cal O}^{\perp}\big)$
from the space
$${\cal Z}_{2;2}(\mu)=
\big\{\big(x,x,b;[v\otimes\ups_1,v\otimes\ups_2]\big)\!\in\!
\Bbb{P}\tilde{E}\!: {\cal D}_1\ups_1\!+\!{\cal D}_1\ups_2\!=\!0,~
\ups_1,\ups_2\!\neq\!0\big\}\approx\Si\times{\cal S}_{2;2}(\mu).$$
Its normal bundle in $\Bbb{P}\tilde{E}$ is 
$${\cal NZ}_{2;2}=T\Si\oplus\ga_{E_2}^*\otimes\ev^*T\PP,
\quad\hbox{where}\quad
E_2=L_1\oplus L_2\lra\bar{\cal V}_2(\mu).$$
With appropriate identifications, 
${\cal D}X\!=\!X$ for all $X\!\in\!\ga_{E_2}^*\otimes\ev^*T\PP$,
where ${\cal D}\!\in\!
     \Ga\big(\Bbb{P}E_2;\ga_{E_2}^*\!\otimes\!\ev^*T\PP\big)$ 
is the section defined in the proof of Lemma~\ref{n2tacnodes}.
Then,
\begin{gather*}
\big|\al^{\tilde{E}}(x,x,b;w,X)-\tilde{\al}_{2;2}(x,b;w,X)\big|
\le C|w|^2\quad\forall(w,X)\!\in\!{\cal NZ}_{2;2;\de},
\qquad\hbox{where}\\
\tilde{\al}_{2;2}(x,b;w,X)=
X\otimes s_x+{\cal D}_{2;2}\otimes s^{(2)}_{g_x,x}(w,\cdot),
\end{gather*}
and ${\cal D}_{2;2}\!\in\!
     \Ga\big(\bar{\cal S}_{2;2}(\mu);\ga_{E_2}^*\!\otimes\!\ev^*T\PP\big)$
is the section defined in Subsection~\ref{str_global_sec}.
With the identification of small neighborhoods of $\lap$ in 
$T\Si\!\lra\!\lap$ and $\Si^2$ used above, 
the coefficients defining ${\cal D}_{2;2}$ 
are $c_1\!=\!0$ and~$c_2\!=\!1$.
By the same argument as in (2) and (4),
${\cal C}_{{\cal Z}_{2;2}}\big(\al^{\tilde{E}\perp}\big)
\!=\!N\big(\tilde{\al}_{2;2}^-)$, where
\begin{gather*}
\tilde{\al}_{2;2}^-\in\Ga\big(\Si\!\times\!\bar{\cal S}_{2;2}(\mu);
\hbox{Hom}(F_2;{\cal O}_2)\big),\quad
\tilde{\al}_{2;2}^-(w)= \pi^{\perp}_{\bar{\nu}}\circ 
     \big\{{\cal D}_{2;2}\otimes s_x^{(2)}(w,\cdot)\big\};\\
F_2=T\Si,\quad
{\cal O}_2=\ga_{\tilde{E}}^*\!\otimes\!
\big({\cal H}_{\Si}^-\!\otimes\!\ev^*T\PP\big)^{\perp}
\approx T\Si^*\otimes\ga_{E_2}^*\otimes
\big({\cal H}_{\Si}^-\!\otimes\!\ev^*T\PP\big)^{\perp}.
\end{gather*}
By Lemmas~\ref{zeros_main} and~\ref{n2tacnodes},
\begin{equation*}\begin{split}
&N\big(\tilde{\al}_{2;2}^-)=
\big\lan c_1(F_2^*)c_1({\cal O}_2)+c_2({\cal O}_2),
\big[\Si\!\times\!\bar{\cal S}_{2;2}(\mu)\big]\big\ran
-{\cal C}_{\tilde{\al}_{2;2}^{\!-~\!\!-\!1}(0)}
\big(\tilde{\al}_{2;2}^{-\perp}\big)\\
&~~=4\big\lan 120a^2\!+\!66a\big(c_1({\cal L}_1^*)\!+\!c_1({\cal L}_2^*)\big)
\!+\!13\big(c_1^2({\cal L}_1^*)\!+\!c_1^2({\cal L}_2^*)\big)
\!+\!13c_1({\cal L}_1^*)c_1({\cal L}_2^*),
\big[\bar{\cal V}_2(\mu)\big]\big\ran-
{\cal C}_{\tilde{\al}_{2;2}^{\!-~\!\!-\!1}(0)}
\big(\tilde{\al}_{2;2}^{-\perp}\big).
\end{split}\end{equation*}
By Proposition~\ref{euler_prp} and an argument similar to
the proof of Lemma~\ref{n2cusps2l},
\begin{equation*}\begin{split}
{\cal C}_{\tilde{\al}_{2;2}^{\!-~\!\!-\!1}(0)}&=
\big\lan c_1\big(F_2^*\otimes{\cal O}_2^{\perp}\big),
\big[\Si\!\times\!\partial\bar{\cal S}_{2;2}(\mu)\big]\big\ran
=28\big|\partial\bar{\cal S}_{2;2}(\mu)\big|\\
&=28\big|{\cal S}_{2;1}(\mu)\big|+ 84\tau_3(\mu)
+28\sum_{[{\cal T}^*]}\sum_{l\in[N]}
\big\lan c({\cal L}_1^*)\!+\!c({\cal L}_2^*),
\big[\bar{\cal U}_{{\cal T}^*(l)}(\mu)\big]\big\ran.
\end{split}\end{equation*}
It follows that
\begin{equation}\label{m2k1_e11}\begin{split}
{\cal C}_{{\cal Z}_{2;2}}\big(\al^{\tilde{E}\perp}\big)=
4\big\lan 120a^2+66a\big(c_1({\cal L}_1^*)\!+\!c_1({\cal L}_2^*)\big)
+7\big(c_1^2({\cal L}_1^*)\!+\!c_1^2({\cal L}_2^*)\big)
+7c_1({\cal L}_1^*)c_1({\cal L}_2^*)\qquad&\\
+6\big(c_1^2(L_1^*)\!+\!c_1^2(L_2^*)\big)+6c_1(L_1^*)c_1(L_2^*),
\big[\bar{\cal V}_2(\mu)\big]\big\ran
-28\big|{\cal S}_{2;1}(\mu)\big|-84\tau_3(\mu)&.
\end{split}\end{equation}
Combining equations~\ref{m2k1_e3}, \ref{m2k1_e5}, \ref{m2k1_e7}, 
\ref{m2k1_e9}, and~\ref{m2k1_e11}, we obtain
\begin{equation*}\begin{split}
n_2^{(1)}(\mu)=4\big\lan 24a^2
+6a\big(c_1({\cal L}_1^*)\!+\!c_1({\cal L}_1^*)\big)
+3c_1({\cal L}_1^*)c_1({\cal L}_1^*)
-\big(c_1^2({\cal L}_1^*)\!+\!c_1^2({\cal L}_1^*)\big)
,\big[\bar{\cal V}_2(\mu)\big]\big\ran&\\
+4\big|{\cal S}_{2;1}(\mu)\big|+12\tau_3(\mu)&.
\end{split}\end{equation*}
The claim then follows by using Lemma~\ref{n2_2cusps}.

\subsection{The Numbers $n^{(2)}_1(\mu)$ and $n^{(1)}_1(\mu)$}

\noindent
We now give topological formulas for the two remaining 
numbers of Corollary~\ref{CR3_str}.
It is possible to obtain the same formulas by going 
through a lengthy computation like in the proof of Lemma~\ref{m2k1}.
Instead we take slightly more geometric approaches.\\

\noindent
Suppose $\bar{\cal Z}$ and $\Si$ are topological spaces,
$L_{\cal Z},V_{\cal Z},E_{\cal Z}\!\lra\!\bar{\cal Z}$ and
$L_{\Si},V_{\Si}\!\lra\!\Si$ are vector bundles, and
$\al_{\cal Z}\!\in\!\Ga\big(\bar{\cal Z};
\hbox{Hom}(E_{\cal Z};L_{\cal Z}^*\otimes V_{\cal Z})\big)$,
$s\!\in\!\Ga\big(\Si;L_{\Si}^*\otimes V_{\Si}\big)$
and $\bar{\nu}\!\in\!\Ga\big(\Si\!\times\!\bar{\cal Z};
                                V_{\Si}\!\otimes V_{\cal Z}\big)$
are sections such that $\bar{\nu}$ does not vanish.
Then we define
\begin{gather*}
\al_{{\cal Z},\bar{\nu}}^s\in\Ga\big(\Si\!\times\!\bar{\cal Z};
\hbox{Hom}\big(E_{\cal Z};
L_{\Si}\otimes L_{\cal Z}^*\otimes (V_{\Si}\otimes V_{\cal Z})^{\perp}
\big)\big),
\quad\hbox{where}\quad
\big(V_{\Si}\otimes V_{\cal Z}\big)^{\perp}=
\big(V_{\Si}\otimes V_{\cal Z}\big)/\Bbb{C}\bar{\nu},\\
\hbox{by}\qquad
\al_{{\cal Z},\bar{\nu}}^s\big(e,w\otimes\ups\big)=
\pi_{\bar{\nu}}^{\perp}\big(\{\al_{\cal Z}(e)\}(\ups)\otimes s(w)\big)
\in  (V_{\Si}\otimes V_{\cal Z})^{\perp}.
\end{gather*}

\begin{lmm}
\label{m1k2_l1}
Suppose $L_{\cal Z},V_{\cal Z},E_{\cal Z}\!\lra\!\bar{\cal Z}$ 
are ms-bundles of rank
one, two, and $\big(2\!-\!\frac{1}{2}\dim\bar{\cal Z}\big)$,
and
$$\al_{\cal Z}\in\Ga\big(\bar{\cal Z};
\hbox{Hom}(E_{\cal Z};L_{\cal Z}^*\otimes V_{\cal Z})\big)$$
is a regular polynomial.
Let $\Si$ be a smooth compact oriented two-manifold, 
$L_{\Si},V_{\Si}\!\lra\!\Si$ smooth vector bundles of rank one and two, 
respectively, and $s\!\in\!\Ga\big(\Si;L_{\Si}^*\otimes V_{\Si}\big)$
a nonvanishing section. 
Then for an open collection of nonvanishing sections
$\bar{\nu}\!\in\!\Ga\big(\Si\!\times\!\bar{\cal Z};
                                V_{\Si}\!\otimes V_{\cal Z}\big)$\\
(1) $\al_{{\cal Z},\bar{\nu}}^s$ is a regular polynomial;\\
(2) if $V_{\cal Z}\!\approx\!\bar{\cal Z}\!\times\!\Bbb{C}^2$,
$N\big(\al_{{\cal Z},\bar{\nu}}^s\big)
=N\big(\al_{\cal Z}\big)\lan 3c_1(L_{\Si}^*)\!+\!2c_1(V_{\Si}),[\Si]\ran$.
\end{lmm}

\noindent
{\it Proof:} (1) The first claim is clear. If 
$V_{\cal Z}\!\approx\!\bar{\cal Z}\!\times\!\Bbb{C}^2$,
we can choose section 
$\bar{\nu}\!\in\!\Ga\big(\Si\!\times\!\bar{\cal Z};
                                V_{\Si}\!\otimes V_{\cal Z}\big)$
that does not intersect $V_{\Si}^+\!\otimes V_{\cal Z}$,
where $V_{\Si}^+\!=\!\hbox{Im}~\!s$.
Then 
$\pi^{\perp}_{\bar{\nu}}\big(V_{\Si}^+\!\otimes V_{\cal Z}\big)$
is a rank-two subbundle of~$\big(V_{\Si}\!\otimes V_{\cal Z}\big)^{\perp}$.
Let
$${\cal O}^+=
L_{\Si}^*\otimes L_{\cal Z}^*\otimes
\pi^{\perp}_{\bar{\nu}}\big(V_{\Si}^+\!\otimes V_{\cal Z}\big),\qquad
{\cal O}^-\equiv
\big(L_{\Si}^*\otimes L_{\cal Z}^*\otimes
(V_{\Si}\!\otimes V_{\cal Z})^{\perp}\big)\big/
{\cal O}^+.$$
We identify ${\cal O}^-$ with a complement of ${\cal O}^+$
in 
$L_{\Si}^*\otimes L_{\cal Z}^*
      \otimes(V_{\Si}\!\otimes V_{\cal Z})^{\perp}$.\\
(2) By definition, $N\big(\al_{\cal Z})$ is the number of zeros
of the affine map
$$\psi_{\al_{\cal Z},\nu}\!:
E_{\cal Z}\lra L_{\cal Z}^*\otimes V_{\cal Z},\qquad
\psi_{\al_{\cal Z},\nu}(b;\ups)=\nu_b+\al_{\cal Z}(\ups)$$
for a generic section 
$\nu\!\in\!\Ga\big(\bar{\cal Z};L_{\cal Z}^*\otimes V_{\cal Z}^*\big)$.
Via the construction preceding the lemma, $\nu$~induces
a section 
$\tilde{\nu}^+\!\in\!\Ga\big(\Si\!\times\!\bar{\cal Z}\!;{\cal O}^+\big)$.
If $\nu$ and 
$\tilde{\nu}^-\!\in\!\Ga\big(\Si\!\times\!\bar{\cal Z}\!;{\cal O}^-\big)$
are generic, $N\big(\al_{{\cal Z},\bar{\nu}}^s\big)$ is the number of
zeros of the affine map
$$\psi_{\al_{{\cal Z},\bar{\nu}}^s,\tilde{\nu}^+ +\tilde{\nu}^-}\!:
E_{\cal Z}\lra {\cal O}^+\oplus{\cal O}^-,\quad
\psi_{\al_{{\cal Z},\bar{\nu}}^s,\tilde{\nu}^+ +\tilde{\nu}^-}
\big(x,b;\ups\big)=
\tilde{\nu}^+_{(x,b)}+\tilde{\nu}^-_{(x,b)}+
\al_{{\cal Z},\bar{\nu}}^s\big(x,b;\ups\big).$$
The solution of the ${\cal O}^+$-part of this equation
is precisely $\Si\!\times\!\psi_{\al_{\cal Z},\nu}^{-1}(0)$.
Thus,
$$N\big(\al_{{\cal Z},\bar{\nu}}^s\big)=
^{\pm}\!\big|
\psi_{\al_{{\cal Z},\bar{\nu}}^s,\tilde{\nu}^+ +\tilde{\nu}^-}^{-1}(0)\big|
=\lan c_1({\cal O}^-),[\Si]\ran N\big(\al_{\cal Z}),$$
as claimed.

\begin{lmm}
\label{m1k2_l2}
Suppose $\bar{\cal Z}$ is an ms-manifold of dimension two, and
$\Si$, $L_{\cal Z}$, $V_{\cal Z}$, $E_{\cal Z}$, $L_{\Si}$, $V_{\Si}$,
$\al_{\cal Z}$ and $s$ are in Lemma~\ref{m1k2_l1}.
Then for an open collection of sections
$\bar{\nu}\!\in\!\Ga\big(\Si\!\times\!\bar{\cal Z};
                                V_{\Si}\!\otimes V_{\cal Z}\big)$
\begin{gather*}
{\cal C}_{\al_{{\cal Z},\bar{\nu}}^{s~\!\!-\!1}(0)}
\big(\al_{{\cal Z},\bar{\nu}}^{s~\!\!\perp}\big)
={\cal C}_{\al_{\cal Z}^{-1}(0)}\big(\al_{\cal Z}^{\perp}\big)
\big\lan 3c_1(L_{\Si}^*)\!+\!2c_1(V_{\Si}),[\Si]\big\ran;\\
N\big(\al_{{\cal Z},\bar{\nu}}^s\big)=
N\big(\al_{\cal Z}\big)
\big\lan 3c_1(L_{\Si}^*)\!+\!2c_1(V_{\Si}),[\Si]\big\ran
+\big\lan c_1(V_{\cal Z}),[\bar{\cal Z}]\big\ran
\big\lan c_1(L_{\Si}^*)\!+\!c_1(V_{\Si}),[\Si]\big\ran.
\end{gather*}
\end{lmm}

\noindent
{\it Proof:} (1) Since $\al_{\cal Z}$ is regular
and $\bar{\cal Z}$ is two-dimensional, 
$\al_{\cal Z}^{-1}(0)$ is a finite set of points.
Thus, we can trivialize the bundles 
$E_{\cal Z}$, $L_{\cal Z}$, and $V_{\cal Z}$ near~$\al_{\cal Z}^{-1}(0)$,
so that $\nu_z\!=\!(0,1)\!\in\! L_{\cal Z}^*\otimes V_{\cal Z}$,
where $\nu$ is as in the proof of Lemma~\ref{m1k2_l1}.
Furthermore, for each $z\!\in\!\al_{\cal Z}^{-1}(0)$, 
there exists $d_z\!\ge\!1$ such that
$$\big|\al_{\cal Z}(u)-(u^{d_z},0)\big|\le\ve(u)|u|^{d_z}
\qquad\forall u\!\in\!\Bbb{C}_{\de},
\qquad\hbox{with}~~\lim_{u\lra0}\ve(u)=0.$$
Then ${\cal C}_z\big(\al_{\cal Z}^{\perp}\big)\!=\!d_z$.\\
(2) By transversality and dimension-counting, 
we can choose $\bar{\nu}\!\in\!\Ga\big(\Si\!\times\!\bar{\cal Z};
                                V_{\Si}\!\otimes V_{\cal Z}\big)$
such~that 
$$\hbox{Im}~\bar{\nu}\cap
\big(s(L_{\Si})\otimes\al_{\cal Z}(E_{\cal Z}\!\otimes L_{\cal Z})\big)
\!=\!\eset\quad\hbox{and}\quad
\hbox{Im}~\bar{\nu}\cap
\big(s(L_{\Si})\otimes\big(V_{\cal Z}|\Si\!\times\!\al_{\cal Z}^{-1}(0)\big)
\big)\!=\!\eset.$$
Then, $\al_{{\cal Z},\bar{\nu}}^{s~\!\!-\!1}(0)\!=\!
                         \Si\!\times\!\al_{\cal Z}^{-1}(0)$.
Furthermore, on a neighborhood of 
$\Si\!\times\!\al_{\cal Z}^{-1}(0)$, we can define a splitting
$$L_{\Si}^*\otimes(V_{\Si}\otimes\Bbb{C}^2)^{\perp}=
{\cal O}^+\oplus{\cal O}^-\approx\Bbb{C}^2\oplus{\cal O}^-,$$
as in the proof of Lemma~\ref{m1k2_l1}.
Let 
$$\tilde{\nu}\in
\Ga\big(\Si\!\times\!{\cal Z};L_{\Si}^*\!\otimes L_{\cal Z}^*\!\otimes
(V_{\Si}\otimes V_{\cal Z})^{\perp}\big)$$
be a nonvanishing section such that on a neighborhood of
$\Si\!\times\!\al_{\cal Z}^{-1}(0)$,
$\tilde{\nu}=\tilde{\nu}^+ +\tilde{\nu}^-$,
with $\tilde{\nu}^+$ as in the proof of Lemma~\ref{m1k2_l1}.
Then,
$$\big({\cal O}^+\oplus{\cal O}^-)\big/\Bbb{C}\tilde{\nu}
\approx \Bbb{C}\oplus{\cal O}^-,
\quad\hbox{and}\quad
\big|\al_{{\cal Z},\bar{\nu}}^{s~\!\perp}-(u^{d_z},0)\big|
\le \ve(u)|u|^{d_z}\quad \forall u\!\in\!\Bbb{C}_{\de}.$$
Thus, by Proposition~\ref{euler_prp},
$${\cal C}_{\Si\times\{z\}}\big(\al_{{\cal Z},\bar{\nu}}^{s~\!\perp}\big)=
d_z\big\lan c_1({\cal O}^-),[\Si]\big\ran
={\cal C}_z\big(\al_{\cal Z}\big)
\big\lan 3c_1(L_{\Si}^*)\!+\!2c_1(V_{\Si}),[\Si]\big\ran.$$
The second claim follows from the first, since
by Lemma~\ref{zeros_main} and Proposition~\ref{euler_prp},
\begin{gather*}
N\big(\al_{\cal Z}\big)=
\big\lan c_1\big(L_{\cal Z}^*\!\otimes V_{\cal Z}\big)\!+\!c_1(E_{\cal Z}),
[\bar{\cal Z}]\big\ran-
{\cal C}_{\al_{\cal Z}^{-1}(0)}\big(\al_{\cal Z}\big);\\
N\big(\al_{{\cal Z},\bar{\nu}}^s\big)=
\big\lan c_2\big(L_{\Si}^*\!\otimes L_{\cal Z}^*\otimes
(V_{\Si}\!\otimes V_{\cal Z})\big)\!+\!
c_1(E^*)c_1\big(L_{\Si}^*\!\otimes L_{\cal Z}^*\otimes
(V_{\Si}\!\otimes V_{\cal Z})\big),\big[\Si\!\times\!\bar{\cal Z}\big]
\big\ran-
{\cal C}_{\al_{{\cal Z},\bar{\nu}}^{s~\!\!-\!1}(0)}
\big(\al_{{\cal Z},\bar{\nu}}^{s~\!\!\perp}\big).
\end{gather*}

\begin{crl}
\label{m1k2_c}
Suppose $\bar{\cal M}$ is an ms-manifold of dimension four,
$L_{\cal M},V_{\cal M}\!\lra\!\bar{\cal M}$ are ms-bundles
of rank one and two, respectively, and 
$\al\!\in\!\Ga\big(\bar{\cal M};\hbox{Hom}(L_{\cal M},V_{\cal Z})\big)$
is a regular polynomial.
If $\Si$ is a compact smooth oriented manifold of dimension~two,
$L_{\Si},V_{\Si}\!\lra\!\bar{\cal M}$ are ms-bundles
of rank one and two, respectively,
and $s\!\in\!\Ga\big(\Si;\hbox{Hom}(L_{\Si},V_{\Si})\big)$
is a nonvanishing section, then
\begin{equation*}\begin{split}
N\big(\al\otimes s\big)=& \big\lan 
\big(c_1(L_{\Si}^*)\!+\!c_1(V_{\Si})\big)
\big(c_1(L_{\cal M}^*)c_1(V_{\cal M})\!+\!c_1^2(V_{\cal M})\big)
-c_1(L_{\Si}^*)c_2(V_{\cal M}),
\big[\Si\!\times\!\bar{\cal M}\big]\big\ran\\
&\qquad- \big\lan c_1(L_{\Si}^*)\!+\!c_1(V_{\Si}),[\Si]\big\ran
\!\!\!\sum_{\al^{-1}(0)=\bigsqcup{\cal Z}_i}\!\!\!
\big\lan c_1(V_{\cal M}),[\bar{\cal Z}_i]\big\ran,
\end{split}\end{equation*}
where the sum is taken over all $\al$-regular subsets ${\cal Z}_i$
in a decomposition of $\al^{-1}(0)$ as in Proposition~\ref{euler_prp}.
\end{crl}

\noindent
{\it Proof:} By Lemma~\ref{zeros_main},
$$N\big(\al\otimes s\big)= \big\lan 
c_3\big(L_{\Si}^*\otimes L_{\cal M}^*\otimes
(V_{\Si}\otimes V_{\cal M})^{\perp}\big),
\big[\Si\!\times\!\bar{\cal M}\big]\big\ran
-{\cal C}_{\Si\times\al^{-1}}\big((s\otimes\al)^{\perp}\big).$$
The last term can be written as the sum of terms
as in the second equation of Lemma~\ref{m1k2_l1}.
On the other hand, by Proposition~\ref{euler_prp}, 
$$\sum_{\al^{-1}(0)=\bigsqcup{\cal Z}_i}N(\al_{{\cal Z}_i})=
\big\lan c_2\big(L_{\cal M}^*\otimes V_{\cal M}\big),
[\bar{\cal M}]\big\ran.$$
Thus, the claim follows from Lemma~\ref{m1k2_l2}.

\begin{lmm}
\label{m1k2}
If $d$ is a positive integer and $\mu$ is
a tuple of $3d\!-\!4$ points in general position in~$\PP$,
$$n_1^{(2)}(\mu)=
12\big\lan 7a^2\!+\!6ac_1({\cal L}^*),
\big[\bar{\cal S}_1(\mu)\big]\big\ran
-12\big\lan 9a^2\!+\!3a\big(c_1({\cal L}_1^*)\!+\!c_1({\cal L}_2^*)\big),
\big[\bar{\cal V}_2(\mu)\big]\big\ran.$$
\end{lmm}

\noindent
{\it Proof:} (1) By Subsection~\ref{str_thm_sec}, 
$n_1^{(2)}(\mu)\!=\!N(\al_{1;1})$, where 
\begin{gather*}
\al_{1;1}\!\in\!\Ga\big( \Si\!\times\!\bar{\cal S}_1(\mu);
\hbox{Hom}(T\Si^{\otimes2}\!\otimes L^{\otimes2},{\cal O})\big),
\quad
{\cal O}={\cal H}_{\Si}^-\!\otimes\ev^*T\PP,\\
\al_{1;1}\big(x,b; v\otimes\ups\big)
= \big({\cal D}^{(2)}\ups\big)\big(s_x^{(2)}v)
\in {\cal H}_{\Si}^-(x)\otimes T_{\ev(b)}\PP.
\end{gather*}
Thus, we can apply Corollary~\ref{m1k2_c} with 
$$\bar{\cal M}=\bar{\cal S}_1(\mu),\quad
L_{\Si}=T\Si^{\otimes2},\quad
V_{\Si}={\cal H}_{\Si}^-,\quad
L_{\cal M}=L^{\otimes2},\quad
V_{\cal M}=\ev^*T\PP,\quad
s\!=\!s^{(2)},\quad \al\!=\!{\cal D}^{(2)}.$$
The first term of Corollary~\ref{m1k2_c} gives
the intersection number on $\bar{\cal S}_1(\mu)$
in the statement of the lemma, since $ac_1(L^*)\!=\!ac_1({\cal L}^*)$.
A decomposition of the zero set of~${\cal D}^{(2)}$
is given in the proof of Lemma~\ref{n2cusps2}.
The only stratum of~$\al^{-1}(0)$ contributing to 
the second term of Corollary~\ref{m1k2_c} is~${\cal S}_{2;2}(\mu)$.
Lemma~\ref{n2tacnodes} reduces this contribution to 
the intersection number on~$\bar{\cal V}_2(\mu)$ of the~lemma.

\begin{lmm}
\label{m1k1} $n_1^{(1)}(\mu)=0$
\end{lmm}

\noindent
{\it Proof:} By Subsection~\ref{str_thm_sec}, 
$n_1^{(1)}(\mu)\!=\!N(\al_1)$, where 
\begin{gather*}
\al_1\!\in\!\Ga\big( \Si\!\times\!\bar{\cal V}_1(\mu);
\hbox{Hom}(T\Si\!\otimes L,{\cal O})\big),
\qquad
{\cal O}={\cal H}_{\Si}^{0,1}\!\otimes\ev^*T\PP,\\
\al_1\big(x,b;v\otimes\ups\big)
= \big({\cal D}\ups\big)\big(s_xv)
\in {\cal H}_{\Si}^{0,1}\otimes T_{\ev(b)}\PP.
\end{gather*}
It will be shown that there exists 
$\bar{\nu}\!\in\!\Ga\big(\PP;{\cal H}_{\Si}^{0,1}\otimes T\PP\big)$
such that the affine map
\begin{equation}\label{m1k1_e1}
\psi_{\al_1,\ev^*\bar{\nu}}\!:  
T\Si\!\otimes L\lra {\cal O},\qquad
\big(x,b;v\otimes\ups\big)\lra 
\bar{\nu}_{\ev(b)}+\al_1\big(x,b;v\otimes\ups\big),
\end{equation}
has no zeros over $\Si\!\times\!\bar{\cal V}_1(\mu)$.
The map 
$$\tilde{\al}\!:\Si\times\Bbb{P}T\PP\lra 
\Bbb{P}\big({\cal H}_{\Si}^{0,1}\otimes T\PP\big)\approx
\Bbb{P}\big(\Bbb{C}^3\otimes T\PP\big),\qquad
(x,\ell)\lra (\hbox{Im}~s_x)\otimes\ell,$$
is an embedding, since $\Si$ is not hyperelleptic.
Let ${\cal W}$ denote the image of 
$\ga_{T\PP}\big|\big(\hbox{Im}~\tilde{\al}\big)$ under the projection map
$$\ga_{T\PP}\lra T\PP,\qquad (q,\ell,v)\lra(q,v).$$
Then ${\cal W}$ is a closed subspace of ${\cal H}_{\Si}^{0,1}\otimes T\PP$ 
and ${\cal W}\!\lra\!\PP$ is a bundle of affine varieties 
of dimension three. 
Thus, by transversality and dimension-counting,
we can choose 
$\bar{\nu}\!\in\!\Ga\big(\PP;{\cal H}_{\Si}^{0,1}\otimes T\PP\big)$
such that the image $\bar{\nu}$ does not intersect~${\cal W}$.
Then the map $\psi_{\al_1,\ev^*\bar{\nu}}$ of~\e_ref{m1k1_e1}
does not vanish.

\end{document}